\documentclass[11pt,a4paper]{article}
\usepackage{setspace}
\usepackage{float}
\usepackage{multirow}
\usepackage{makecell}
\usepackage{amsthm}
\usepackage{mathtools}
\usepackage{libertine}
\usepackage[utf8]{inputenc}
\usepackage[margin=1in]{geometry}
\usepackage{amsmath,amssymb}
\allowdisplaybreaks
\usepackage{multicol}
\usepackage[shortlabels]{enumitem}
\usepackage{siunitx}
\usepackage{cancel}
\usepackage{graphicx}
\usepackage{pgfplots}
\usepackage{listings}
\usepackage{tikz}
\usepackage{hyperref}
\usepackage{xcolor}
\usepackage{epsfig,  subcaption,  graphics}
\usepackage{geometry}
\usepackage{caption}
\usepackage{apacite}
\usepackage{enumitem}
\newtheorem{theorem}{{\bf Theorem}}
\newtheorem{lemma}[theorem]{{\bf Lemma}}
\newtheorem{remark}{\bf Remark}

\newtheorem{proposition}[theorem]{{\bf Proposition}}

\pgfplotsset{width=10cm,compat=1.9}
\usepgfplotslibrary{external}
\author{
    Bihan Chatterjee, Sigrún Andradóttir, and Hayriye Ayhan
    \\[1em] 
    \small H. Milton Stewart School of Industrial and Systems Engineering, \\
    \small Georgia Institute of Technology
    \\[0.5em]
    \small \text{bchatterjee9@gatech.edu}, \text{sa@gatech.edu}, \text{hayriye.ayhan@isye.gatech.edu}
}

\date{}

\begin{document}
\title{Collaboration versus Specialization in Service Systems with Impatient Customers

}
\maketitle
\begin{abstract}
We study tandem queueing systems in which servers work more efficiently in teams than on their own and customers are impatient in that they may leave the system while waiting for service. Our goal is to determine the server assignment policy that maximizes the long-run average throughput. We show that when each server is equally skilled at all tasks, the optimal policy has all the servers working together at all times. We also provide a complete characterization of the optimal policy for Markovian systems with two stations and two servers when each server's efficiency may be task dependent. We show that the throughput is maximized under the policy which assigns one server to each station (based on their relative skill at that station) unless station 2 has no work (in which case both servers work at station 1) or the number of customers in the buffer reaches a threshold whose value we characterize (in which case both servers work at station 2).
We study how the optimal policy varies with the level of server synergy (including no synergy) and also 
compare the optimal policy for systems with different customer abandonment rates (including no abandonments). Finally, we investigate the case where the synergy among collaborating servers can be task-dependent and provide numerical results.
\end{abstract}

\section{Introduction}
In the past decade, the dynamic server assignment problem has received a lot of attention, resulting in a  significant body of research addressing the question of how flexible servers should be assigned to different  tasks in order to optimize system performance. In this paper we consider tandem queueing systems where collaborating servers work more efficiently than on their own and customers may abandon the system while waiting for service. Our objective is to determine the dynamic server assignment policy that maximizes the long-run average  throughput.

More specifically, we consider a queueing network with $N \geq 1$ stations and $M \geq 1$ flexible servers. There is an infinite amount of raw material in front of station 1, infinite room for departing customers after station $N$, a finite buffer between stations $j$ and $j+1$, for $j \in \{1, \ldots, N-1\}$, whose size is $B_j$, and the system is operating under
manufacturing blocking. At any time there can be at most one customer at each station, and each server can work on at most one customer. The servers are flexible, in that they are cross-trained and allowed to switch between stations with negligible travel time. We assume that the service rate of server $i \in \{1 , \ldots, M\}$ at station $j \in \{1, \ldots, N\}$ is $\mu_{ij}$ and, without loss of generality, we assume the mean service requirement of each customer at each station is 1. When multiple servers work on a single customer, there is synergy between them. In particular, if servers $i_1, \ldots, i_k \in \{1, \ldots, M\}$ work together on a customer at station $j \in \{1, \ldots, N\}$, the service rate is equal to $\gamma \sum_{r=1}^k \mu_{i_rj}$, where $\gamma \geq 1$ is the synergistic factor. Additionally, any customer may abandon the system while waiting to be processed by any of stations 2 through $N$.

For systems where each server is equally skilled at all tasks, we show that the optimal throughput is achieved when all  the servers work together in one team at all times. We also completely characterize the optimal server assignment policy for systems with exponential service requirements and abandonment times, $N=2$ stations, and $M=2$ servers when the service rates of each server at different stations are allowed to be heterogeneous. Further, we study some properties of the optimal policy and determine how it depends on  the magnitude $\gamma$ of the synergy  among the servers as well as on the abandonment rate $\theta$ of the customers. Additionally, we propose a policy for systems with station (or task) dependent synergy and numerically verify its optimality. We also determine how the policy depends on the synergistic factors at the two stations.

The earliest works on server assignment problems (e.g., \cite{AA2005,CLV,HR,AAD2001,AAD2007,DHM}) assumed the service rates to be additive, i.e., the service rate at a station is simply the sum over the rates of all servers working there (corresponding to $\gamma=1$ in our model). However, in practice it is often observed that when multiple servers work together on a single customer, their collaboration could be synergistic and therefore the assumption that the service rates are additive is restrictive. Optimal server assignment for systems with synergistic servers has been investigated in \cite{SHD}. The authors analyze when it is better to take advantage of synergy among servers, rather than exploiting the servers’ special skills, to achieve the best possible system throughput. In \cite{WAA}, a more generalized system is considered in which the synergy among collaborating servers can be task-dependent.

Another factor that arises in many applications, such as call centers and hospital emergency rooms, is that customers can get impatient and leave the system while waiting for service.
Queues with abandonments have been well studied in the literature (see, e.g., \cite{BF,BB,JD,BDY,FPG}). However, abandonment is a difficult phenomenon to analyze exactly and there are only a few papers that study optimal server assignment policies for queues with abandonments. For example, in \cite{DGM} an optimal server assignment policy is provided for a system with two classes of customers and a single server to serve both classes when both classes have equal service rates. The authors allow the customer abandonment rates to be class dependent and derive the optimal policy under the two most common cost/reward structures, namely to maximize rewards per service or minimize holding costs per customer per unit  time. This result was extended in \cite{SHF} to $K$ classes of customers with heterogeneous service rates. In \cite{BT}, patients waiting for treatment in a hospital emergency department were modeled as parallel queues with impatient customers.
In \cite{WMB}, the authors studied tandem queueing systems of two stations with the second station having finite buffer capacity. The authors develop a recursive algorithm to compute the throughput under various staffing strategies when the customers can get impatient and leave while waiting for service at either station. A similar system has been considered in \cite{ZC2016}, but with infinite buffer capacity at the second station. The authors consider non-collaborative 
homogeneous servers and rewards for each service completion. 
Reference \cite{ZC2023} generalizes the setting by allowing new customers to arrive at both stations and considers holding costs and abandonment costs instead of rewards. Customers arriving at station 1
 may go through both phases of service and customers arriving at station 2 may go through only one phase of service. In \cite{ZC2016} and \cite{ZC2023}, the authors provide conditions under which there exists an optimal policy where the servers need not be split between
 stations. They also specify when priority rules are optimal for the single-server model.
In \cite{AA}, a complete characterization of the dynamic assignment policy of collaborative heterogeneous servers in tandem queues in the presence of customer abandonments was considered for the first time. The authors provide the policy that maximizes the long-run average throughput. They show that unlike the optimal policy without abandonments, here it is not necessarily optimal to use the entire buffer space, suggesting that there is a trade-off between using servers where they are effective and avoiding abandonments.  To the best of our knowledge, systems with both synergy among collaborating servers and customer abandonments have not been studied yet and that is the goal of this paper. 

The outline of the paper is as follows. In Section \ref{sec2}, we show that for systems with generalist servers, the optimal throughput is achieved when all servers work together at all times. In Section \ref{sec3}, we completely characterize the optimal policy for Markovian systems with two stations and two servers.  In Section \ref{sec4}, we study how the optimal policy varies as the synergy between collaborating servers increases or when the abandonment rate of the customers changes. Section \ref{Numerical Results same gamma} provides numerical experiments to validate our theoretical findings and investigate the conditions under which the simple expedite policy serves as an effective heuristic. In Section \ref{sec5} we investigate the case when the synergy can be task dependent and provide numerical results. Section \ref{Section Conclusion} summarizes our findings. Finally, proofs of some technical results are given in the Appendix, and the proof of our main theorem, together with the details of some other proofs and results, are included in the Supplementary Materials.
\section{Systems with Generalist Servers}\label{sec2}
In this section we consider systems with generalists servers (i.e., each server is equally skilled at all tasks). Here, the service rate of each server at each station can be expressed as a product of two constants, one representing the server's speed and the other representing the difficulty of the task at the station. That is $\mu_{ij}=\mu_i \delta_j$ $\forall i \in \{1, \ldots, M\}$ and $j \in \{1, \ldots,N\}$. For each $j=1, \ldots, N$, let the service requirement  of customer $k \geq 1 $ at station $j$ be $S_{k,j}$. Assume $S_{k,j}$ are independent and identically distributed (iid) (but not necessarily exponentially distributed). And independent of the $S_{k,j}$, let $Y_{k,j}$ for $j=2,\ldots,N$ be the patience time of customer $k$ at station $j$, that is the maximum time the customer will wait for service at that station before abandoning. 
Let $\Pi$ be the set of all server assignment policies under consideration and let $D_{\pi}(t)$ denote the number of departures after service completion at station $N$ under policy $\pi \in \Pi$ by time $t \geq 0$. Define
\begin{align}\label{theroughput}
    g_{\pi}=\limsup_{t \to \infty} \frac{\mathbb E[D_{\pi}(t)]}{t}
\end{align}
as the long-run average throughput under policy $\pi$. Our goal is to solve the following optimization problem
\begin{align}\label{opt problem}
    \max_{\pi \in \Pi} g_\pi.
\end{align}

 The `expedite policy' assigns all  servers to a single team that
will follow each customer from the first to the last station and only starts work on a new customer once all work on the previous customer has been completed. The next theorem states that when the servers are generalists, the expedite policy is the optimal policy. This is not surprising because in this case there is no trade-off between server's special skills and synergy and the expedite policy  takes full advantage of the synergy among the servers while preventing any abandonment of customers. 
This result generalizes Corollary 5.1 in \cite{AA} to $N>2$ stations.
\begin{theorem}\label{thm generalist}
Assume that for all $ t \geq 0$, if there is a customer in service at station $j$ at time $t$, then the expected remaining service requirement at station $j$ of that customer is bounded by a scalar $1 \leq \bar{S} < \infty$.
 When $\mu_{ij}=\mu_i \delta_j$ for all $i=1, \ldots, M$ and $j=1, \ldots, N$ and $\gamma \geq 1 $, then the expedite policy is optimal with long-run average throughput $ g_{\pi}=\frac{  \gamma \sum_{i=1}^M\mu_i}{\sum_{j=1}^N \frac{1}{\delta_j}}$.
\end{theorem}
\begin{proof}
Let $A_{\pi}(t)$, $Q_{\pi}(t)$, and  $U_\pi(t)$ be the number of customers that have entered the system by time $t$, the number of customers in the system at time $t$, and the number of  abandonments by time $t$, respectively. Then
\begin{align*}
    A_\pi(t)=Q_\pi(t)+D_\pi(t)+U_\pi(t).
\end{align*}
From \eqref{theroughput}, since $Q_{\pi}(t) \leq \sum_{j=1}^{N-1}B_j+N$ for all policies $\pi$, we have that
\begin{align}\label{limsupeq}
    g_\pi=\limsup_{t \to \infty}\bigg[ \frac{\mathbb E[A_\pi(t)]}{t}-\frac{\mathbb E[U_{\pi}(t)]}{t}-\frac{\mathbb E[Q_{\pi}(t)]}{t}\bigg]= \limsup_{t \to \infty}\bigg[ \frac{\mathbb E[A_\pi(t)]}{t}-\frac{\mathbb E[U_{\pi}(t)]}{t} \bigg].
\end{align}

Our model is equivalent to a model where the service requirements of successive customers at station $j \in \{1, \ldots, N\}$ are iid with mean $\frac{1}{\delta_j}$ and service rates depend only on the server with $\mu_{ij}=\mu_i$ $\forall i \in \{ 1, \ldots. M\}$. Let $S_k=\sum_{j=1}^N \frac{S_{k,j}}{\delta_j}$ be the total service requirement of customer $k$ and 
let $W_\pi(t) = \sum_{k=1}^{A_\pi(t)}S_k$. 
 Let $W_{\pi,p}(t)$ be the total work $\mathit{p}$erformed by time $t$ under policy $\pi$, $W_{\pi,a}(t)$ be the total unfinished service of the $\mathit{a}$bandoned customers by time $t$, and $W_{\pi,r}(t)$ be the total $\mathit{r}$emaining service requirement at time $t$ for the customers that entered the system by time $t$ and that have not abandoned. We have
 \begin{align}\label{eq1 without lim}
   E[W_\pi(t)]=\mathbb E[W_{\pi,p}(t)]+  \mathbb E[W_{\pi,a}(t)]+\mathbb E[W_{\pi,r}(t)].  
 \end{align}
 Since for any policy $\pi$, $\mathbb E[W_{\pi,r}(t)] \leq (N+\sum_{j=1}^{N-1}B_j)\bar{S} \times \sum_{j=1}^N \frac{1}{\delta_j}$ we have
\begin{align}\label{remaining work}
    \lim_{t \to \infty }\frac{\mathbb E[W_{\pi,r}(t)]}{t}=0.
\end{align} 

Let $\pi^*$ be the expedite policy. We have $W_{\pi^*,p}(t)= \gamma \sum_{i=1}^M\mu_it$. The number of abandonments under this policy at any time $t$ is $U_{\pi^*}(t)=0$. 
This yields  $\mathbb E[W_{\pi^*,a}(t)]=0$
and hence using \eqref{limsupeq}-\eqref{remaining work}
\begin{align}\label{g pi*}
     g_{\pi^*}=\limsup_{t \to \infty}\bigg[ \frac{\mathbb E[A_{\pi^*}(t)]}{t} \bigg]
\end{align}
and
\begin{align}\label{pi* ratio}
   \lim_{t \to \infty} \frac{\mathbb E[W_{\pi^*}(t)]}{t}=  \lim_{t \to \infty}\frac{\mathbb E[W_{\pi^*,p}(t)]}{t}=\gamma \sum_{i=1}^M\mu_i.
\end{align}
For all $n \geq 0$, let $Z_n=(S_{n,1}, \ldots, S_{n,N},Y_{n,2},\ldots, Y_{n,N})$.
Since the event $\{A_{\pi}(t)=n\}$ for any policy $\pi$ is completely determined by the random vectors $Z_1, \ldots, Z_{n}$ (and is independent of $Z_{n+1},Z_{n+2}, \ldots$), $A_{\pi}(t)$ is a stopping time for the sequence of vectors $\{Z_n\}$. And $A_{\pi}(t) \leq C(t)+1$ for all $t \geq 0$, where $C(t)$ is the number of customers departing station 1 by time $t$ if all servers work at station 1 at all times and there is unlimited room for completed customers after station 1. It follows from the Elementary Renewal Theorem that $\lim_{t \to \infty} \frac{\mathbb E[C(t)]}{t}=\gamma \sum_{i=1}^M \mu_i \delta_1 <\infty$ . And since $C(t)$ is a non-decreasing process we have $\mathbb E[A_{\pi}(t)] < \infty$. Therefore, from  Wald's identity,
\begin{align}\label{eq2}
    \mathbb E[W_\pi(t)]=\mathbb E\bigg[\sum_{k=1}^{A_{\pi}(t)}S_k\bigg]=\mathbb E[A_\pi(t)] \times \sum_{j=1}^N \frac{1}{\delta_j}.
\end{align}
In particular we have 
\begin{align*}
    \mathbb E[W_{\pi^*}(t)]=\mathbb E\bigg[\sum_{k=1}^{A_{\pi^*}(t)}S_k\bigg]=\mathbb E[A_{\pi^*} (t)] \times \sum_{j=1}^N \frac{1}{\delta_j}.
\end{align*}
That is, \eqref{g pi*}-\eqref{pi* ratio} imply that 
\begin{align*}
    \gamma \sum_{i=1}^M\mu_i=  \lim_{t \to \infty} \frac{\mathbb E[W_{\pi^*}(t)]}{t}=\lim_{t \to \infty} \frac{\mathbb E[A_{\pi^*} (t)]}{t} \times \sum_{j=1}^N \frac{1}{\delta_j}=g_{\pi^*}\times \sum_{j=1}^N \frac{1}{\delta_j}.
\end{align*}
Therefore,
\begin{align*}
    g_{\pi^*}=\frac{  \gamma \sum_{i=1}^M\mu_i}{\sum_{j=1}^N \frac{1}{\delta_j}}.
\end{align*}

For any other policy $\pi$ we have 
\begin{align}\label{pi work performed}
    W_{\pi,p}(t)\leq \gamma \sum_{i=1}^M\mu_it.
\end{align}
 Next we show $\mathbb E[W_{\pi,a}(t)] \leq \mathbb E[U_{\pi}(t)] \times \sum_{j=1}^N\frac{1}{\delta_j}$. 
We have 
\begin{align*}
    W_{\pi,a}(t)=\sum_{i=1}^{U_{\pi}(t)}\bigg[ \sum_{j=K_i}^N \frac{S_{R_i,j}}{\delta_j}\bigg]
\end{align*}
where $R_i$ is the $i^{th}$ abandoning customer who abandons while waiting for service at station $K_i$.
We have
\begin{align}\label{pi unfinished work}
     \mathbb E[W_{\pi,a}(t)]&=\sum_{n=1}^{\infty}\sum_{k_1,\ldots,k_n=2}^N\mathbb E \bigg[\sum_{i=1}^{n}\bigg[ \sum_{j=k_i}^N \frac{S_{R_i,j}}{\delta_j}\bigg]\bigg| U_{\pi}(t)=n, K_1=k_1,\ldots,K_n=k_n\bigg] \times \nonumber \\
     &\hspace{7cm}P(U_{\pi}(t)=n,K_1=k_1,\ldots,K_n=k_n) \nonumber\\
     &\stackrel{(*)}{=}\sum_{n=1}^{\infty}\sum_{k_1,\ldots,k_n=2}^N\mathbb E \bigg[\sum_{i=1}^{n}\bigg[ \sum_{j=k_i}^N \frac{S_{1,j}}{\delta_j}\bigg]\bigg]P(U_{\pi}(t)=n,K_1=k_1,\ldots,K_n=k_n) \nonumber\\
     &\leq \sum_{n=1}^{\infty}\sum_{k_1,\ldots,k_n=2}^Nn\bigg(\sum_{j=1}^N\frac{1}{\delta_j}\bigg)P(U_{\pi}(t)=n,K_1=k_1,\ldots,K_n=k_n) \nonumber\\
&=\bigg(\sum_{j=1}^N\frac{1}{\delta_j}\bigg)\sum_{n=1}^{\infty}nP(U_{\pi}(t)=n) \nonumber\\
 &=E[U_{\pi}(t)] \times \sum_{j=1}^N\frac{1}{\delta_j},
\end{align}
where $(*)$ follows since $S_{R_i,j}$ for $1 \leq i \leq U_{\pi}(t)$ and  $K_i \leq j \leq N$ are service requirements that are not executed because the corresponding customers have already abandoned.

Now from \eqref{eq1 without lim}, \eqref{eq2}, \eqref{pi work performed}, and \eqref{pi unfinished work}
\begin{align*}
    \frac{\mathbb E[A_\pi(t)]}{t} \times \sum_{j=1}^N \frac{1}{\delta_j} & \leq \gamma \sum_{i=1}^M \mu_i+\frac{\mathbb E[U_\pi(t)]}{t}\times \sum_{j=1}^N\frac{1}{\delta_j}+\frac{\mathbb E[W_{\pi,r}(t)]}{t}.
\end{align*}
It follows from \eqref{remaining work} that
\begin{align*}
     \limsup_{t \to \infty}\left[\frac{\mathbb E[A_\pi(t)]}{t}-\frac{\mathbb E[U_\pi(t)]}{t} \right]\times \sum_{j=1}^N\frac{1}{\delta_j} &\leq \gamma \sum_{i=1}^M \mu_i,
\end{align*}
and \eqref{limsupeq} yields
\begin{align*}
     g_\pi \times \sum_{j=1}^N\frac{1}{\delta_j}&\leq \gamma \sum_{i=1}^M \mu_i.
\end{align*}
Therefore $\pi^*$ is the optimal policy.
\end{proof}
\section{Optimal policy for 2 Stations and 2 Servers}\label{sec3}
In this section we completely characterize the optimal server assignment policy with $N=2$ stations, $M=2$ servers, and independent exponential service requirements at the stations. The customers waiting in line to be served by station 2 (after being served at station 1) may abandon the system after an exponential amount of time with parameter $\theta>0$. We set $B_1=B$ for notational convenience. We assume, without loss of generality, that $\mu_{11}\mu_{22} \geq \mu_{12} \mu_{21}$. Let $\Sigma_j=\mu_{1j}+\mu_{2j}$ and $\Sigma'_j=\gamma \Sigma_j$ for $j=1,2$. We let 
\begin{align*}
    \lambda=\gamma-1.
\end{align*} 
When both servers work in station $j$ their combined service rate is $\Sigma'_j$. 
 As in Section \ref{sec2}, our goal is to solve the optimization problem \eqref{opt problem}.
In order to avoid trivial cases, we assume $\Sigma_j >0$  for  $j=1,2$. We also assume that if $\lambda=0$, $\sum_{j=1}^2 \mu_{ij} >0$ for $i=1,2$. This is because $\lambda=0$ and $\sum_{j=1}^2 \mu_{ij} =0$ would imply that the server cannot process customers on their own and there are no synergistic effects for collaboration. Hence  we
can simply exclude the server from the system. It can be easily shown (take $\mu_i=0, \mu_j=\mu_{j1}$, $\delta_1=1$, and $\delta_2=\mu_{j2}/\mu_{j1}$ in Theorem \ref{thm generalist}) that the
expedite policy (that assigns all servers to the same job until completion) would be optimal in this case.



 Let $X_{\pi}(t)$ be the number of customers that have been processed at station 1 and are waiting to be processed or being processed at station 2 at time $t$ under policy $\pi \in \Pi$, where we have $X_{\pi}(t) \in \{0,1, \ldots, B+2\} \eqqcolon S$. From now on, we let $\Pi$ be the set of all stationary deterministic Markovian policies corresponding to the state space $S$. We define the action space $A \coloneqq \{a_{\sigma_1  \sigma_2}: \sigma_i \in \{1,2\}, \forall i=1,2\}$, where $\sigma_i=j \in \{1,2\}$ means that server $i$ is assigned to station $j$. For any given $\pi$, $\{X_{\pi}(t)\}$ is a birth-death process. Since $B < \infty$ and the sum of all transition rates out of each state $s \in S$ is bounded by $\sum_{j=1}^2 \Sigma'_j +(B+2)\theta< \infty$, $\{X_{\pi}(t)\}$ is uniformizable with a finite uniformization constant, say $q$, satisfying $\sum_{j=1}^2 \Sigma'_j +(B+2)\theta \leq q < \infty$. Thus, the continuous-time optimization problem \eqref{opt problem} can be translated into a
discrete-time Markov decision problem by uniformization (i.e., modeling the state transitions
at the event times of a Poisson process with rate $q$; see for example \cite{AAD2001}). Note that the chosen value of $q$ has no impact on
our results.

We next provide for all $s \in S$ and $a \in A$, the immediate rewards, $r(s,a)$, and the transition probability $p(j|s,a)$ of going to state $j$ from state $s$ in one step when action $a$ is chosen in state $s$.
We have for all $a \in A$
\begin{align*}
    r(0,a)=0
\end{align*}
and for all $s=1, \ldots, B+2$,
\begin{equation} \label{rewards general}
 \left.\begin{array}{r@{\;}l}
     r(s, a_{11})=&0,\\
    r(s,a_{21})=&\mu_{12},\\
    r(s,a_{12})=&\mu_{22}, \\
    r(s,a_{22})=&\Sigma'_2.
\end{array} \right\}
\end{equation}
Furthermore,
\begin{align}
    p(j|s,a_{12})=\begin{cases}
        \frac{\mu_{11}}{q} \quad & \text{for $s=0, \ldots, B+1$, $j=s+1$},\\
        \frac{\mu_{22}+(s-1)\theta}{q} & \text{for $s=1, \ldots, B+2$, $j=s-1$},\\
        1-\frac{\mu_{11}+\mu_{22}+(s-1)\theta}{q} & \text{for $s=1, \ldots, B+1$, $j=s$},\\
        1-\frac{\mu_{22}+(s-1)\theta}{q} & \text{for $s=B+2$, $j=s$},\\
        1-\frac{\mu_{11}}{q} & \text{for $s=0$, $j=s$,}\\
        0 &\text{otherwise;}
    \end{cases}
\end{align}
\begin{align}
    p(j|s,a_{21})=\begin{cases}
        \frac{\mu_{21}}{q} \quad & \text{for $s=0, \ldots, B+1$, $j=s+1$},\\
        \frac{\mu_{12}+(s-1)\theta}{q} & \text{for $s=1, \ldots, B+2$, $j=s-1$},\\
        1-\frac{\mu_{21}+\mu_{12}+(s-1)\theta}{q} & \text{for $s=1, \ldots, B+1$, $j=s$},\\
        1-\frac{\mu_{12}+(s-1)\theta}{q} & \text{for $s=B+2$, $j=s$},\\
        1-\frac{\mu_{21}}{q} & \text{for $s=0$, $j=s$,}\\
        0 &\text{otherwise;}
    \end{cases}
\end{align}
\begin{align}
    p(j|s,a_{11})=\begin{cases}
        \frac{\Sigma'_1}{q} \quad & \text{for $s=0, \ldots, B+1$, $j=s+1$,}\\
        \frac{s\theta}{q} & \text{for $s=1, \ldots, B+2$, $j=s-1$},\\
        1-\frac{\Sigma'_1+s\theta}{q} & \text{for $s=1, \ldots, B+1$, $j=s$},\\
        1-\frac{s\theta}{q} & \text{for $s=B+2$, $j=s$,}\\
        1-\frac{\Sigma'_1}{q} & \text{for $s=0$, $j=s$},\\
        0 &\text{otherwise;}
    \end{cases}
\end{align}
and 
\begin{align}
    p(j|s,a_{22})=\begin{cases}
        \frac{\Sigma'_2+(s-1)\theta}{q} & \text{for $s=1, \ldots, B+2$, $j=s-1$,}\\
        1-\frac{\Sigma'_2+(s-1)\theta}{q} & \text{for $s=1, \ldots, B+2$, $j=s$,}\\
        1 & \text{for $s=0$, $j=s$,}\\
        0 &\text{otherwise.}
    \end{cases}
\end{align}

  We next characterize the dynamic server assignment policy that maximizes the long-run average throughput. For $n=1, \ldots, B+2$, we define the decision rule $d_n$ as follows
\begin{align}\label{d_n}
   d_n(s)= \begin{cases}
    a_{11} \ \text{for } s=0,\\
    a_{12} \ \text{for } s=1, \ldots, n-1, \\
    a_{22} \ \text{for } s=n, \ldots, B+2.
    \end{cases}
\end{align}
The next proposition gives an expression for the long-run average throughput $g_n$ under policy $d_n$. 
\begin{proposition}\label{prop gain}
We have for $n=1,\ldots, B+2$
    \begin{align}\label{g_n}
        g_n=\frac{\beta(n)}{\Delta(n)},
    \end{align}
    where 
\begin{align}\label{beta(n)}
    \beta(n)=[\Sigma'_2+(n-1)\theta]\Sigma'_1\mu_{22}\alpha(n)+\Sigma'_1 \Sigma'_2\mu_{11}^{n-1} 
\end{align}
and
\begin{align}\label{Delta(n)}
    \Delta(n)=[\Sigma'_2+(n-1)\theta][f(1,n)+\Sigma'_1\alpha(n)]+\Sigma'_1\mu_{11}^{n-1}
\end{align}
with 
\begin{align*}
     \alpha(n)&=\sum_{k=2}^{n} \mu_{11}^{k-2}f(k,n) \mbox{ and }\\
     f(k,n)&=\prod_{j=k}^{n-1}[\mu_{22}+(j-1)\theta] \text{ for } k=1,2, \ldots,n,
\end{align*}
and with the convention that summation over an empty set is 0 and the multiplication over
an empty set is 1.
\end{proposition}
\begin{proof}
For the policy $d_n$, states $n+1, \ldots , B+2$ are transient. For the recurrent states we have the corresponding  birth  rates given as
\begin{align*}
    \kappa_i=\begin{cases}
        \Sigma_1' \quad &\text{ for } i=0,\\
        \mu_{11} &\text{ for } i=1,2, \ldots, n-1,\\
        0 &\text{ for } i=n,
    \end{cases}
\end{align*}
and death rates given as 
\begin{align*}
    \nu_i=\begin{cases}
        \mu_{22}+(i-1)\theta \quad &\text{ for } i=1,2, \ldots,n-1,\\
        \Sigma_2'+(i-1)\theta &\text{ for } i=n.
    \end{cases}
\end{align*}

We know that for birth-death processes with finite state space $\{0,1, \ldots, n\}$, the stationary probabilities are given by 
\begin{align}\label{p_i}
    p_i=\frac{\prod_{j=0}^{i-1}\kappa_j}{\left(\prod_{j=1}^i\nu_i\right)\left(\sum_{k=0}^n\frac{\prod_{j=0}^{k-1}\kappa_j}{\prod_{j=1}^k\nu_j }\right)} \quad \quad 0 \leq i \leq n.
\end{align}
Therefore the long-run average throughput under policy $d_n$ is given by 
\begin{align*}
    g_n
    =&\sum_{i=1}^{n-1}\mu_{22}p_i+\Sigma'_2p_n.
\end{align*}
Using \eqref{p_i} and some algebra we get
\begin{align*}
    g_n
    =&\frac{[\Sigma'_2+(n-1)\theta]\Sigma'_1 \mu_{22}\bigg[\sum_{l=1}^{n-1}\mu_{11}^{l-1}\left(\prod_{j=l+1}^{n-1}[\mu_{22}+(j-1)\theta]\right)\bigg]+\Sigma'_1 \Sigma'_2\mu_{11}^{n-1}
    }{[\Sigma'_2+(n-1)\theta]\bigg[\prod_{j=1}^{n-1}[\mu_{22}+(j-1)\theta]+\Sigma'_1 \bigg(\sum_{l=1}^{n-1}\mu_{11}^{l-1} \prod_{j=l+1}^{n-1}[\mu_{22}+(j-1)\theta]\bigg)\bigg]+\Sigma'_1 \mu_{11}^{n-1}}\\
    =&\frac{\beta(n)}{\Delta(n)},
\end{align*}
where the last equality follows from the change of variables $k=l+1$.   
\end{proof}
The  next proposition characterizes the difference $g_n-g_{n-1}$. The sign of this quantity will be essential in proving our results. Note that $g_n-g_{n-1}$ has the same sign as $\tau'(n)$.
\begin{proposition}\label{tau' proposition}
We have
\begin{align*}
    g_1=\frac{\tau'(1)}{\Sigma'_1+\Sigma'_2}
\end{align*}
and for $n=2, \ldots, B+2$,
\begin{align*}
    g_n-g_{n-1}=\frac{\Sigma'_1 \mu_{11}^{n-2} \tau'(n)}{\Delta(n) \Delta(n-1)},
\end{align*}
where
\begin{align*}
    \tau'(1)=\Sigma'_1 \Sigma'_2
\end{align*}
and for $n=2, \ldots, B+2$,
\begin{align}\label{tau'}
    \tau'(n)=&[\Sigma'_2+(n-1)\theta][\Sigma'_2+(n-2)\theta]\mu_{22}f(1,n-1) -[\Sigma'_2+(n-1)\theta][\Sigma'_1\alpha(n)(\Sigma'_2-\mu_{22})+\Sigma'_2f(1,n)] \nonumber \\
     &  +[\Sigma'_2+(n-2)\theta]\mu_{11}[\Sigma'_2f(1,n-1)+\Sigma'_1\alpha(n-1)(\Sigma'_2-\mu_{22})].
\end{align} 
\end{proposition}
\begin{proof}
From Proposition \ref{prop gain} we immediately have
\begin{align*}
    g_1=\frac{\Sigma'_1 \Sigma'_2}{\Sigma'_1+\Sigma'_2}=\frac{\tau'(1)}{\Sigma'_1+\Sigma'_2}.
\end{align*}
Next we consider $n=2, \ldots, B+2$. 
We will use the following equalities in deriving the expression for $\tau'(n)$:
\begin{align}\label{f recursion}
    f(k,n)=f(k,n-1)[\mu_{22}+(n-2)\theta], \quad k=1, \ldots,n-1,
\end{align}
and
\begin{align}\label{eq5}
     \alpha(n)
     &=\sum_{k=2}^{n-1} \mu_{11}^{k-2}f(k,n)+\mu_{11}^{n-2}f(n,n) \nonumber \\
     &=\sum_{k=2}^{n-1} \mu_{11}^{k-2}f(k,n-1)[\mu_{22}+(n-2)\theta]+\mu_{11}^{n-2}\nonumber\\
     &=[\mu_{22}+(n-2)\theta]\alpha(n-1)+\mu_{11}^{n-2}.
\end{align}
Using these and some algebra, we have
\begin{align*}
     &\beta(n)\Delta(n-1)-\beta(n-1)\Delta(n) \nonumber \\
      =&[\Sigma'_2+(n-1)\theta]\Sigma'_1\mu_{22}\alpha(n)[\Sigma'_2+(n-2)\theta][f(1,n-1)+\Sigma_1'\alpha(n-1)] \nonumber\\
    &-[\Sigma'_2+(n-2)\theta]\Sigma'_1\mu_{22}\alpha(n-1)[\Sigma'_2+(n-1)\theta][f(1,n)+\Sigma_1'\alpha(n)] \nonumber\\
     &+[\Sigma'_2+(n-1)\theta]\Sigma'_1 \mu_{22}\alpha(n)\Sigma'_1\mu_{11}^{n-2} 
      -\Sigma'_1\Sigma'_2\mu_{11}^{n-2}[\Sigma'_2+(n-1)\theta][f(1,n)+\Sigma'_1\alpha(n)] \nonumber\\
      &-[\Sigma'_2+(n-2)\theta]\Sigma'_1 \mu_{22} \alpha(n-1)\Sigma'_1 \mu_{11}^{n-1} 
    +\Sigma'_1\Sigma'_2 \mu_{11}^{n-1}[\Sigma'_2+(n-2)\theta][f(1,n-1)+\Sigma'_1\alpha(n-1)]\\
    =&\Sigma'_1[\Sigma'_2+(n-1)\theta][\Sigma'_2+(n-2)\theta]\mu_{22}f(1,n-1)\mu_{11}^{n-2}\\
    &+\Sigma'_1 \mu_{11}^{n-2}[\Sigma'_2+(n-1)\theta]\bigg[-\Sigma'_1\alpha(n)[\Sigma'_2-\mu_{22}]-\Sigma'_2f(1,n) \bigg]\\
    &+\Sigma'_1\mu_{11}^{n-2}[\Sigma'_2+(n-2)\theta]\mu_{11}\bigg[\Sigma'_2f(1,n-1)+\Sigma'_1\alpha(n-1)(\Sigma'_{2}-\mu_{22})\bigg]\\
     =&\Sigma'_1\mu_{11}^{n-2}\bigg[[\Sigma'_2+(n-1)\theta][\Sigma'_2+(n-2)\theta]\mu_{22}f(1,n-1)-[\Sigma'_2+(n-1)\theta][\Sigma'_1\alpha(n)(\Sigma'_2-\mu_{22})+\Sigma'_2f(1,n)]\\
     & \quad \quad \quad +[\Sigma'_2+(n-2)\theta]\mu_{11}[\Sigma'_2f(1,n-1)+\Sigma'_1\alpha(n-1)(\Sigma'_2-\mu_{22})]\bigg ]\\
     =&\Sigma'_1\mu_{11}^{n-2} \tau'(n).
\end{align*}
Therefore, for $n=2, \ldots, B+2$
\begin{align}\label{g difference}
    g_n-g_{n-1}=\frac{ \beta(n)\Delta(n-1)-\beta(n-1)\Delta(n)}{\Delta(n)\Delta(n-1)}=\frac{\Sigma'_1\mu_{11}^{n-2} \tau'(n)}{\Delta(n)\Delta(n-1)}.
\end{align}
\end{proof}
The next proposition provides a useful property of the $\tau'(\cdot)$ function.
\begin{proposition}\label{tau property}
    For $n \in \{ 2, \ldots, B+2 \}$ if $\tau'(n) \geq 0$, then $\tau'(k) \geq 0$ for all $1 \leq k \leq n-1$. Similarly, if $\tau'(n) \leq 0$, then $\tau'(k) \leq 0 $ for all $n+1 \leq k \leq B+2$.
\end{proposition}
\begin{proof}

For $n \geq 2$,
\begin{align}\label{tau'(n+1)}
    \tau'(n+1)=&[\Sigma'_2+n\theta][\Sigma'_2+(n-1)\theta]\mu_{22}f(1,n)-[\Sigma'_2+n\theta][\Sigma'_1\alpha(n+1)(\Sigma'_2-\mu_{22})+\Sigma'_2f(1,n+1)] \nonumber\\
     & \quad \quad \quad +[\Sigma'_2+(n-1)\theta]\mu_{11}[\Sigma'_2f(1,n)+\Sigma'_1\alpha(n)(\Sigma'_2-\mu_{22})].
\end{align}
Using some algebra (details in Supplementary Material \ref{SMProp4}), we have for $n\geq 2$  
\begin{align}\label{tau(n+1)}
    \tau'(n+1)=[\mu_{22}+(n-1)\theta]\tau'(n)-\epsilon'(n),
\end{align}
where
\begin{align}\label{epsilon'(n)}
    \epsilon'(n)=&\theta [\Sigma'_2-\mu_{22}][\Sigma'_2-\mu_{22}+\theta]\mu_{22}f(1,n-1) \nonumber \\
    &+\theta\Sigma'_1\mu_{11}^{n-1}[\Sigma'_2-\mu_{22}]+\theta[\mu_{22}+(n-1)\theta][\mu_{22}+(n-2)\theta]\Sigma'_1 [\Sigma'_2-\mu_{22}]\alpha(n-1)] \nonumber\\
    &+\theta[\mu_{22}+(n-1)\theta][\mu_{22}+(n-2)\theta][\Sigma'_2-\mu_{22}]f(1,n-1) \nonumber\\
    &    +\Sigma'_1[\Sigma'_2-\mu_{22}]\mu_{11}^{n-2}\theta[\mu_{22}+(n-1)\theta]\nonumber \\
    &+\theta [\Sigma'_2-\mu_{22}]\mu_{11}[\Sigma'_2f(1,n-1)+\Sigma'_1 (\Sigma'_2-\mu_{22})\alpha(n-1)]\\
    \geq& 0 \nonumber.
\end{align}

Also note that $\tau'(1)>0$. Hence, for $n \in \{ 2, \ldots, B+2 \}$ if $\tau'(n) \geq 0$, then $\tau'(k) \geq 0$ for all $1 \leq k \leq n-1$. And if $\tau'(n) \leq 0$, then $\tau'(k) \leq 0 $ for all $n+1 \leq k \leq B+2$.
\end{proof}
Define 
\begin{align}\label{Ndef}
    N=\max{\{n \in \{1, \ldots,B+2\}: \tau'(n) \geq 0 \}}.
\end{align}
This implies $g_N \geq g_{N-1} \geq \ldots \geq g_1$ and $g_N > g_{N+1}\ldots \geq g_{B+2}$.

From \eqref{rewards general}-\eqref{d_n}, the reward vector $r_{d_N}$ and and the transition probability matrix $P_{d_N}$ corresponding to the decision rule $d_{N}$ are given by 
\begin{align}\label{r_{d_N}}
    r_{d_N}(s)=\begin{cases}
        0 \quad &\text{ for } s=0,\\
        \mu_{22} &\text{for }s=1.\ldots ,N-1,\\
        \Sigma_2'&\text{for } s=N, \ldots, B+2;
    \end{cases}
\end{align}
\begin{align}\label{P_{d_N}}
    P_{d_N}(s,s')=\begin{cases}
        \frac{\Sigma'_1}{q} \quad &\text{for } s=0,s'=1,\\
        \frac{\mu_{11}}{q} &\text{for } s=1,\ldots,N-1,s'=s+1,\\
        \frac{\mu_{22}+(s-1)\theta}{q} &\text{for } s=1, \ldots, N-1, s'=s-1,\\
        \frac{\Sigma'_2+(s-1)\theta}{q} &\text{for } s=N,\ldots,B+2, s'=s-1,\\
        1-\frac{\Sigma'_1}{q} &\text{for } s=0,s'=0,\\
        1-\frac{\mu_{11}+\mu_{22}+(s-1)\theta}{q} &\text{for } s=1, \ldots,N-1,s'=s,\\
        1-\frac{\Sigma'_2+(s-1)\theta}{q} &\text{for } s=N, \ldots,B+2,s'=s,\\
        0 & \text{otherwise}.
    \end{cases}
\end{align}
Due to the abandonments, the Markov decision problem is unichain and therefore we can use the unichain policy iteration algorithm. In order to use the algorithm we let $h_N$  be such that 
\begin{align}\label{eq7}
    r_{d_N}-g_Ne+(P_{d_N}-I)h_N=0
\end{align}
with $h_N(0)=0$, where $e$ is a column vector of ones and $I$ is the identity matrix. The next two lemmas characterize this vector $h_N$. Lemma \ref{lemma h_N1} gives an expression of $h_N(1)$ and the differences $h_N(s)-h_N(s-1)$ in terms of the gain $g_N$. It will also be  helpful to have a closed-form expression of the difference, as given in Lemma \ref{lemma h_N2}. Note that Lemma \ref{lemma h_N2} implies that for all $s=1, \ldots, B+2$, $h_N(s)-h_N(s-1) \geq 0$. The proofs are included in the Appendix.
\begin{lemma}\label{lemma h_N1}
We have 
\begin{itemize}
    \item \begin{align}\label{lemma h_N1(i)}
        h_N(1)=\frac{qg_N}{\Sigma'_1}
    \end{align}
    \item For $s=1, \ldots, N$,
    \begin{align}\label{lemma h_N1(ii)}
        h_N(s)-h_N(s-1)=q \bigg[\sum_{i=0}^{s-2} \frac{g_N-\mu_{22}}{\mu_{11}^{i+1}}f(s-i,s)+\frac{g_N}{\Sigma'_1\mu_{11}^{s-1}}f(1,s) \bigg].
    \end{align}
    \item For $s=N,\ldots,B+2$,
    \begin{align}\label{lemma h_N1(iii)}
          h_N(s)-h_N(s-1)=\frac{q[\Sigma'_2-g_N]}{\Sigma'_2+(s-1)\theta}.
    \end{align}
\end{itemize}
\end{lemma}
\begin{lemma}\label{lemma h_N2}
For $s=1, \ldots, N$,
\begin{align}\label{lemma h_N2(i)}
    &h_N(s)-h_N(s-1) \nonumber\\
    =&\frac{q}{\Delta(N)} \bigg[ \Sigma'_1\mu_{11}^{N-s}[\Sigma'_2-\mu_{22}]\alpha(s)+[\Sigma'_2+(N-1)\theta]\mu_{22}\sum_{k=s+1}^{N} \mu_{11}^{k-s-1}f(k,N)f(1,s)+\Sigma'_2\mu_{11}^{N-s}f(1,s)\bigg ].
\end{align}
For $s=N, \ldots,B+2$,
\begin{align}\label{lemma h_N2(ii)}
     h_N(s)-h_N(s-1)=\frac{q[\Sigma'_2+(N-1)\theta]\bigg(\Sigma'_2f(1,N)+\Sigma'_1[\Sigma'_2-\mu_{22}]\alpha(N)\bigg)}{\Delta(N)[\Sigma'_2+(s-1)\theta]}.
\end{align}
\end{lemma}

Another result that will be useful is given in the following lemma. The proof is included in the Appendix.

\begin{lemma}\label{Lemma alpha}
    For $ s=1,2, \ldots, N-1$, we have 
    \begin{align}\label{alphaineq}
     \alpha(N) \geq \mu_{11}^{N-s} \alpha(s).
\end{align}
\end{lemma}
We next state and prove the main result of this section which completely characterizes the optimal server assignment policy. While Proposition \ref{tau property} implies that $d_N$ (as defined in \eqref{Ndef}) maximizes throughput within the restricted class of policies $\{d_n\}$,  Theorem \ref{thm2S2S} establishes that it is, in fact,  optimal across all policies.
 \begin{theorem}\label{thm2S2S}
     For $N$ as defined in \eqref{Ndef}, $d_N$ is the optimal policy.
 \end{theorem}
 The complete proof of the Theorem is included in the Supplementary Material \ref{SMTheorem1}.  We use policy iteration for unichain models (see \cite{Puterman}). We set the initial decision rule $\nu_0$ equal to $d_N$ and then show that the next decision rule computed by the policy iteration algorithm is equal to $d_N$. That is we show that 
\begin{align*}
    \nu_1(s) \in \arg \max_{a \in A} \bigg\{r(s,a)+\sum_{j \in S} p(j|s,a)h_N(j)\bigg\} 
\end{align*}
is equal to $d_N(s)$, $\forall s \in S$. We show this by showing the difference 
\begin{align*}
    \Gamma(s,a)=r(s,d_N(s))+\sum_{j \in S} p(j|s,d_N(s))h_N(j)-r(s,a)-\sum_{j \in S} p(j|s,a)h_N(j)
\end{align*}
is non negative for all $s \in S$ and $a \in A$ using Lemmas \ref{lemma h_N1} and \ref{lemma h_N2}. The result then follows from Theorem 8.6.6 in \cite{Puterman}.
\section{Properties of the Optimal Policy}\label{sec4}
In this section we provide some conditions under which the expedite policy is the optimal policy. We also determine how the optimal buffer size $N$ changes with the synergistic factor $\gamma$ and the abandonment rate $\theta$.

The next proposition gives a lower bound on $\gamma$ that guarantees that the expedite policy is optimal. Note that from our assumption that $\mu_{11}\mu_{22}-\mu_{21}\mu_{12} \geq 0$, we have $\gamma \geq 1$. The result states that when the synergistic factor reaches a threshold, then irrespective of the abandonment rate, the servers should always work in teams to take advantage of the synergy.
\begin{proposition}\label{lambda bound}
\begin{enumerate}
    \item[(i)]  If 
    \begin{align*}
        \gamma> \frac{\bigg(\Sigma_2+\frac{\mu_{11}\mu_{22}-\mu_{21}\mu_{12}}{\Sigma_1}-\theta \bigg) + \sqrt{\bigg(\Sigma_2+\frac{\mu_{11}\mu_{22}-\mu_{21}\mu_{12}}{\Sigma_1}-\theta\bigg)^2+4\theta \mu_{22}}}{2 \Sigma_2},
    \end{align*}
    we have $N=1$, i.e., the expedite policy is the optimal policy.
    \item[(ii)]   If  $\gamma \geq 1+ \frac{\mu_{11}\mu_{22}-\mu_{21}\mu_{12}}{\Sigma_1\Sigma_2} $,  we have $N=1$, irrespective of the value of $\theta$.
\end{enumerate}
\end{proposition}
\begin{proof}
\begin{enumerate}
    \item[(i)] We have 
    \begin{align}
    \tau'(2) 
    &=[\Sigma'_2+\theta]\Sigma'_2\mu_{22}-[\Sigma'_2+\theta][\Sigma'_1(\Sigma'_2-\mu_{22})+\Sigma'_2\mu_{22}]+\mu_{11}(\Sigma'_2)^2 \nonumber \\
     &=-[\Sigma'_2+\theta][\Sigma'_1(\Sigma'_2-\mu_{22})]+\mu_{11}(\Sigma'_2)^2 \nonumber\\
   &= -\Sigma'_2 \Sigma'_1(\Sigma'_2-\mu_{22})+(\Sigma'_2)^2\mu_{11} -\theta \Sigma'_1(\Sigma'_2-\mu_{22}) \label{tau'(2.2)}\\
   &=-\Sigma'_2 \Sigma'_1(\lambda \Sigma_2+\mu_{12})+(\Sigma'_2)^2\mu_{11} -\theta \Sigma'_1(\Sigma'_2-\mu_{22}) \nonumber\\
   &=\Sigma'_2(1+\lambda)[-\lambda\Sigma_1\Sigma_2-\mu_{12}\Sigma_1+\Sigma_2\mu_{11}] -\theta \Sigma'_1(\Sigma'_2-\mu_{22}) \nonumber\\
   &=\Sigma'_2(1+\lambda)[-\lambda\Sigma_1\Sigma_2+\mu_{11}\mu_{22}-\mu_{21}\mu_{12}] -\theta \Sigma'_1(\Sigma'_2-\mu_{22}). \label{tau'(2)}
\end{align}
Hence we have $\tau'(2)<0$ and $N=1$ if and only if  
    \begin{align*}
        -\Sigma_1(\Sigma_2)^2 \gamma^3+\Sigma_2(\Sigma_1\Sigma_2+\mu_{11}\mu_{22}-\mu_{21}\mu_{12}-\theta \Sigma_1)\gamma^2+\theta\Sigma_1\mu_{22}\gamma <0.
    \end{align*}
    Let $G(\gamma)=-\Sigma_1(\Sigma_2)^2 \gamma^2+\Sigma_2(\Sigma_1\Sigma_2+\mu_{11}\mu_{22}-\mu_{21}\mu_{12}-\theta \Sigma_1)\gamma+\theta\Sigma_1\mu_{22}$. The solutions of $G(\gamma)=0$ are 
    \begin{align*}
        \frac{\bigg(\Sigma_2+\frac{\mu_{11}\mu_{22}-\mu_{21}\mu_{12}}{\Sigma_1}-\theta \bigg)\pm \sqrt{\bigg(\Sigma_2+\frac{\mu_{11}\mu_{22}-\mu_{21}\mu_{12}}{\Sigma_1}-\theta\bigg)^2+4\theta \mu_{22}}}{2 \Sigma_2}.
    \end{align*}
    Since the smaller solution is negative, we have that $G(\gamma)<0$ if and only if
    \begin{align*}
        \gamma> \frac{\bigg(\Sigma_2+\frac{\mu_{11}\mu_{22}-\mu_{21}\mu_{12}}{\Sigma_1}-\theta \bigg) + \sqrt{\bigg(\Sigma_2+\frac{\mu_{11}\mu_{22}-\mu_{21}\mu_{12}}{\Sigma_1}-\theta\bigg)^2+4\theta \mu_{22}}}{2 \Sigma_2}.
    \end{align*}

\item[(ii)]From \eqref{tau'(2)},
 if $\lambda > \frac{\mu_{11}\mu_{22}-\mu_{21}\mu_{12}}{\Sigma_1\Sigma_2} $, then $\tau'(2)<0$. And if $\lambda = \frac{\mu_{11}\mu_{22}-\mu_{21}\mu_{12}}{\Sigma_1\Sigma_2} >0$,
 then our assumption that  $\Sigma_1,\Sigma_2>0$ in Section \ref{sec3} implies  $\theta \Sigma'_1(\Sigma'_2-\mu_{22})>0$ and hence $\tau'(2)<0$. Finally if $\lambda = \frac{\mu_{11}\mu_{22}-\mu_{21}\mu_{12}}{\Sigma_1\Sigma_2} =0$, from our assumption in Section \ref{sec3} that $\mu_{11}, \mu_{22}>0$, we have that $\mu_{12}>0$. As a result  $\theta \Sigma'_1(\Sigma'_2-\mu_{22})>0$ and again $\tau'(2)<0$. Hence we have $N=1$.
 \end{enumerate}
\end{proof}
Similarly, we also have a lower bound on the abandonment rate $\theta$ that guarantees that the expedite policy is optimal. Note that the bound in Proposition \ref{thetabound} will be negative, and hence the expedite policy is optimal for all $\theta>0$, if
\begin{align*}
    0>& \Sigma'_2 \mu_{11}-\Sigma_1' (\Sigma'_2-\mu_{22})\\
    =&(1+\lambda) \bigg[ (\mu_{22}+\mu_{12})\mu_{11}-(\mu_{11}+\mu_{21})(\lambda \Sigma_2+\mu_{12}) \bigg]\\
    =&(1+\lambda) \bigg[ -\lambda\Sigma_1 \Sigma_2+\mu_{11}\mu_{22}-\mu_{21}\mu_{12}\bigg].
\end{align*}
That is the bound can be negative only if $\lambda >\frac{\mu_{11}\mu_{22}-\mu_{21}\mu_{12}}{\Sigma_1\Sigma_2} $, and hence  Proposition \ref{thetabound} is consistent with Proposition \ref{lambda bound} in this case.

\begin{remark}
    Since $\Sigma_1\Sigma_2 = \mu_{11}\mu_{12} + \mu_{11}\mu_{22} + \mu_{21}\mu_{12} + \mu_{21}\mu_{22} \geq \mu_{11}\mu_{22} - \mu_{21}\mu_{12}$, from Proposition \ref{lambda bound} (ii) we see that if $\gamma \geq 2$, the expedite policy is always optimal.
\end{remark}

 \begin{proposition}\label{thetabound}
 \begin{enumerate}
     \item[(i)] If  $\theta > \frac{\Sigma'_2[\Sigma'_2\mu_{11} - \Sigma'_1(\Sigma'_2-\mu_{22})]}{\Sigma'_1(\Sigma'_2-\mu_{22}) }$, we have $N=1$, i.e., the expedite policy is the optimal policy.
     \item[(ii)] If $\theta > \frac{\Sigma_2(\mu_{11}\mu_{22}-\mu_{21}\mu_{12})}{\Sigma_1 \mu_{12}}$, we have $N=1$, irrespective of the value of $\gamma$.
 \end{enumerate}
       
   \end{proposition}
\begin{proof}
\begin{enumerate}
    \item[(i)] Follows from the proof of Proposition \ref{Prop N_TD theta} in Section \ref{Proposed policy} by setting $\gamma_1=\gamma_2=\gamma$.
    \item[(ii)] Let 
    \begin{align*}
        F(\gamma)=&\frac{\Sigma'_2[\Sigma'_2\mu_{11} - \Sigma'_1(\Sigma'_2-\mu_{22})]}{\Sigma'_1(\Sigma'_2-\mu_{22}) }
        =\frac{(\Sigma_2')^2\mu_{11}}{\Sigma_1'(\Sigma_2'-\mu_{22})}-\Sigma_2'
        =\frac{\gamma \Sigma_2^2\mu_{11}}{\Sigma_1(\gamma \Sigma_2-\mu_{22})}-\gamma \Sigma_2.
    \end{align*}
    Taking derivative with respect to $\gamma$ yields
    \begin{align*}
        F'(\gamma)=&\frac{\Sigma_1(\gamma \Sigma_2-\mu_{22})\Sigma_2^2 \mu_{11}-\gamma \Sigma_2^3\mu_{11} \Sigma_1}{[\Sigma_1(\gamma \Sigma_2-\mu_{22})]^2}-\Sigma_2
        =-\frac{\Sigma_2^2 \mu_{11}\mu_{22}}{\Sigma_1(\gamma \Sigma_2-\mu_{22})^2}-\Sigma_2<0
    \end{align*}
    That is $F(\gamma)$ is decreasing and the maximum value is attained at $\gamma=1$. Hence
    \begin{align*}
        \frac{\Sigma_2(\mu_{11}\mu_{22}-\mu_{21}\mu_{12})}{\Sigma_1 \mu_{12}} \geq \frac{\Sigma'_2[\Sigma'_2\mu_{11} - \Sigma'_1(\Sigma'_2-\mu_{22})]}{\Sigma'_1(\Sigma'_2-\mu_{22}) }.
    \end{align*}
\end{enumerate}
\end{proof}

Note that if $\Sigma'_2-\mu_{22}=0 $, from \eqref{tau'}-\eqref{f recursion} we have
\begin{align}\label{Sigma'=mu22}
\tau'(n)=&[\mu_{22}+(n-1)\theta][\mu_{22}+(n-2)\theta]\mu_{22}f(1,n-1) -[\mu_{22}+(n-1)\theta]\mu_{22}f(1,n) \nonumber \\
     & +[\mu_{22}+(n-2)\theta]\mu_{11}\mu_{22}f(1,n-1) \nonumber\\
     =&[\mu_{22}+(n-2)\theta]\mu_{11}\mu_{22}f(1,n-1).
\end{align}
Therefore $\tau'(n) \geq 0$ and from \eqref{Ndef} we have $N=B+2$, consistent with the result in Proposition \ref{thetabound}(i).

The next proposition compares the optimal buffer sizes for two systems with different synergistic factors $\gamma$ and $\gamma'$. In order to state the result, we modify our notation. In addition to $\lambda=\gamma-1$, we let 
\begin{align*}
    \lambda'=\gamma'-1
\end{align*} 
Furthermore, consistent with \eqref{tau'} and \eqref{Ndef}, we define for $\xi \in \{\gamma, \gamma'\}$
\begin{align*}
    \tau'_{\xi}(n)=&[\xi\Sigma_2+(n-1)\theta][\xi\Sigma_2+(n-2)\theta]\mu_{22}f(1,n-1) -[\xi\Sigma_2+(n-1)\theta][\xi\Sigma_1\alpha(n)(\xi\Sigma_2-\mu_{22})+\xi\Sigma_2f(1,n)] \nonumber \\
     &  +[\xi\Sigma_2+(n-2)\theta]\mu_{11}[\xi\Sigma_2f(1,n-1)+\xi\Sigma_1\alpha(n-1)(\xi\Sigma_2-\mu_{22})]
\end{align*}  
and 
\begin{align*}
    N_{\xi}=\max{\{n\in \{1,\ldots,B+2\}: \tau_{\xi}'(n) \geq 0\}}.
\end{align*}
The result of Proposition \ref{thm lamba_1< lambda_2} says that as the synergy between the servers increases, the buffer size used by the optimal policy decreases. This is intuitive because it implies that servers should take advantage of the greater synergistic factor by working in a team more often.

\begin{proposition}\label{thm lamba_1< lambda_2}
For $\gamma \geq \gamma' \geq 1$, fixed $\mu_{11}, \mu_{21}, \mu_{12}, \mu_{22}$, $\theta$,  and for $n=2, \ldots, B+2$ we have $\tau_{\gamma}'(n) \geq 0$ implies $\tau'_{\gamma'}(n) \geq 0$. That is the optimal buffer size satisfies $N_{\gamma} \leq N_{\gamma'}$.    
\end{proposition}
\begin{proof}
We first consider the case when $\gamma'=1$.

Some algebra, as shown in Supplementary Material \ref{SMProp11}, yields, 
\begin{align}\label{tau'(n)delta(n)}
    \tau_{\gamma}'(n)=(1+\lambda) \tau'_1(n)+\delta_{\lambda}(n),
\end{align}
where
\begin{align}\label{delta_lambda(n)}
    \delta_{\lambda}(n)
    =&-(1+\lambda)\lambda\Sigma_2 \delta'(n)-(1+\lambda)\lambda\Sigma_2^2(n-2)\theta f(1,n-1) \nonumber\\
    &-(1+\lambda)\lambda \Sigma_2[\lambda \Sigma_2+(n-2)\theta]\Sigma_1 [\alpha(n)-\mu_{11}\alpha(n-1)]-(1+\lambda)\lambda \Sigma_2\theta \Sigma_1\alpha(n) \nonumber\\
    &-\lambda(n-1)(n-2)\theta^2\mu_{22}f(1,n-1)
\end{align}
and 
\begin{align*}
  \delta'(n)=  \Sigma_1 \alpha(n)\mu_{12}-\mu_{11}\Sigma_2 f(1,n-1)-\Sigma_1 \alpha(n-1)\mu_{12}\mu_{11}+\Sigma_2\Sigma_1\alpha(n)-\Sigma_2\Sigma_1 \alpha(n-1)\mu_{11}.
\end{align*}
We want to show that $\delta_{\lambda}(n) \leq 0$. Since from \eqref{alphaineq} we have  $\alpha(n) \geq \mu_{11}\alpha(n-1)$,
it is sufficient to show that $\delta'(n) \geq 0$.
We show this by induction on $n$. For $n=2$ we have 
\begin{align*}
    \delta'(2)=&  \Sigma_1 \alpha(2)\mu_{12}-\mu_{11}\Sigma_2 f(1,1)-\Sigma_1 \alpha(1)\mu_{12}\mu_{11}+\Sigma_2\Sigma_1\alpha(2)-\Sigma_2\Sigma_1 \alpha(1)\mu_{11}\\
    =&\Sigma_1 \mu_{12}-\mu_{11}\Sigma_2+\Sigma_2\Sigma_1\\
    =&\Sigma_1\mu_{12}+\Sigma_2\mu_{21} \geq 0.
\end{align*}
Next we assume that $\delta'(n) \geq 0$ and show that it implies  $\delta'(n+1) \geq 0$. We have 
\begin{align*}
    \delta'(n+1)=&  \Sigma_1 \alpha(n+1)\mu_{12}-\mu_{11}\Sigma_2 f(1,n)-\Sigma_1 \alpha(n)\mu_{12}\mu_{11}+\Sigma_2\Sigma_1\alpha(n+1)-\Sigma_2\Sigma_1 \alpha(n)\mu_{11}.
\end{align*}
Now, from \eqref{f recursion}-\eqref{eq5},
\begin{align*}
   & \delta'(n+1)\\
   =& \Sigma_1 \mu_{12}[\mu_{22}+(n-1)\theta]\alpha(n)+\Sigma_1 \mu_{12}\mu_{11}^{n-1}-\mu_{11}\Sigma_2[\mu_{22}+(n-2)\theta]f(1,n-1)\\
    &-\Sigma_1\mu_{12}\mu_{11}[\mu_{22}+(n-2)\theta]\alpha(n-1)-\Sigma_1 \mu_{12}\mu_{11}^{n-1}\\
    &+\Sigma_2 \Sigma_1 [\mu_{22}+(n-1)\theta]\alpha(n)+\Sigma_2 \Sigma_1\mu_{11}^{n-1}-\Sigma_2 \Sigma_1 \mu_{11}[\mu_{22}+(n-2)\theta]\alpha(n-1) -\Sigma_2 \Sigma_1 \mu_{11}^{n-1}\\
    =& \Sigma_1 \mu_{12}[\mu_{22}+(n-1)\theta]\alpha(n)-\mu_{11}\Sigma_2[\mu_{22}+(n-2)\theta]f(1,n-1)-\Sigma_1\mu_{12}\mu_{11}[\mu_{22}+(n-2)\theta]\alpha(n-1)\\
    &+\Sigma_2 \Sigma_1 [\mu_{22}+(n-1)\theta]\alpha(n)-\Sigma_2 \Sigma_1 \mu_{11}[\mu_{22}+(n-2)\theta]\alpha(n-1) \\
    =& [\mu_{22}+(n-2)\theta] \bigg[ \Sigma_1 \mu_{12} \alpha(n)-\mu_{11}\Sigma_2 f(1,n-1)-\Sigma_1 \mu_{12}\mu_{11} \alpha(n-1)+\Sigma_2 \Sigma_1 \alpha(n)-\Sigma_2 \Sigma_1 \mu_{11} \alpha(n-1)\bigg]\\
    &+\theta \Sigma_1 \mu_{12} \alpha(n)+\theta \Sigma_2 \Sigma_1 \alpha(n)\\
    =& [\mu_{22}+(n-2)\theta] \delta'(n)+\theta \Sigma_1 \mu_{12} \alpha(n)+\theta \Sigma_2 \Sigma_1 \alpha(n) \geq 0,
\end{align*}
where the inequality follows from the induction assumption. 
Hence we have $\delta'(n) \geq 0$ for all $n \geq 2$ and therefore $\delta_{\lambda}(n) \leq 0$.
As a result from \eqref{tau'(n)delta(n)},
we have shown that $\tau'_{\gamma}(n) \geq 0$ implies $\tau'_1(n) \geq 0$ and therefore $N_{\gamma} \leq N_1$.

Next we consider the general case where the synergy for the first system is $\gamma$ and for the second system is $\gamma'$, where $\gamma \geq \gamma'$.
 From \eqref{tau'(n)delta(n)} we have
\begin{align*}
      \tau_{\gamma'}'(n)=(1+\lambda') \tau'_1(n)+\delta_{\lambda'}(n).
\end{align*}
After solving for $\tau'_1(n)$, \eqref{tau'(n)delta(n)} yields
\begin{align*}
    \tau'_{\gamma}(n)=&(1+\lambda)\bigg[\frac{\tau'_{\gamma'}(n)-\delta_{\lambda'}(n)}{(1+\lambda')} \bigg] +\delta_{\lambda}(n) =\frac{(1+\lambda)}{(1+\lambda')}\tau'_{\gamma'}(n)+\delta_{\lambda}(n)-\frac{(1+\lambda)}{(1+\lambda')}\delta_{\lambda'}(n).
\end{align*}
In order to show $\tau_{\gamma}'(n) \geq 0$ implies $\tau'_{\gamma'}(n) \geq 0$, it is sufficient to show that $\delta_{\lambda}(n)-\frac{(1+\lambda)}{(1+\lambda')}\delta_{\lambda'}(n) \leq 0.$
We have from \eqref{delta_lambda(n)} that 
\begin{align*}
     &\delta_{\lambda}(n)-\frac{(1+\lambda)}{(1+\lambda')}\delta_{\lambda'}(n)\\
     =&-(1+\lambda)(\lambda-\lambda')\Sigma_2 \delta'(n) -(1+\lambda)(\lambda-\lambda')\Sigma_2\Sigma_2(n-2)\theta f(1,n-1)\\
&-(1+\lambda)\Sigma_2\bigg[(\lambda^2-(\lambda')^2) \Sigma_2+(\lambda-\lambda')(n-2)\theta \bigg]\Sigma_1 [\alpha(n)-\mu_{11}\alpha(n-1)]\\
    &-(1+\lambda)(\lambda-\lambda') \Sigma_2\theta \Sigma_1\alpha(n)-\bigg[\lambda-\frac{1+\lambda}{1+\lambda'}\lambda'\bigg](n-1)(n-2)\theta^2\mu_{22}f(1,n-1)\leq 0,
\end{align*}
where in the last inequality we have used $\delta'(n) \geq 0$ , $\alpha(n) \geq \mu_{11}\alpha(n-1)$, and $\lambda \geq \lambda'$.

Hence, $\tau_{\gamma}'(n) \geq 0$ implies $\tau'_{\gamma'}(n) \geq 0$ and therefore $N_{\gamma} \leq N_{\gamma'}$. 
\end{proof}
We conclude this section by stating a proposition which compares the optimal buffer size as the abandonment rate changes (including no abandonments). We let $N^{(\theta)}$ denote the optimal buffer size  when the abandonment rate is $\theta$. The following proposition states that the optimal buffer size decreases as the abandonment rate increases. This result is expected since with abandonments the work done in station 1 is wasted. Therefore, in order to maximize the throughput, the optimal policy will try to reduce the number of abandonments by reducing the number of customers waiting in line to start service in station 2.
\begin{proposition}\label{N decreases with theta}
    For $\theta_1 \geq \theta_2 \geq 0$, fixed $\mu_{11}, \mu_{21}, \mu_{12}, \mu_{22}$, $\gamma$, we have the optimal buffer size $N^{(\theta_1)} \leq N^{(\theta_2)}$.
\end{proposition}
\begin{proof}
First we consider the scenario where $\theta_1>\theta_2=0$.
Here we consider two cases.
\begin{itemize}
    \item Case 1 : $\gamma-1=\lambda >\frac{\mu_{11}\mu_{22}-\mu_{21}\mu_{12}}{\Sigma_1 \Sigma_2} $. In this case, Proposition \ref{lambda bound} yields that $N^{(\theta_1)}=1$. Moreover, we have
$$\gamma=1+\lambda > 1+ \frac{\mu_{11}\mu_{22}-\mu_{21}\mu_{12}}{\Sigma_1 \Sigma_2} =\frac{\mu_{11}}{\mu_{11}+\mu_{21}}+\frac{\mu_{22}}{\mu_{22}+\mu_{12}}.$$
From Theorem 3.1 in \cite{SHD} we have that in the absence of abandonments a policy is optimal if and only if it
is nonidling and fully collaborative, i.e., $N^{(0)}=1$. We have shown that $N^{(\theta_1)}=N^{(0)}$.
\item Case 2 : $ \gamma-1=\lambda \leq \frac{\mu_{11}\mu_{22}-\mu_{21}\mu_{12}}{\Sigma_1 \Sigma_2}$. In this case

$$\gamma=1+\lambda \leq \frac{\mu_{11}}{\mu_{11}+\mu_{21}}+\frac{\mu_{22}}{\mu_{22}+\mu_{12}}$$
and hence from Theorem 3.1 in \cite{SHD} we have that in the absence of abandonments $N^{(0)}=B+2$. Hence from definition of $N$ we have $N^{(\theta_1)} \leq N^{(0)}$ for $\theta_1>0$.

\end{itemize}

Next we consider the scenerio $\theta_1 \geq \theta_2>0$. Note that if $\Sigma'_2-\mu_{22} =0$, we have from \eqref{epsilon'(n)} that $\epsilon'(n)=0$ and therefore \eqref{tau(n+1)} and \eqref{Ndef} imply that $N^{(\theta_1)}=B+2$ for any $\theta_1$ and the result holds trivially. Also, if $\mu_{11}\mu_{22}-\mu_{21}\mu_{12}=0$ or $\lambda \geq \frac{\mu_{11}\mu_{22}-\mu_{21}\mu_{12}}{\Sigma_1\Sigma_2} $, we have from Proposition \ref{lambda bound} that $N^{(\theta_1)}=1$ for all $\theta_1>0$ and again the result holds trivially. Therefore we assume that $\Sigma'_2-\mu_{22}>0$, $\mu_{11}\mu_{22}-\mu_{21}\mu_{12}>0$ and $\lambda < \frac{\mu_{11}\mu_{22}-\mu_{21}\mu_{12}}{\Sigma_1\Sigma_2} $ . We write $\tau'(n,\theta)$ to specify the $\theta$ used in the expression for $\tau'(n)$ and show that $\tau'(n,\theta)$ has a single positive root $\theta(n)$. For $n \geq 2$ we have from \eqref{tau(n+1)}, \eqref{epsilon'(n)},
and $\eqref{tau'(2)}$ that
\begin{align}\label{tau'(n,0)}
    \tau'(n,0)=\mu_{22}^{n-2}\tau'(2,0)
    >0,
\end{align}
where the last inequality follows from our assumption that $\lambda < \frac{\mu_{11}\mu_{22}-\mu_{21}\mu_{12}}{\Sigma_1\Sigma_2} $. Since the highest power of $\theta$ in $f(1,n-1)$ and $\alpha(n-1)$ is $n-3$, we have from \eqref{epsilon'(n)} that  the highest power of $\theta$ in $\epsilon'(n-1)$ is $n-1$ and from \eqref{tau'(2.2)}
\begin{align*}
    \tau'(2,\theta)=-\Sigma'_2 \Sigma'_1(\Sigma'_2-\mu_{22})+(\Sigma'_2)^2\mu_{11} -\theta \Sigma'_1(\Sigma'_2-\mu_{22}),
\end{align*}
it follows from \eqref{tau(n+1)} that the highest power of $\theta$ in $\tau'(n,\theta)$ is $\theta^{n-1}$ and it has a negative coefficient. Hence, 
\begin{align*}
    \lim_{\theta \to \infty} \tau'(n, \theta)=-\infty,
\end{align*}
and therefore for each $n \geq 2$, $\tau'(n,\theta)$ has a positive root. We show that it has exactly one positive root. We do this by showing that $g_n(\theta)$, the gain under policy $d_n$ when the abandonment rate is $\theta$, has the following two properties
\begin{enumerate}
    \item[(1)] $g_n(\theta)$ is non-increasing in $\theta$ for $n \geq 1$. 
    \item[(2)] As $\theta$ increases, $g_n(\theta)$ decreases at least as rapidly as $g_{n-1}(\theta)$ for $n \geq 2$.
\end{enumerate}
Assuming these properties, suppose there exist two positive roots of $\tau'(n,\theta)$, say $ \theta^1(n)$  and $ \theta^2(n)$, with  $0 < \theta^1(n)< \theta^2(n)$, that is, $\tau'(n, \theta^1(n))=\tau'(n, \theta^2(n))=0$. From Proposition \ref{tau' proposition} this implies that $g_n(\theta^1(n))=g_{n-1}(\theta^1(n))$ and $g_n(\theta^2(n))=g_{n-1}(\theta^2(n))$. However, from, property $(2)$ above, if $g_n(\theta^1(n))=g_{n-1}(\theta^1(n))$  and since $\theta^1(n)< \theta^2(n)$ we must have $g_n(\theta^2(n))\leq g_{n-1}(\theta^2(n))$. Further, the inequality should be strict because otherwise we would have $g_n(\theta)=g_{n-1}(\theta)$ $\forall \theta^1(n) \leq \theta \leq \theta^2(n)$, which cannot be true since $\tau'(n,\theta)$ is a polynomial in $\theta$ and hence cannot satisfy $\tau'(n,\theta )=0$ $\forall \theta^1(n) \leq \theta \leq \theta^2(n)$. 

We have shown that there exists a single $\theta(n) >0$ such that $\tau'(n, \theta(n))=0$. From \eqref{tau(n+1)} this implies that 
\begin{align*}
       \tau'(n+1,\theta(n))=[\mu_{22}+(n-1)\theta]\tau'(n,\theta(n))-\epsilon'(n,\theta(n))<0,
\end{align*}
where the last inequality follows from \eqref{epsilon'(n)} and $\Sigma'_2-\mu_{22}>0$. 
This together with the fact that $\tau'(n,0)>0$ from \eqref{tau'(n,0)} implies that $\theta(n+1)<\theta(n)$. For convenience, define $\theta(1)=\infty$ and $\theta(B+3)=0$. Therefore for $\theta_1 \geq \theta_2>0$ we can have the following scenerios :
\begin{enumerate}
      \item [(i)] $\exists n \in \{1, \ldots, B+2\}$ such that $\theta(n+1) < \theta_2 \leq \theta_1 \leq \theta(n)$.\\
    In this case we have $N^{(\theta_1)}=N^{(\theta_2)}=n$.
    \item[(ii)] $\exists n \in \{2, \ldots,B+2\}$ such that $\theta_2 \leq \theta(n) < \theta_1$.\\
    In this case we have $N^{(\theta_1)}<n$ and $N^{(\theta_2)} \geq n$.
\end{enumerate}
Hence, in both cases we have $N^{(\theta_1)} \leq N^{(\theta_2)}$. 

It remains to show that properties $(1)-(2)$ hold for $g_n(\theta)$. We start by showing property $(1)$. Note that from Proposition \ref{tau' proposition}, $g_1(\theta)=\frac{\Sigma'_1 \Sigma'_2}{\Sigma'_1+\Sigma'_2}$, that is, $g_1(\theta)$ is independent of $\theta$. Now consider two systems, System 1 and System 2,  where both the systems are operating under the policy $d_n$ but the abandonment rate for System 1 is $\theta^1$ and that of System 2 is $\theta^2$ with $\theta^1 > \theta^2$. Consider sample paths of both systems starting in the same state and generated
using (the same) independent Poisson processes with rates $\Sigma'_1, \Sigma'_2, (n-1)\theta^1$,
and thinning. Due to the higher abandonment rate of
System 1, the number of customers in System 2 will be always greater than or equal to the number of customers in System 1. In particular, when both systems are in the same state, every abandonment from
(service completion at station 1 in) System 2 is associated with a simultaneous abandonment
from (service completion at station 1 in) System 1. Furthermore, every service completion at station 2 in System 1 is associated with a simultaneous service completion at station 2 in
System 2. This shows that more customers complete service at station 2 in System 2 and
$g_n(\theta)$ is non-increasing in $\theta$.

Next we show property $(2)$. Consider two systems, System 1 and System 2, now with the same abandonment rate $\theta$, but System 1 operates under policy $d_n$ and System 2 under policy $d_{n-1}$. Consider sample paths of both systems generated
using (the same) independent Poisson processes with rates $\Sigma'_1, \Sigma'_2, (n-1)\theta^1$,
and thinning and starting from the same state $s \leq n-1$. Both systems evolve identically until state $n-1$ is reached. If there is an abandonment
in state $n-1$ both systems experience the abandonment, but it is possible that in state
$n-1$, $(i)$ System 1 (but not System 2) can transition to state $n$ due to a service completion at
station 1 or $(ii)$ System 2 (but not necessarily System 1) will transition to state $n-2$ due to a
service completion at station 2. If this happens, then from this point onward, the number of
customers in System 1 will dominate the number of customers in System 2. As a result  every abandonment in System 2 is associated with a simultaneous abandonment in System 1. However, when the two systems are not in the same state, there may be abandonments in System 1
that do not occur in System 2. Some of these additional abandonments in System 1 may be
replaced by additional service completions at station 1 in System 1 but not in System 2 (if
both systems are in state $n-1$); others may not be replaced. As $\theta$ increases, these unreplaced additional abandonments also increase and since all service completions at
station 1 lead to either abandonments or service completions at station 2, this implies that
as $\theta$ increases, $g_n(\theta)$ decreases at least as rapidly as $g_{n-1}(\theta)$. This proves the two properties completing the proof.
\end{proof}
\section{Numerical Results}\label{Numerical Results same gamma}

In this section, we provide a numerical study to illustrate the theoretical results derived in Sections \ref{sec3} and \ref{sec4}. Specifically, we examine the behavior of the long-run average throughput $g_n$ as a function of $n$, see how the optimal buffer size $N$ changes as a function of the synergistic factor $\gamma$ and abandonment rate $\theta$, analyze the distribution of the optimal buffer sizes across a wide range of parameter sets, and evaluate the relative performance of the optimal policy compared to the simple expedite policy for different domains of the system parameters.

\begin{figure}[H]
\centering

\begin{minipage}{0.45\linewidth}
    \centering
    \includegraphics[width=\linewidth]{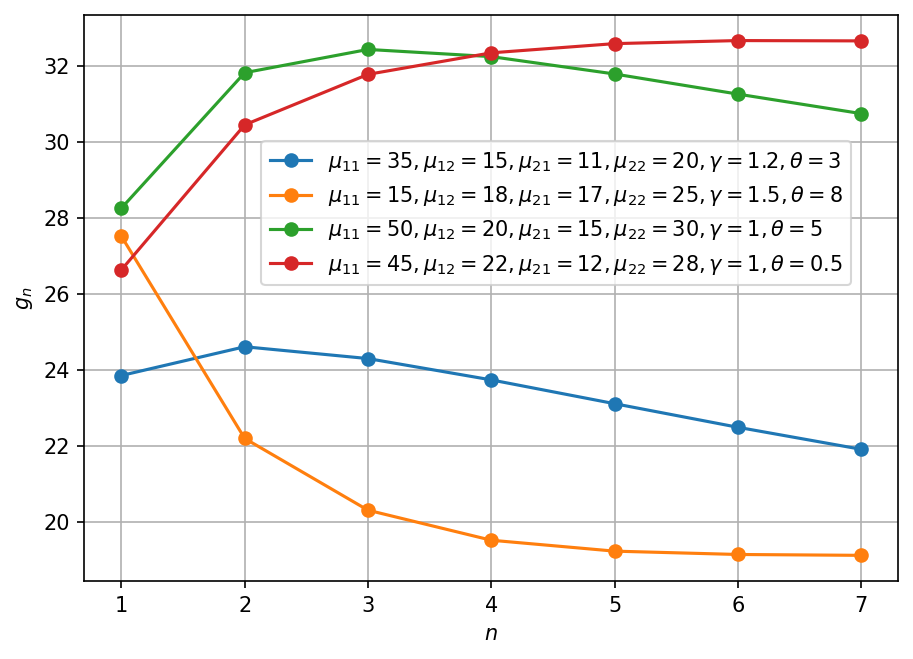}
    \caption{Behavior of the long-run average throughput $g_n$ as a function of $n$ $(B=5)$.}
    \label{g_n vs n}
\end{minipage}
\hfill
\begin{minipage}{0.45\linewidth}
    \centering
    \includegraphics[width=\linewidth]{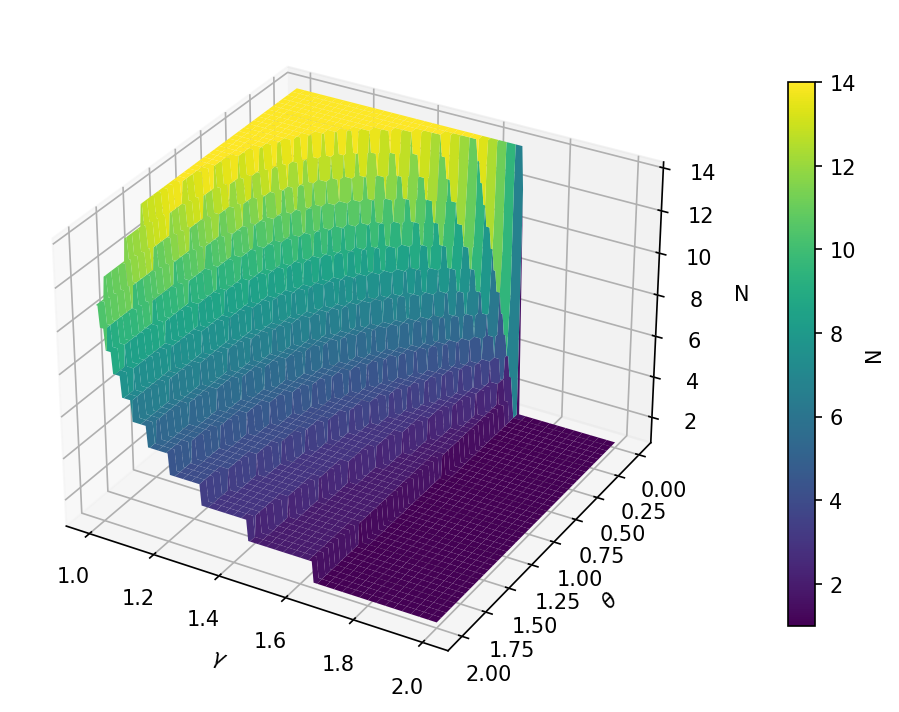}
    \caption{Behavior of $N$ as a function of $\gamma \in [1,2]$ and $\theta \in [0,2]$.}
    \label{N as a function of gamma and theta}
\end{minipage}

\end{figure}

First, we investigate how the system throughput $g_n$ varies with respect to $n$. In Proposition \ref{tau property}, we showed that  there exists an $N$ such that  $g_N \geq g_{N-1} \geq \cdots \geq g_1$ and $g_N > g_{N+1} \geq \cdots \geq g_{B+2}$, i.e, the long-run average throughput $g_n$ for policy $d_n$ (as defined in \eqref{d_n}) first increases with respect to the buffer size, reaches a maximum, and then decreases. Figure \ref{g_n vs n} illustrates this behavior for different sets of parameters. Note that $N$ can be equal to 1, in which case the throughput always decreases in $n$, or $N$ can  be equal to $B+2$ in which case the throughput always increases in $n$ (as shown in Figure \ref{g_n vs n}).

We next examine how the optimal buffer size behaves as a function of the synergistic factor $\gamma$ and abandonment rate $\theta$. From Propositions \ref{thm lamba_1< lambda_2} and \ref{N decreases with theta}, we know that the optimal buffer size $N$ decreases with respect to both $\theta$ and $\gamma$. We observe this behavior in Figure \ref{N as a function of gamma and theta}, where  $\mu_{11}=10,\mu_{12}=2,\mu_{21}=2,\mu_{22}=13$, and $B=12$. Over the plotted range of ($\gamma$, $\theta$), the optimal buffer size decreases from the upper boundary $B+2$ (here, 
14) for small $\gamma$ and $ \theta$ to the lower boundary 1 for sufficiently large $\gamma$ and $\theta$.
We also observe that when $\theta$ is close to 0, $N$ exhibits a sharp drop once $\gamma$ crosses a threshold. However, the drop is more gradual as $\theta$ increases.
To understand the prevalence of different optimal policies in a general setting, we conducted a numerical experiment consisting of 10,000 randomly generated parameter instances from the ranges given in Table \ref{rangetable same gamma}. For each instance, we computed the optimal buffer size $N$ that maximizes the throughput. The frequency of different $N$ values is given in Table \ref{Nfrequency same gamma}

\begin{table}[ht]
    \centering
    \begin{minipage}[t]{0.3\textwidth}
        \vspace{0pt} 
        \centering
        \begin{tabular}{ |c|c| } 
            \hline
            Parameter & Range\\
            \hline \hline
            $\mu_{ij}$ & [0,400]\\
            \hline
            $\gamma$ & [1,2]\\
            \hline
            $\theta$ & [0,10]\\
            \hline
            $B$ & $\{0,\ldots,100\}$\\
            \hline
        \end{tabular}
        \caption{Ranges of Parameters.}
        \label{rangetable same gamma}
    \end{minipage}%
    \hfill 
    \begin{minipage}[t]{0.68\textwidth}
        \vspace{0pt} 
        \centering
        \footnotesize 
        \scalebox{0.9}{%
        \begin{tabular}{ |c|c| } 
        \hline
        $ N$ value & Frequency\\
        \hline \hline
        1 & 7362\\
        \hline
        2 &338\\
        \hline
        3 & 354\\
        \hline
        4 & 348\\
        \hline
        5 & 256\\
        \hline
        6&182\\
        \hline
        7&158\\
        \hline
        8&128\\
        \hline
        9&104\\
        \hline
        10&82 \\
        \hline
        11&69\\
        \hline
        12&54\\
        \hline
        13&47\\
        \hline
        14&46\\
        \hline
        15&44\\
        \hline
        16&30\\
        \hline
        17&29\\
        \hline
        18&27\\
        \hline
        19&25\\
        \hline
        20&19\\
        \hline
        21&19\\
        \hline
        22&20\\
        \hline
        23&17\\
        \hline
        24&13\\
        \hline
        25&15\\
        \hline
        26&12\\
        \hline
        27&8\\
        \hline
        28&5\\
        \hline
        29&10\\
        \hline
        30&5\\
        \hline
        \end{tabular}
        \quad
        \begin{tabular}{ |c|c|}
        \hline
        $ N$ Value & Frequency\\
        \hline \hline
        31&6\\
        \hline
        32&5\\
        \hline
        33&5\\
        \hline
        34&7\\
        \hline
        35&5\\
        \hline
        36&6\\
        \hline
        37&8\\
        \hline
        38&2\\
        \hline
        39&4\\
        \hline
        40&3\\
        \hline
        41&5\\
        \hline
        42&4\\
        \hline
        43&8\\
        \hline
        44&3\\
        \hline
        45&7\\
        \hline
        46&3\\
        \hline
        47&2\\
        \hline
        48&4\\
        \hline
        49&4\\
        \hline
        50&2\\
        \hline
        51&3\\
        \hline
        52&2\\
        \hline
        53&4\\
        \hline
        54&2\\
        \hline
        55&3\\
        \hline
        56&2\\
        \hline
        57&1\\
        \hline
        58&3\\
        \hline
        59&3\\
        \hline
        60&1\\
        \hline
        \end{tabular}
        \quad
        \begin{tabular}{ |c|c|}
        \hline
        $ N$ Value & Frequency\\
        \hline \hline
        61&2\\
        \hline
        62&5\\
        \hline
        63&1\\
        \hline
        64&1\\
        \hline
        66&1\\
        \hline
        67&3\\
        \hline
        69&1\\
        \hline
        70&4\\
        \hline
        71&4\\
        \hline
        72&1\\
        \hline
        73&1\\
        \hline
        74&2\\
        \hline
        75&1\\
        \hline
        76&1\\
        \hline
        79&3\\
        \hline
        80&5\\
        \hline
        82&2\\
        \hline
        83&3\\
        \hline
        85&1\\
        \hline
        88&1\\
        \hline
        91&3\\
        \hline
        92&3\\
        \hline
        94&2\\
        \hline
        96&2\\
        \hline
        97&1\\
        \hline
        99&2\\
        \hline
        100&1\\
        \hline
        \end{tabular}}
        \caption{Frequency of $N$ values.}
        \label{Nfrequency same gamma}
    \end{minipage}
\end{table}
The results given in Table \ref{Nfrequency same gamma} indicate that the expedite policy ($N=1$) is optimal in a significant proportion of the generated cases $(73.62\%)$. This could be due in part to the fact that the possible values of $ N$ for a buffer size of $B$ are $\{1,2,\ldots,B+2\}$, and therefore for the system to have a high $N$ value, the randomly chosen $B$ should also be high.

Given the high frequency of the expedite policy's optimality observed in Table \ref{Nfrequency same gamma}, it is practically useful to determine the ``penalty" for using the expedite policy when it is not optimal. We compare the optimal throughput $g_N$  to the throughput of the expedite policy $g_1$. We perform 10,000 random replications of the parameters from different parameter ranges as given in Table \ref{Expedite vs optimal throughput} and estimate the relative performance of the expedite policy via the ratio
\begin{align*}
    \frac{\sum_{i=1}^{10,000}{g_1^{(i)}}}{\sum_{i=1}^{10,000} g_N^{(i)}},
\end{align*}
where $g_1^{(i)}$ and $g_N^{(i)}$ correspond to the throughputs of the expedite and the optimal policies, respectively, for the $i^{th}$ parameter instance.

\begin{center}
 \begin{table}[t]
\centering
\renewcommand{\arraystretch}{1.2}
\setlength{\tabcolsep}{8pt}
\scalebox{0.72}{%
\begin{tabular}{|c|c|c|c|c|c|c|}
\hline
\multirow{3}{*}{$\gamma$}
& \multicolumn{6}{c|}{Ranges of Parameters} \\
\cline{2-7}
& \multicolumn{2}{c|}{\thead{$0 \leq \mu_{11},\mu_{22} \leq 400 $,\\ $ 0 \leq \mu_{12},\mu_{21}\leq 400$}} & \multicolumn{2}{c|}{\thead{$200 \leq \mu_{11},\mu_{22} \leq 400 $,\\ $ 0 \leq \mu_{12},\mu_{21}\leq 200$}} & \multicolumn{2}{c|}{\thead{$200 \leq \mu_{11},\mu_{22} \leq 400 $,\\ $ 0 \leq \mu_{12},\mu_{21}\leq 20$}} \\
\cline{2-7}
& $ 0 \leq \theta \leq 10$ & $ 0 \leq \theta \leq 1$ & $ 0 \leq \theta \leq 10$ & $ 0 \leq \theta \leq 1$ & $ 0 \leq \theta \leq 10$ & $ 0 \leq \theta \leq 1$ \\
\hline\hline

1.0 & 0.856  & 0.839&  0.722 & 0.701 &0.588 & 0.570\\
\hline
1.1 & 0.909  & 0.896  & 0.782  & 0.761 &0.638 & 0.618 \\
\hline
1.2 & 0.947  & 0.938  & 0.839 & 0.820 &0.686 & 0.665\\
\hline
1.3 & 0.971  & 0.966 &  0.892  & 0.876&0.733 & 0.712\\
\hline
1.4 & 0.986  & 0.983 & 0.936 & 0.924 &0.780 & 0.758\\
\hline
1.5 & 0.994  & 0.992  & 0.967  & 0.960 &0.826 & 0.805\\
\hline
1.6 & 0.997  & 0.997 & 0.986 & 0.982&0.871 & 0.851\\
\hline
1.7 & 0.999  & 0.998 & 0.995 &  0.994 &0.914 & 0.897\\
\hline
1.8 & 0.999 & 0.999  & 0.999 &  0.998 &0.955 & 0.942\\
\hline
1.9 & 0.999 & 0.999 & 0.999 &  0.999 &0.990 & 0.985\\
\hline
2.0 & 1.0 & 1.0 &  1.0 & 1.0 &1.0 & 1.0\\
\hline

\end{tabular}}

\captionof{table}{Relative performance of the expedite policy vs. the optimal policy $(B=100).$
}
\label{Expedite vs optimal throughput}
\end{table}
\end{center}


Table~\ref{Expedite vs optimal throughput} highlights the patterns in the relative performance of the expedite policy. Firstly, we observe that for each fixed choice of the service-rate ranges (i.e., within each of the three $\mu_{ij}$ range blocks), the column with $0\leq \theta \leq 1$ is consistently lower than the corresponding column with
$0\leq \theta \leq 10$ across all values of $\gamma$. Thus, restricting the abandonment rate $\theta$ to the smaller range  makes the expedite policy marginally less competitive relative to the optimal policy. This behavior is expected from Proposition \ref{N decreases with theta}. Further, we can observe the dependence of performance on the degree of server specialization. In parameter regimes where the service rates across stations are broadly comparable (e.g., $\mu_{11}, \mu_{22}$ drawn from a similar range as $\mu_{12},\mu_{21}$), the expedite policy performs very well, even for moderate values of the synergistic factor $\gamma$. In contrast, when servers are highly specialized (characterized by large $\mu_{11}, \mu_{22}$ as compared to $\mu_{12}, \mu_{21}$), the performance of the expedite policy deteriorates substantially, especially for low to moderate values of $\gamma$. In these regimes, the expedite policy may divert capacity away from tasks for which servers are significantly more efficient, leading to pronounced throughput losses. This reinforces the practical relevance of computing or approximating the optimal policy in such systems. While the expedite policy offers simplicity and performs well in many cases, Table~\ref{Expedite vs optimal throughput} shows that its use can entail substantial throughput penalties in structurally asymmetric systems.

\section{Task Dependent Synergy}\label{sec5}
In this section, we generalize the models in Sections \ref{sec2}, \ref{sec3}, and \ref{sec4} by allowing different synergistic factors for tasks processed
at different stations. That is,  collaboration is more effective at some stations than others. 
We first show in Section \ref{sec5.1} that in the case of systems with $N$ stations and $M$ generalist servers, the expedite policy is the optimal policy. Next for systems with 2 servers and 2 stations and non-generalist servers, we propose a dynamic server assignment policy in Section \ref{Proposed policy} and provide numerical results that suggest that the proposed policy is indeed the policy that achieves the maximum throughput in Section \ref{Numerical Results}.
\subsection{Generalist Servers with $N$ Stations and $M$ Servers}\label{sec5.1}
We consider the same setting as in Section \ref{sec2} but now the synergy among the servers is allowed to be task (station) dependent. We show that even under this generalization, the expedite policy continues to be optimal.

\begin{theorem}
Assume that for all $ t \geq 0$, if there is a customer in service at station $j$ at time $t$, then the expected remaining service requirement at station $j$ of that customer is bounded by a scalar $1 \leq \bar{S} < \infty$.
 When $\mu_{ij}=\mu_i \delta_j$ for all $i=1, \ldots, M$ and $j=1, \ldots, N$ and $\gamma_j \geq 1 $ is the synergistic factor at station $j=1, \ldots,N$, then the expedite policy is optimal with long-run average throughput $ g_{\pi^*}=\frac{ \sum_{i=1}^M\mu_i}{\sum_{j=1}^N \frac{1}{\delta_j\gamma_j}}$.
\end{theorem}

\begin{proof}
Consider a model $\mathcal{H}$ where the service rate of server $i$ at station $j$ is $\mu_{ij}=\mu_i\delta_j\gamma_j$ and there is no synergy among the servers. This is equivalent to a model where the service requirements of successive customers at station $j \in \{1, \ldots, N\}$ are iid with mean $\frac{1}{\delta_j\gamma_j}$, the service rates depend only on the server with $\mu_{ij}=\mu_i$ $\forall i \in \{ 1, \ldots. M\}$, and $\gamma=1$. Now suppose the servers can also choose to operate at full efficiency at station $j$ with $\mu_{ij}=\mu_i$ or partial efficiency with  $\mu_{ij}=\mu_i/\gamma_j$. Let $\pi_+^*$ be the policy where all servers always work in a single team at full efficiency, and the team follows each customer from the first to the last station and only starts work on a new customer once all work on the previous customer has been completed. Since \eqref{pi work performed} is still satisfied (with $\gamma=1$) for $\mathcal H$ (because partial efficiency can only reduce the total work performed), we have from the proof of Theorem \ref{thm generalist} with $\gamma=1$ and $\delta_j$ replaced by $\delta_j \gamma_j$ for $j \in \{1, \ldots, N\}$  that $\pi_+^*$ is the optimal policy in model $\bar{\mathcal{H}}$ with the long-run average throughput $g_{\pi^*_+}^{\mathcal{H}}=\frac{  \sum_{i=1}^M\mu_i}{\sum_{j=1}^N \frac{1}{\delta_j\gamma_j}}.$
That is, if $\pi'$ is any other policy in $\mathcal H$ and $g_{\pi'}^{\mathcal{H}}$ is the long-run average throughput under that policy we have $g^{\mathcal{H}}_{\pi^*_+} \geq g^{\mathcal{H}}_{\pi'}.$

Now consider our original model where the service rate of server $i$ at station $j$ is $\mu_{ij}=\mu_i \delta_j$ and the synergistic factor at station $j$ is $\gamma_j$. Let $\pi^*$ be the expedite policy. Since all servers work in a single team at all times, the long-run average throughput for the expedite policy in this model will be the same as the long-run average throughput for the $\pi_{+}^*$ policy in model $\mathcal{H}$.
%
That is we have $  g_{\pi^*}=g^{\mathcal{H}}_{\pi^*_+}.$

Now for any policy $\pi$ in the original model we can get a policy $\pi'$ in our modified model ${\mathcal{H}}$ such that $g^{\mathcal{H}}_{\pi'}=g_{\pi}.$
This is because when a server $i$ is working on his own at a station $j$ in $\pi$, we can make the server $i$ work at station $j$ with partial efficiency in the modified model ${\mathcal{H}}$. And since we have shown that for all $\pi'$
\begin{align*}
    g_{\pi^*}=g^{\mathcal{H}}_{\pi^*_+} \geq g^{\mathcal{H}}_{\pi'}=g_{\pi},
\end{align*}
we have that the expedite policy is the optimal policy for our original model. 
\end{proof}

\subsection{Proposed Policy for 2 Stations and 2 Servers}\label{Proposed policy}
We consider the same settings as in Section \ref{sec3} except that the synergy is task (station) dependent, that is, the synergistic factor for station 1 is $\gamma_1$ and that for station 2 is $\gamma_2$. For $j=1,2$, let
\begin{align*}
    \bar \Sigma'_j=\gamma_j(\mu_{1j}+\mu_{2j})=\gamma_j\Sigma_j
\end{align*}
and adapt \eqref{tau'} and \eqref{Ndef} as follows
\begin{align}\label{bar{tau}}
    \bar{\tau}'(n)=&[\bar\Sigma'_2+(n-1)\theta][\bar\Sigma'_2+(n-2)\theta]\mu_{22}f(1,n-1) -[\bar\Sigma'_2+(n-1)\theta][\bar\Sigma'_1\alpha(n)(\bar\Sigma'_2-\mu_{22})+\bar\Sigma'_2f(1,n)] \nonumber \\
     & +[\bar\Sigma'_2+(n-2)\theta]\mu_{11}[\bar\Sigma'_2f(1,n-1)+\bar\Sigma'_1\alpha(n-1)(\bar\Sigma'_2-\mu_{22})];
\end{align}
\begin{align*}
    \bar{N}=\max{\{n \in \{1,2, \ldots, B+2\}: \bar{\tau}'(n) \geq 0\}}.
\end{align*}
Our proposed policy is very similar to the optimal policy $d_N$ defined in Section \ref{sec3}.
In particular, we define decision rule $d_{\bar{N}}$ to be $d_n$ as in \eqref{d_n} with $n=\bar{N}$.
The immediate reward $\bar{r}_{d_{\bar N}}$, transition probability matrix $\bar{P}_{d_{\bar N}}$ and the gain $\bar{g}_{\bar{N}}$ under this policy are given as in \eqref{r_{d_N}}, \eqref{P_{d_N}}, and \eqref{g_n}, respectively, with $\Sigma'_i$ replaced by $\bar\Sigma'_i$ for $i=1,2$ and $N$ replaced by $\bar N$.
The next proposition gives a lower bound on $\theta$, that guarantees $\bar{N}=1$.
 \begin{proposition}\label{Prop N_TD theta}
If  $\theta > \frac{\bar\Sigma'_2[\bar\Sigma'_2\mu_{11} - \bar\Sigma'_1(\bar\Sigma'_2-\mu_{22})]}{\bar\Sigma'_1(\bar\Sigma'_2-\mu_{22}) }$, we have $\bar N=1$.
   \end{proposition}
\begin{proof}
From \eqref{bar{tau}}
we have
\begin{align*}
    \bar{\tau}'(2)=&[\bar\Sigma'_2+\theta]\bar\Sigma'_2\mu_{22}-[\bar\Sigma'_2+\theta][\bar\Sigma'_1(\bar\Sigma'_2-\mu_{22})+\bar\Sigma'_2\mu_{22}] +\mu_{11}(\bar\Sigma'_2)^2\\
     =&-[\bar\Sigma'_2+\theta]\bar\Sigma'_1(\bar\Sigma'_2-\mu_{22}) +\mu_{11}(\bar\Sigma'_2)^2 \\
     =&-\bar\Sigma'_2 \bar\Sigma'_1(\bar\Sigma'_2-\mu_{22})+(\bar\Sigma'_2)^2\mu_{11} -\theta \bar\Sigma'_1(\bar\Sigma'_2-\mu_{22}) .
\end{align*}
When $\theta > \frac{\bar\Sigma'_2[\bar\Sigma'_2\mu_{11} - \bar\Sigma'_1(\bar\Sigma'_2-\mu_{22})]}{\bar\Sigma'_1(\bar\Sigma'_2-\mu_{22}) }$, we have $\bar{\tau}'(2) <0$. 
Exactly as in the proof of Proposition \ref{tau property}, it can be shown that  if $\bar{\tau}'(n) \geq 0$, then $\bar{\tau}'(k) \geq 0$ for all $1 \leq k \leq n-1$. Similarly if $\bar{\tau}'(n) \leq 0$, then $\bar{\tau}'(k) \leq 0 $ for all $n+1 \leq k \leq B+2$. Therefore, $\bar{N}=1$.
\end{proof}
Note that if $\bar\Sigma'_2-\mu_{22}=0$ (i.e. $\mu_{12}=0$ and $\gamma=1$), from \eqref{bar{tau}} we have
\begin{align*}
    \bar\tau'(n)=&[\mu_{22}+(n-1)\theta][\mu_{22}+(n-2)\theta]\mu_{22}f(1,n-1) -[\mu_{22}+(n-1)\theta]\mu_{22}f(1,n) \nonumber \\
     & +[\mu_{22}+(n-2)\theta]\mu_{11}\mu_{22}f(1,n-1) \\
     \geq& 0,
\end{align*}
where the inequality follows from \eqref{Sigma'=mu22}, and hence we have $\bar N=B+2$, consistent with the result of Proposition \ref{Prop N_TD theta}.
\subsection{Numerical Results for 2 Stations and 2 Servers}\label{Numerical Results}
\begin{figure}[H]
    \centering
    \begin{minipage}[c]{0.45\textwidth}
        \centering
        \begin{tabular}{ |c|c| } 
            \hline
            Parameter & Range\\
            \hline \hline
            $\mu_{11},\mu_{22},\mu_{12},\mu_{21}$ & [0,400]\\
            \hline
            $\gamma_1, \gamma_2$ & [1,2]\\
            \hline
            $\theta$ & [0,10]\\
            \hline
            $B$ & $\{0,1,\ldots,100\}$\\
            \hline
        \end{tabular}
        \captionof{table}{Ranges of Parameters.}
        \label{rangetable}
    \end{minipage}
    \hfill 
    \begin{minipage}[c]{0.45\textwidth}
        \centering
        \includegraphics[height=4.5cm, keepaspectratio]{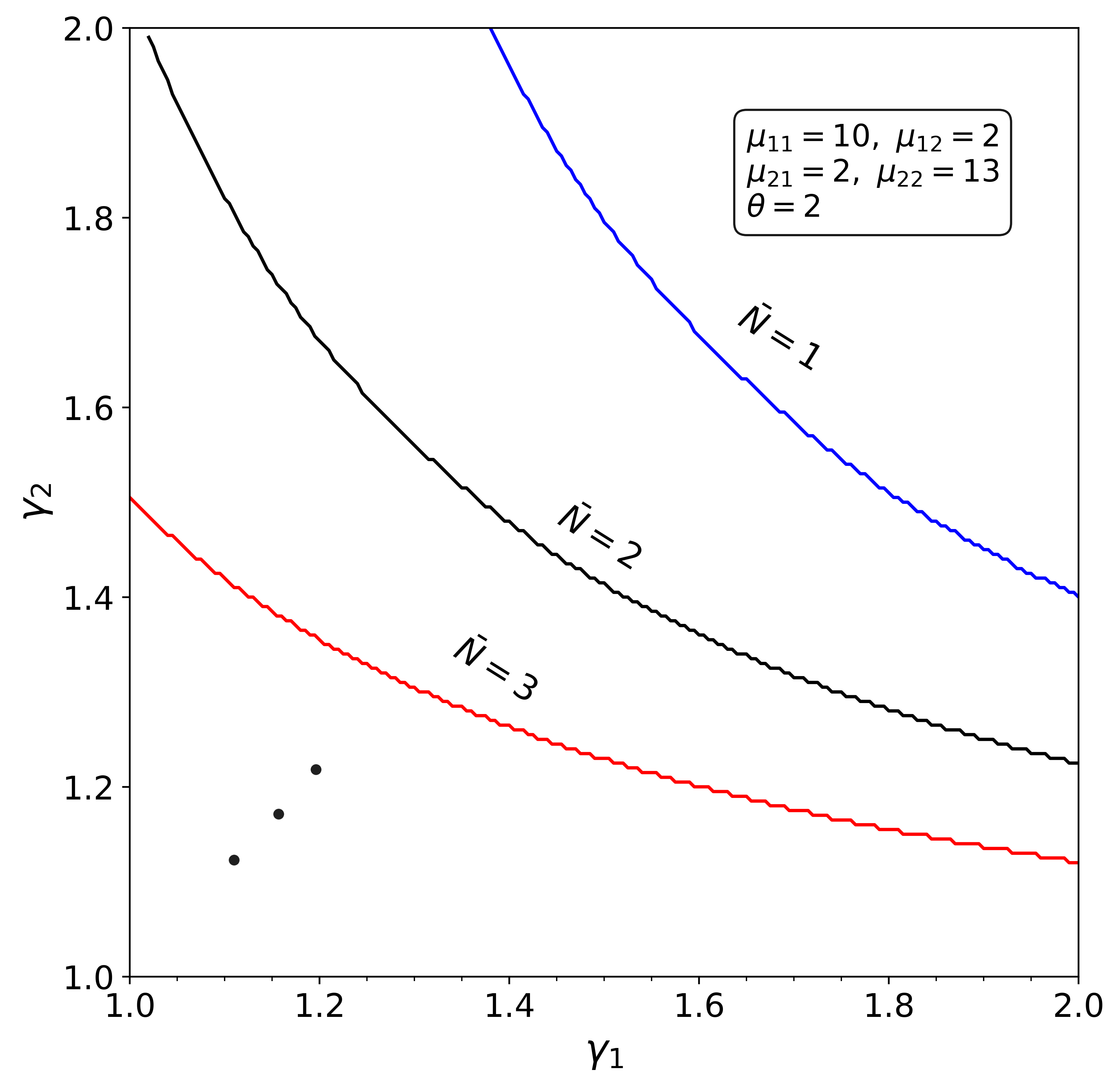}
        \caption{Domains of the optimal buffer size $\bar N$ as functions of $\gamma_1$ and $\gamma_2$.}
        \label{fig1}
    \end{minipage}
\end{figure}

In this section, we present numerical results that suggest that the server assignment policy $d_{\bar N}$ appears to be optimal with task dependent synergy.
Let $h_{\bar N}$ be the vector solving
\begin{align*}
    r_{d_{\bar N}}-\bar{g}_{\bar{N}}e+(\bar{P}_{d_{\bar{N}}}-I)h_{\bar N}=0
\end{align*}
with $h_{\bar N}(0)=0$. As in Section \ref{sec3}, in order to show that $d_{\bar N}$ is the optimal policy using policy iteration for unichain models, we set the initial decision rule $\nu_0$ equal to $d_{\bar N}$ and then show that the next decision rule computed by the policy iteration algorithm is equal to $d_{\bar
N}$. That is we show that 
\begin{align*}
    \nu_1(s) \in \arg \max_{a \in A} \bigg\{r(s,a)+\sum_{j \in S} p(j|s,a)h_{\bar N}(j)\bigg\}
\end{align*}
is equal to $d_{\bar N}(s)$ $ \forall s \in S$. We show this by showing the difference 
\begin{align}\label{difference TD}
    \bar{\Gamma}(s,a)=r(s,d_{\bar N}(s))+\sum_{j \in S} p(j|s,d_{\bar N}(s))h_{\bar N}(j)-r(s,a)-\sum_{j \in S} p(j|s,a)h_{ \bar N}(j)
\end{align}
is non negative, $\forall s \in S$ and $a \in A$.
The difference \eqref{difference TD} was calculated for all states $s$ and actions $a$ for 10,000 replications where the parameters have been chosen independently and uniformly at  random from the ranges given in Table \ref{rangetable}. 
All the numerically computed values of  $\bar \Gamma$ are non-negative for all $s\in S$ and $a \in A$ and for all 10,000 replications 
suggesting that $d_{\bar N}$ is the optimal policy. 

 The frequency of the observed $\bar N$ values for the parameter ranges in Table \ref{rangetable} is given in Supplementary Material \ref{SMTable}. We see that for more than 77.02\% of the cases, the expedite policy is optimal and the frequency of the $\bar N$ value decreases for higher $\bar N$ values, just as in Table \ref{Nfrequency same gamma}. 

 Figure \ref{fig1} shows the domains of $\bar N$ as  functions of $\gamma_1$ and $\gamma_2$ for parameter values $\mu_{11}=10, \mu_{12}=2, \mu_{21}=2, \mu_{22}=13$ and $\theta=2$. It partitions the plane $\gamma_1 \geq 1,  \gamma_2 \geq 1$ into regions based on the value of $\bar N$. The figure suggests that as the values of the synergistic factors $\gamma_1$ and $\gamma_2$ increase, the buffer size used by the proposed optimal policy decreases, as expected from Proposition \ref{thm lamba_1< lambda_2}. This extends  Theorem \eqref{thm2S2S}, which characterizes the behavior along the line $\gamma_1=\gamma_2$, to the region $\gamma_1,\gamma_2 \geq 1$.
\section{Conclusion} \label{Section Conclusion}
In this paper, we considered tandem queueing lines with flexible synergistic servers, finite buffers, and impatient customers. We showed that for systems where each server is equally skilled at all tasks, the expedite policy achieves the maximum long-run average throughput.  For two-station, two-server Markovian systems, we show that the optimal policy does not necessarily use the full buffer. Instead, there exists a threshold depending on the synergistic factor, the abandonment rate, and the service rates of the servers, such that when the number of customers waiting to be processed or being processed at station 2 reaches the threshold the optimal policy switches from having the servers working at their preferred stations to both of them working at station 2. We further showed that as the synergistic factor and abandonment rate increase, the buffer space used by the optimal policy decreases and ultimately the expedite policy becomes the optimal policy. 
Finally, we investigate systems where the synergy can be task dependent, propose an optimal policy, and numerically verify its optimality.
\section*{Acknowledgement}
This work was supported by NSF under grant 2127778. The second author was also supported by NSF under grant 2348409.
\bibliographystyle{apacite}
\bibliography{BibData}
\section*{Appendix}
\begin{proof}[Proof of Lemma \ref{lemma h_N1}]
From \eqref{eq7}, we have for $s=0$,
\begin{align}\label{eq8}
    &0-g_N-\frac{\Sigma'_1}{q}h_N(0)+\frac{\Sigma'_1}{q}h_N(1)=0.
\end{align}
Next, if $N\geq2$, then for $s=1,\ldots N-1$,
\begin{align}\label{eq10}
    &\mu_{22}-g_N+\frac{\mu_{22}+(s-1)\theta}{q}h_N(s-1)-\frac{\mu_{11}+\mu_{22}+(s-1)\theta}{q}h_N(s)+\frac{\mu_{11}}{q}h_N(s+1)=0.
\end{align}
And for $s=N,\ldots, B+2$,
\begin{align}\label{eq11}
   & \Sigma'_2-g_N+\frac{\Sigma'_2+(s-1)\theta}{q}h_N(s-1)-\frac{\Sigma'_2+(s-1)\theta}{q}h_N(s)=0.
\end{align}
From \eqref{eq8} and $h_N(0)=0$, we have 
\begin{align}\label{h_N(1)}
    h_N(1)=\frac{g_Nq}{\Sigma'_1}.
\end{align}
Therefore, \eqref{lemma h_N1(i)} and \eqref{lemma h_N1(ii)} hold for $s=1 \leq N$.  
Next, if $N \geq 2$, we prove by induction that \eqref{lemma h_N1(ii)} holds for $s=2,\ldots, N$.
Plugging $s=1$ into \eqref{eq10} yields 
\begin{align*}
    \frac{\mu_{11}}{q}[h_N(2)-h_N(1)]-\frac{\mu_{22}}{q}h_N(1)&=g_N-\mu_{22}.
\end{align*}
That is 
\begin{align*}
    h_N(2)-h_N(1)&=q\bigg[\frac{g_N-\mu_{22}}{\mu_{11}}+\frac{\mu_{22}g_N}{\Sigma'_1\mu_{11}} \bigg]=q\bigg[\frac{g_N-\mu_{22}}{\mu_{11}}+\frac{g_Nf(1,2)}{\Sigma'_1\mu_{11}} \bigg].
\end{align*}
This proves \eqref{lemma h_N1(ii)} for  $s=2 \leq N$. Now assume that \eqref{lemma h_N1(ii)} holds for $k=2, \ldots,s$ and that $s+1 \leq N$. We show that this implies \eqref{lemma h_N1(ii)} for $k=s+1$. From \eqref{eq10} 
\begin{align*}
    \frac{\mu_{11}}{q}[h_N(s+1)-h_N(s)]=g_N-\mu_{22}+\frac{\mu_{22}+(s-1)\theta}{q}[h_N(s)-h_N(s-1)].
\end{align*}
Using the induction hypothesis, we obtain
\begin{align*}
    h_N(s+1)-h_N(s)
    &=\frac{q}{\mu_{11}}\bigg[g_N-\mu_{22}+[\mu_{22}+(s-1)\theta]\bigg[\sum_{i=0}^{s-2} \frac{g_N-\mu_{22}}{\mu_{11}^{i+1}}f(s-i,s)+\frac{g_N}{\Sigma'_1\mu_{11}^{s-1}}f(1,s) \bigg]\bigg]\\
     &=q\bigg[\frac{g_N-\mu_{22}}{\mu_{11}}+\bigg[\sum_{i=0}^{s-2} \frac{g_N-\mu_{22}}{\mu_{11}^{i+2}}f(s-i,s+1)+\frac{g_N}{\Sigma'_1\mu_{11}^{s}}f(1,s+1) \bigg]\bigg].
\end{align*}
Letting $i'=i+1$, we have
\begin{align*}
     h_N(s+1)-h_N(s)
    &=q\bigg[\frac{g_N-\mu_{22}}{\mu_{11}}+\sum_{i'=1}^{s-1} \frac{g_N-\mu_{22}}{\mu_{11}^{i'+1}}f(s-i'+1,s+1)+\frac{g_N}{\Sigma'_1\mu_{11}^{s}}f(1,s+1) \bigg]\\
     &=q\bigg[\sum_{i'=0}^{s-1} \frac{g_N-\mu_{22}}{\mu_{11}^{i'+1}}f(s-i'+1,s+1)+\frac{g_N}{\Sigma'_1\mu_{11}^{s}}f(1,s+1) \bigg].
\end{align*}
Thus \eqref{lemma h_N1(ii)} holds for $s+1$. This completes the proof of \eqref{lemma h_N1(ii)} for $s=1, \ldots, N$.

Finally from \eqref{eq11}, we have 
for $s=N,\ldots, B+2$,
\begin{align*}
    h_N(s)-h_N(s-1)=\frac{q[\Sigma'_2-g_N]}{\Sigma'_2+(s-1)\theta}.
\end{align*}
\end{proof}
\begin{proof}[Proof of Lemma \ref{lemma h_N2}]
For $s=1$, we have using \eqref{g_n}, \eqref{beta(n)}, \eqref{h_N(1)}, and $h_N(0)=0$ that
\begin{align*}
    h_N(1)-h_N(0)=\frac{g_Nq}{\Sigma'_1}=\frac{\beta(N)q}{\Delta(N)\Sigma'_1}=\frac{q}{\Delta(N)}\bigg[[\Sigma'_2+(N-1)\theta]\mu_{22}\sum_{k=2}^{N} \mu_{11}^{k-2}f(k,N)+ \Sigma'_2\mu_{11}^{N-1} \bigg].
\end{align*}
Hence \eqref{lemma h_N2(i)} holds for $s=1 \leq N$. 

For $s=2,\ldots, N$, we have using \eqref{g_n} that
\begin{align*}
    g_N-\mu_{22}=&\frac{[\Sigma'_2+(N-1)\theta]\Sigma'_1\mu_{22}\alpha(N)+\Sigma'_1 \Sigma'_2\mu_{11}^{N-1}}{[\Sigma'_2+(N-1)\theta][f(1,N)+\Sigma'_1\alpha(N)]+\Sigma'_1\mu_{11}^{N-1}}-\mu_{22}\\
    =&\frac{-[\Sigma'_2+(N-1)\theta]\mu_{22}f(1,N)+\Sigma'_1 \mu_{11}^{N-1}[\Sigma'_2-\mu_{22}]}{\Delta(N)}
\end{align*}
and 
\begin{align*}
    \frac{g_N}{\Sigma'_1}
    =&\frac{[\Sigma'_2+(N-1)\theta]\mu_{22}\sum_{k=2}^{N} \mu_{11}^{k-2}f(k,N)+\Sigma'_2\mu_{11}^{N-1}}{\Delta(N)}.
\end{align*}
Therefore the first term  in \eqref{lemma h_N1(ii)} satisfies
\begin{align*}
    &\sum_{i=0}^{s-2} \frac{g_N-\mu_{22}}{\mu_{11}^{i+1}}f(s-i,s)\\
   =& \frac{\Sigma'_1\mu_{11}^{N-1}[\Sigma'_2-\mu_{22}]}{\Delta(N)}\sum_{i=0}^{s-2} \frac{1}{\mu_{11}^{i+1}}f(s-i,s)-\frac{[\Sigma'_2+(N-1)\theta]\mu_{22}f(1,N)}{\Delta(N)}\sum_{i=0}^{s-2} \frac{1}{\mu_{11}^{i+1}}f(s-i,s).
\end{align*}
The first term in the above satisfies
\begin{align*}
   & \frac{\Sigma'_1\mu_{11}^{N-1}[\Sigma'_2-\mu_{22}]}{\Delta(N)}\sum_{i=0}^{s-2} \frac{1}{\mu_{11}^{i+1}}f(s-i,s)=\frac{\Sigma'_1\mu_{11}^{N-s}[\Sigma'_2-\mu_{22}]}{\Delta(N)}\sum_{i=0}^{s-2}\mu_{11}^{s-2-i} f(s-i,s).
\end{align*}
Letting $s-i=i'$ yields
\begin{align*}
  \frac{\Sigma'_1\mu_{11}^{N-1}[\Sigma'_2-\mu_{22}]}{\Delta(N)}\sum_{i=0}^{s-2} \frac{1}{\mu_{11}^{i+1}}f(s-i,s) 
    =& \frac{\Sigma'_1\mu_{11}^{N-s}[\Sigma'_2-\mu_{22}]}{\Delta(N)}\alpha(s).
\end{align*}
Therefore,
\begin{align}\label{eq13}
     &\sum_{i=0}^{s-2} \frac{g_N-\mu_{22}}{\mu_{11}^{i+1}}f(s-i,s)
     = \frac{\Sigma'_1\mu_{11}^{N-s}[\Sigma'_2-\mu_{22}]}{\Delta(N)}\alpha(s)-\frac{[\Sigma'_2+(N-1)\theta]\mu_{22}f(1,N)}{\Delta(N)}\sum_{i=0}^{s-2} \frac{1}{\mu_{11}^{i+1}}f(s-i,s).
\end{align}
From \eqref{g_n}- \eqref{beta(n)}, the second term in \eqref{lemma h_N1(ii)} satisfies
\begin{align*}
    \frac{g_N}{\Sigma'_1\mu_{11}^{s-1}}f(1,s) 
    =&\frac{[\Sigma'_2+(N-1)\theta]\mu_{22}}{\Delta(N)}\sum_{k=2}^s \mu_{11}^{k-1-s}f(k,N)f(1,s)\\
    &+\frac{[\Sigma'_2+(N-1)\theta]\mu_{22}\sum_{k=s+1}^{N} \mu_{11}^{k-2}f(k,N)+\Sigma'_2\mu_{11}^{N-1}}{\Delta(N)\mu_{11}^{s-1}}f(1,s).
\end{align*}
For the first term in the above, letting $k'=s-k$ yields
\begin{align*}
    \frac{[\Sigma'_2+(N-1)\theta]\mu_{22}}{\Delta(N)}\sum_{k=2}^s \mu_{11}^{k-1-s}f(k,N)f(1,s)=&\frac{[\Sigma'_2+(N-1)\theta]\mu_{22}}{\Delta(N)}\sum_{k'=0}^{s-2} \frac{1}{\mu_{11}^{k'+1}}f(s-k',N)f(1,s)\\
    =&\frac{[\Sigma'_2+(N-1)\theta]\mu_{22}f(1,N)}{\Delta(N)}\sum_{k'=0}^{s-2} \frac{1}{\mu_{11}^{k'+1}}f(s-k',s),
\end{align*}
where in the last equality we have used $f(s-k',N)f(1,s)=f(1,N)f(s-k',s)$. Hence 
\begin{align}\label{eq14}
     \frac{g_N}{\Sigma'_1\mu_{11}^{s-1}}f(1,s) 
     =&\frac{[\Sigma'_2+(N-1)\theta]\mu_{22}f(1,N)}{\Delta(N)}\sum_{k'=0}^{s-2} \frac{1}{\mu_{11}^{k'+1}}f(s-k',s)\nonumber\\
     &+\frac{[\Sigma'_2+(N-1)\theta]\mu_{22}\sum_{k=s+1}^{N} \mu_{11}^{k-s-1}f(k,N)+\Sigma'_2\mu_{11}^{N-s}}{\Delta(N)}f(1,s).
\end{align}
Putting \eqref{eq13} and \eqref{eq14} in \eqref{lemma h_N1(ii)} we get for $s=2,\ldots, N$,
\begin{align*}
    &h_N(s)-h_N(s-1)\\
    =&q \bigg[ \frac{\Sigma'_1\mu_{11}^{N-s}[\Sigma'_2-\mu_{22}]}{\Delta(N)}\alpha(s)-\frac{[\Sigma'_2+(N-1)\theta]\mu_{22}f(1,N)}{\Delta(N)}\sum_{i=0}^{s-2} \frac{1}{\mu_{11}^{i+1}}f(s-i,s)\\
   &+\frac{[\Sigma'_2+(N-1)\theta]\mu_{22}f(1,N)}{\Delta(N)}\sum_{k'=0}^{s-2} \frac{1}{\mu_{11}^{k'+1}}f(s-k',s)\nonumber\\
     &+\frac{[\Sigma'_2+(N-1)\theta]\mu_{22}\sum_{k=s+1}^{N} \mu_{11}^{k-s-1}f(k,N)+\Sigma'_2\mu_{11}^{N-s}}{\Delta(N)}f(1,s) \bigg]\\
     =&\frac{q}{\Delta(N)} \bigg[ \Sigma'_1\mu_{11}^{N-s}[\Sigma'_2-\mu_{22}]\alpha(s)+[\Sigma'_2+(N-1)\theta]\mu_{22}\sum_{k=s+1}^{N} \mu_{11}^{k-s-1}f(k,N)f(1,s)+\Sigma'_2\mu_{11}^{N-s}f(1,s)\bigg ].
\end{align*}
We have shown that \eqref{lemma h_N2(i)} holds for $s=1, \ldots, N$.

Finally for $s=N, \ldots,B+2$ we have from \eqref{lemma h_N1(iii)} that
\begin{align*}
    h_N(s)-h_N(s-1)=\frac{q[\Sigma'_2-g_N]}{\Sigma'_2+(s-1)\theta}.
\end{align*}
Now, from \eqref{g_n}-\eqref{Delta(n)}, we have
\begin{align*}
    \Sigma'_2-g_N=\Sigma'_2-\frac{\beta(N)}{\Delta(N)}
    =&\Sigma'_2-\frac{[\Sigma'_2+(N-1)\theta]\Sigma'_1\mu_{22}\alpha(N)+\Sigma'_1 \Sigma'_2\mu_{11}^{N-1}}{[\Sigma'_2+(N-1)\theta][f(1,N)+\Sigma'_1\alpha(N)]+\Sigma'_1\mu_{11}^{N-1}}\\
    =&\frac{[\Sigma'_2+(N-1)\theta]\bigg(\Sigma'_2f(1,N)+\Sigma'_1[\Sigma'_2-\mu_{22}]\alpha(N)\bigg)}{\Delta(N)}.
\end{align*}
Hence,
\begin{align*}
      h_N(s)-h_N(s-1)=\frac{q[\Sigma'_2+(N-1)\theta]\bigg(\Sigma'_2f(1,N)+\Sigma'_1[\Sigma'_2-\mu_{22}]\alpha(N)\bigg)}{\Delta(N)[\Sigma'_2+(s-1)\theta]}.
\end{align*}
\end{proof}
\begin{proof}[Proof of Lemma \ref{Lemma alpha}]
We have to show that for $s=1, \ldots, N-1$, $\sum_{k=2}^N\mu_{11}^{k-2}f(k,N) \geq \mu_{11}^{N-s}\sum_{k=2}^s \mu_{11}^{k-2}f(k,s).$
That is 
\begin{align*}
   & \sum_{k=2}^N\mu_{11}^{k}f(k,N) \geq \sum_{k=2}^s \mu_{11}^{N-s+k}f(k,s)
    \iff   \sum_{k=N-s+2}^N\mu_{11}^{k}f(k,N)-\sum_{k=2}^s \mu_{11}^{N-s+k}f(k,s)+ \sum_{k=2}^{N-s+1}\mu_{11}^{k}f(k,N) \geq 0.
\end{align*}
Setting $i=N-k$ in the first term and $i=s-k$ in the second term yields
\begin{align*}
    & \sum_{k=2}^N\mu_{11}^{k}f(k,N) \geq \sum_{k=2}^s \mu_{11}^{N-s+k}f(k,s)  \iff \sum_{i=0}^{s-2}\mu_{11}^{N-i}\bigg[f(N-i,N)-f(s-i,s) \bigg]+\sum_{k=2}^{N-s+1}\mu_{11}^{k}f(k,N) \geq 0.
\end{align*}
Therefore it is sufficient to show that $f(N-i,N) \geq f(s-i,s)$,    which is true since for $s=1, \ldots, N-1$, we have
\begin{align*}
    f(N-i,N)=\prod_{j=N-i}^{N-1}[\mu_{22}+(j-1)\theta] \geq \prod_{j=s-i}^{s-1}[\mu_{22}+(j-1)\theta]=f(s-i,s).
\end{align*}
\end{proof}
\setcounter{section}{0} 
\renewcommand{\thesection}{A} 
\newpage
\section{Supplementary Materials}
\subsection{Details of the Proof of Proposition 4}\label{SMProp4}
Eq. \eqref{tau'(n+1)} states that for $n \geq 2$,
\begin{align*}
    \tau'(n+1)=&[\Sigma'_2+n\theta][\Sigma'_2+(n-1)\theta]\mu_{22}f(1,n)-[\Sigma'_2+n\theta][\Sigma'_1\alpha(n+1)(\Sigma'_2-\mu_{22})+\Sigma'_2f(1,n+1)] \nonumber\\
     & \quad \quad \quad +[\Sigma'_2+(n-1)\theta]\mu_{11}[\Sigma'_2f(1,n)+\Sigma'_1\alpha(n)(\Sigma'_2-\mu_{22})].
\end{align*}
The second term satisfies
\begin{align}\label{tau'(n+1)eq1}
   & [\Sigma'_2+n\theta][\Sigma'_1\alpha(n+1)(\Sigma'_2-\mu_{22})+\Sigma'_2f(1,n+1)] \nonumber\\
   =& [\Sigma'_2+n\theta][\Sigma'_1(\Sigma'_2-\mu_{22})(\mu_{22}+(n-1)\theta)\alpha(n)+\Sigma'_1 (\Sigma'_2-\mu_{22})\mu_{11}^{n-1}+\Sigma'_2f(1,n)[\mu_{22}+(n-1   )\theta]] \nonumber\\
   =&[\mu_{22}+(n-1)\theta][\Sigma'_2+n\theta][\Sigma'_1(\Sigma'_2-\mu_{22})\alpha(n)+\Sigma'_2f(1,n)]+\Sigma'_1\mu_{11}^{n-1}(\Sigma'_2-\mu_{22})[\Sigma'_2+n\theta] \nonumber\\
   =&[\mu_{22}+(n-1)\theta][\Sigma'_2+(n-1)\theta][\Sigma'_1(\Sigma'_2-\mu_{22})\alpha(n)+\Sigma'_2f(1,n)]+\Sigma'_1\mu_{11}^{n-1}(\Sigma'_2-\mu_{22})[\Sigma'_2+n\theta] \nonumber\\
   &+\theta[\mu_{22}+(n-1)\theta][\Sigma'_1 (\Sigma'_2-\mu_{22})\alpha(n)+\Sigma'_2f(1,n)].
\end{align}
The third term satisfies
\begin{align}\label{tau'(n+1)eq2}
    &[\Sigma'_2+(n-1)\theta]\mu_{11}[\Sigma'_2f(1,n)+\Sigma'_1\alpha(n)(\Sigma'_2-\mu_{22})] \nonumber\\
    =&[\Sigma'_2+(n-2)\theta]\mu_{11}[\Sigma'_2f(1,n)+\Sigma'_1\alpha(n)(\Sigma'_2-\mu_{22})]+\theta \mu_{11}[\Sigma'_2f(1,n)+\Sigma'_1\alpha(n)(\Sigma'_2-\mu_{22})] \nonumber\\
    =&[\Sigma'_2+(n-2)\theta]\mu_{11}[\Sigma'_2f(1,n-1)[\mu_{22}+(n-2)\theta]+\Sigma'_1 (\Sigma'_2-\mu_{22})[\mu_{22}+(n-2)\theta]\alpha(n-1) \nonumber\\
    &+\Sigma'_1 (\Sigma'_2-\mu_{22})\mu_{11}^{n-2}]
    +\theta \mu_{11}[\Sigma'_2f(1,n)+\Sigma'_1\alpha(n)(\Sigma'_2-\mu_{22})] \nonumber\\
    =&[\Sigma'_2+(n-2)\theta]\mu_{11}[\Sigma'_2f(1,n-1)[\mu_{22}+(n-2)\theta]+\Sigma'_1 (\Sigma'_2-\mu_{22})[\mu_{22}+(n-2)\theta]\alpha(n-1)] \nonumber\\
    &+\Sigma'_1 (\Sigma'_2-\mu_{22})\mu_{11}^{n-2}[\Sigma'_2+(n-2)\theta]\mu_{11}+\theta \mu_{11}[\Sigma'_2f(1,n)+\Sigma'_1\alpha(n)(\Sigma'_2-\mu_{22})] \nonumber\\
     =&[\Sigma'_2+(n-2)\theta]\mu_{11}[\mu_{22}+(n-2)\theta][\Sigma'_2f(1,n-1)+\Sigma'_1 (\Sigma'_2-\mu_{22})\alpha(n-1)] \nonumber\\
    &+\Sigma'_1 (\Sigma'_2-\mu_{22})\mu_{11}^{n-1}[\Sigma'_2+(n-2)\theta]+\theta \mu_{11}[\Sigma'_2f(1,n)+\Sigma'_1\alpha(n)(\Sigma'_2-\mu_{22})] \nonumber\\
     =&[\Sigma'_2+(n-2)\theta]\mu_{11}[\mu_{22}+(n-1)\theta][\Sigma'_2f(1,n-1)+\Sigma'_1 (\Sigma'_2-\mu_{22})\alpha(n-1)] \nonumber\\
    &+\Sigma'_1 (\Sigma'_2-\mu_{22})\mu_{11}^{n-1}[\Sigma'_2+(n-2)\theta]+\theta \mu_{11}[\Sigma'_2f(1,n)+\Sigma'_1\alpha(n)(\Sigma'_2-\mu_{22})] \nonumber\\
    &-\theta [\Sigma'_2+(n-2)\theta]\mu_{11}[\Sigma'_2f(1,n-1)+\Sigma'_1 (\Sigma'_2-\mu_{22})\alpha(n-1)].
\end{align}
Finally, the first term satisfies
\begin{align}\label{tau'(n+1)eq3}
    &[\Sigma'_2+n\theta][\Sigma'_2+(n-1)\theta]\mu_{22}f(1,n) \nonumber\\
    =&[\Sigma'_2+(n-2)\theta][\Sigma'_2+(n-1)\theta]\mu_{22}f(1,n)+2\theta[\Sigma'_2+(n-1)\theta]\mu_{22}f(1,n) \nonumber\\
    =&[\Sigma'_2+(n-2)\theta][\Sigma'_2+(n-1)\theta]\mu_{22}f(1,n-1)[\mu_{22}+(n-2)\theta]+2\theta[\Sigma'_2+(n-1)\theta]\mu_{22}f(1,n) \nonumber\\
    =&[\Sigma'_2+(n-2)\theta][\Sigma'_2+(n-1)\theta]\mu_{22}f(1,n-1)[\mu_{22}+(n-1)\theta]+2\theta[\Sigma'_2+(n-1)\theta]\mu_{22}f(1,n) \nonumber\\
    &-\theta [\Sigma'_2+(n-2)\theta][\Sigma'_2+(n-1)\theta]\mu_{22}f(1,n-1).
\end{align}
Combining the first terms in  \eqref{tau'(n+1)eq1}-\eqref{tau'(n+1)eq3} yields $[\mu_{22}+(n-1)\theta]\tau'(n)$. Therefore, writing
\begin{align*}
    \tau'(n+1)=[\mu_{22}+(n-1)\theta]\tau'(n)-\epsilon'(n)
\end{align*}
leads to 
\begin{align*}
&\epsilon'(n)\\
   = &-2\theta[\Sigma'_2+(n-1)\theta]\mu_{22}f(1,n)+\theta [\Sigma'_2+(n-2)\theta][\Sigma'_2+(n-1)\theta]\mu_{22}f(1,n-1) \\
    &+\Sigma'_1\mu_{11}^{n-1}(\Sigma'_2-\mu_{22})[\Sigma'_2+n\theta]+\theta[\mu_{22}+(n-1)\theta][\Sigma'_1 (\Sigma'_2-\mu_{22})\alpha(n)+\Sigma'_2f(1,n)]\\
    &-\Sigma'_1 (\Sigma'_2-\mu_{22})\mu_{11}^{n-1}[\Sigma'_2+(n-2)\theta]-\theta \mu_{11}[\Sigma'_2f(1,n)+\Sigma'_1\alpha(n)(\Sigma'_2-\mu_{22})]\\
    &+\theta [\Sigma'_2+(n-2)\theta]\mu_{11}[\Sigma'_2f(1,n-1)+\Sigma'_1 (\Sigma'_2-\mu_{22})\alpha(n-1)]\\    
    =&-2\theta[\Sigma'_2+(n-1)\theta]\mu_{22}f(1,n)+\theta [\Sigma'_2+(n-2)\theta][\Sigma'_2+(n-1)\theta]\mu_{22}f(1,n-1) \\
    &+2\theta\Sigma'_1\mu_{11}^{n-1}(\Sigma'_2-\mu_{22})+\theta[\mu_{22}+(n-1)\theta][\Sigma'_1 (\Sigma'_2-\mu_{22})\alpha(n)+\Sigma'_2f(1,n)]\\
    &-\theta \mu_{11}[\Sigma'_2f(1,n)+\Sigma'_1\alpha(n)(\Sigma'_2-\mu_{22})]\\
    &+\theta [\Sigma'_2+(n-2)\theta]\mu_{11}[\Sigma'_2f(1,n-1)+\Sigma'_1 (\Sigma'_2-\mu_{22})\alpha(n-1)]\\
    =&-2\theta[\Sigma'_2+(n-1)\theta]\mu_{22}f(1,n)+\theta [\Sigma'_2+(n-2)\theta][\Sigma'_2+(n-1)\theta]\mu_{22}f(1,n-1) \\
    &+2\theta\Sigma'_1\mu_{11}^{n-1}(\Sigma'_2-\mu_{22})+\theta[\mu_{22}+(n-1)\theta][\Sigma'_1 (\Sigma'_2-\mu_{22})\alpha(n)+\Sigma'_2f(1,n)]\\
    &-\theta \mu_{11}[\mu_{22}+(n-2)\theta][\Sigma'_2f(1,n-1)+\Sigma'_1\alpha(n-1)(\Sigma'_2-\mu_{22})]-\theta \mu_{11}^{n-1}\Sigma'_1 (\Sigma'_2-\mu_{22})\\
    &+\theta [\Sigma'_2+(n-2)\theta]\mu_{11}[\Sigma'_2f(1,n-1)+\Sigma'_1 (\Sigma'_2-\mu_{22})\alpha(n-1)]\\
    =&-2\theta[\Sigma'_2+(n-1)\theta]\mu_{22}f(1,n-1)[\mu_{22}+(n-2)\theta]+\theta [\Sigma'_2+(n-2)\theta][\Sigma'_2+(n-1)\theta]\mu_{22}f(1,n-1) \\
    &+\theta\Sigma'_1\mu_{11}^{n-1}(\Sigma'_2-\mu_{22})+\theta[\mu_{22}+(n-1)\theta][\Sigma'_1 (\Sigma'_2-\mu_{22})\alpha(n)+\Sigma'_2f(1,n)]\\
    &+\theta (\Sigma'_2-\mu_{22})\mu_{11}[\Sigma'_2f(1,n-1)+\Sigma'_1 (\Sigma'_2-\mu_{22})\alpha(n-1)]\\
    =&-\theta[\Sigma'_2+(n-1)\theta]\mu_{22}f(1,n-1)[\mu_{22}+(n-2)\theta]+\theta (\Sigma'_2-\mu_{22})[\Sigma'_2+(n-1)\theta]\mu_{22}f(1,n-1) \\
    &+\theta\Sigma'_1\mu_{11}^{n-1}(\Sigma'_2-\mu_{22})+\theta[\mu_{22}+(n-1)\theta][\Sigma'_1 (\Sigma'_2-\mu_{22})\alpha(n)+\Sigma'_2f(1,n)]\\
    &+\theta (\Sigma'_2-\mu_{22})\mu_{11}[\Sigma'_2f(1,n-1)+\Sigma'_1 (\Sigma'_2-\mu_{22})\alpha(n-1)]\\
     =&-\theta[\Sigma'_2+(n-1)\theta]\mu_{22}f(1,n-1)[\mu_{22}+(n-2)\theta]+\theta (\Sigma'_2-\mu_{22})[\Sigma'_2+(n-1)\theta]\mu_{22}f(1,n-1) \\
    &+\theta\Sigma'_1\mu_{11}^{n-1}(\Sigma'_2-\mu_{22})+\theta[\mu_{22}+(n-1)\theta][\mu_{22}+(n-2)\theta][\Sigma'_1 (\Sigma'_2-\mu_{22})\alpha(n-1)]\\
    &+\theta[\mu_{22}+(n-1)\theta][\mu_{22}+(n-2)\theta]\Sigma'_2f(1,n-1) +\Sigma'_1(\Sigma'_2-\mu_{22})\mu_{11}^{n-2}\theta[\mu_{22}+(n-1)\theta]\\
    &+\theta (\Sigma'_2-\mu_{22})\mu_{11}[\Sigma'_2f(1,n-1)+\Sigma'_1 (\Sigma'_2-\mu_{22})\alpha(n-1)]\\
    =&-\theta[\Sigma'_2+(n-1)\theta]\mu_{22}f(1,n-1)[\mu_{22}+(n-2)\theta]+\theta (\Sigma'_2-\mu_{22})[\Sigma'_2+(n-1)\theta]\mu_{22}f(1,n-1) \\
    &+\theta\Sigma'_1\mu_{11}^{n-1}(\Sigma'_2-\mu_{22})+\theta[\mu_{22}+(n-1)\theta][\mu_{22}+(n-2)\theta]\Sigma'_1 (\Sigma'_2-\mu_{22})\alpha(n-1)\\
    &+\theta[\mu_{22}+(n-1)\theta][\mu_{22}+(n-2)\theta](\Sigma'_2-\mu_{22})f(1,n-1)\\
    &+\theta[\mu_{22}+(n-1)\theta][\mu_{22}+(n-2)\theta]\mu_{22}f(1,n-1) +\Sigma'_1(\Sigma'_2-\mu_{22})\mu_{11}^{n-2}\theta[\mu_{22}+(n-1)\theta]\\
    &+\theta (\Sigma'_2-\mu_{22})\mu_{11}[\Sigma'_2f(1,n-1)+\Sigma'_1 (\Sigma'_2-\mu_{22})\alpha(n-1)]\\
     =&-\theta(\Sigma'_2-\mu_{22})\mu_{22}f(1,n-1)[\mu_{22}+(n-2)\theta]+\theta (\Sigma'_2-\mu_{22})[\Sigma'_2+(n-1)\theta]\mu_{22}f(1,n-1) \nonumber\\
    &+\theta\Sigma'_1\mu_{11}^{n-1}(\Sigma'_2-\mu_{22})+\theta[\mu_{22}+(n-1)\theta][\mu_{22}+(n-2)\theta]\Sigma'_1 (\Sigma'_2-\mu_{22})\alpha(n-1) \nonumber\\
    &+\theta[\mu_{22}+(n-1)\theta][\mu_{22}+(n-2)\theta](\Sigma'_2-\mu_{22})f(1,n-1) \nonumber\\
    &    +\Sigma'_1(\Sigma'_2-\mu_{22})\mu_{11}^{n-2}\theta[\mu_{22}+(n-1)\theta] \nonumber\\
    &+\theta (\Sigma'_2-\mu_{22})\mu_{11}[\Sigma'_2f(1,n-1)+\Sigma'_1 (\Sigma'_2-\mu_{22})\alpha(n-1)] \nonumber\\
    =&\theta (\Sigma'_2-\mu_{22})[\Sigma'_2-\mu_{22}+\theta]\mu_{22}f(1,n-1) \nonumber \\
    &+\theta\Sigma'_1\mu_{11}^{n-1}(\Sigma'_2-\mu_{22})+\theta[\mu_{22}+(n-1)\theta][\mu_{22}+(n-2)\theta]\Sigma'_1 (\Sigma'_2-\mu_{22})\alpha(n-1) \nonumber\\
    &+\theta[\mu_{22}+(n-1)\theta][\mu_{22}+(n-2)\theta](\Sigma'_2-\mu_{22})f(1,n-1) \nonumber\\
    &    +\Sigma'_1(\Sigma'_2-\mu_{22})\mu_{11}^{n-2}\theta[\mu_{22}+(n-1)\theta] \nonumber\\
    &+\theta (\Sigma'_2-\mu_{22})\mu_{11}[\Sigma'_2f(1,n-1)+\Sigma'_1 (\Sigma'_2-\mu_{22})\alpha(n-1)]. \nonumber
\end{align*}
\subsection{Proof of Theorem \ref{thm2S2S}}\label{SMTheorem1}
\begin{proof}
As mentioned at the end of Section \ref{sec3}, we employ policy iteration for unichain models, as detailed in \cite{Puterman}. By initializing the decision rule as $\nu_0 = d_N$, we demonstrate that the subsequent decision rule generated by the algorithm, remains equal to $d_N$. Specifically, we prove that$$\nu_1(s) \in \arg \max_{a \in A} \bigg\{r(s,a)+\sum_{j \in S} p(j|s,a)h_N(j)\bigg\}$$yields $d_N(s)$ for all $s \in S$. To achieve this, we establish that the difference function$$\Gamma(s,a)=r(s,d_N(s))+\sum_{j \in S} p(j|s,d_N(s))h_N(j)-r(s,a)-\sum_{j \in S} p(j|s,a)h_N(j)$$is non-negative for all $s \in S$ and $a \in A$. The conclusion then follows directly from Theorem 8.6.6 in \cite{Puterman}. We use policy iteration for unichain models (see \cite{Puterman}). We set the initial decision rule $\nu_0$ equal to $d_N$ and then show that the next decision rule computed by the policy iteration algorithm is equal to $d_N$. That is we show that 
\begin{align*}
    \nu_1(s) \in \arg \max_{a \in A} \bigg\{r(s,a)+\sum_{j \in S} p(j|s,a)h_N(j)\bigg\} 
\end{align*}
is equal to $d_N(s)$, $\forall s \in S$. We show this by showing the difference 
\begin{align*}
    \Gamma(s,a)=r(s,d_N(s))+\sum_{j \in S} p(j|s,d_N(s))h_N(j)-r(s,a)-\sum_{j \in S} p(j|s,a)h_N(j)
\end{align*}
is non negative for all $s \in S$ and $a \in A$. The result then follows from Theorem 8.6.6 in \cite{Puterman} .

We consider three separate cases; Case 1: $s=0$; Case 2: $s=1, \ldots, N-1$; and Case 3: $s=N, \ldots, B+2$.

\item[Case 1:] $s=0$.

In this case we have $d_N(s)=a_{11}$.
   
    \item[ Case 1, Part (i)] :  $a=a_{12}$.
     
     From Lemma \ref{lemma h_N1}, we have
    \begin{align*}
        \Gamma(0,a_{12})=&\frac{\Sigma'_1}{q}h_N(1)+\left( 1-\frac{\Sigma'_1}{q}\right)h_N(0)-\frac{\mu_{11}}{q}h_N(1)-\left( 1-\frac{\mu_{11}}{q}\right)h_N(0)\\
        =&\frac{h_N(1)}{q}[\Sigma'_1-\mu_{11}]=\frac{g_N}{\Sigma'_1}[\Sigma'_1-\mu_{11}] \geq 0.
    \end{align*}
    
    
    \item[Case 1, Part (ii) :] $a=a_{21}$.
    
    From Lemma \ref{lemma h_N1}, we have
    \begin{align*}
        \Gamma(0,a_{21})=&\frac{\Sigma'_1}{q}h_N(1)+\left( 1-\frac{\Sigma'_1}{q}\right)h_N(0)-\frac{\mu_{21}}{q}h_N(1)-\left( 1-\frac{\mu_{21}}{q}\right)h_N(0)\\
        =&\frac{h_N(1)}{q}[\Sigma'_1-\mu_{21}]=\frac{g_N}{\Sigma'_1}[\Sigma'_1-\mu_{21}]
        \geq 0.
    \end{align*}
   \item[Case 1, Part (iii) :]  $a=a_{22}$.

  From Lemma \ref{lemma h_N1}, we have
    \begin{align*}
        \Gamma(0,a_{22})=&\frac{\Sigma'_1}{q}h_N(1)+\left( 1-\frac{\Sigma'_1}{q}\right)h_N(0)=\frac{\Sigma'_1}{q}\frac{g_Nq}{\Sigma'_1}=g_N \geq0.
    \end{align*}
    
    Thus we have proved $\Gamma(s,a) \geq 0$ for all policies $a$ in Case 1.
    
\item[Case 2:] $s=1,\ldots, N-1$.
    
   In this case we have $d_N(s)=a_{12}$. 
    \item[Case 2, Part (i) :] $a=a_{22}$.
    \begin{align}\label{Gamma(s,a_22)}
        &\Gamma(s,a_{22}) \nonumber\\
        =&\mu_{22}+\frac{\mu_{22}+(s-1)\theta}{q}h_N(s-1)+\left(1-\frac{\mu_{11}+\mu_{22}+(s-1)\theta}{q} \right)h_N(s)+\frac{\mu_{11}}{q}h_N(s+1) \nonumber\\
        &-\Sigma'_2-\frac{\Sigma'_2+(s-1)\theta}{q}h_N(s-1)-\left(1-\frac{\Sigma'_2+(s-1)\theta}{q} \right)h_N(s) \nonumber\\
        =&\mu_{22}-\Sigma'_2+\frac{\mu_{22}}{q}[h_N(s-1)-h_N(s)]-\frac{\mu_{11}}{q}[h_N(s)-h_N(s+1)]-\frac{\Sigma'_2}{q}[h_N(s-1)-h_N(s)] \nonumber\\
        =&\mu_{22}-\Sigma'_2+\frac{\Sigma'_2-\mu_{22}}{q}[h_N(s)-h_N(s-1)]+\frac{\mu_{11}}{q}[h_N(s+1)-h_N(s)].
        \end{align}
       We have, from \eqref{eq10} in proof of Lemma \ref{lemma h_N1}, that
       \begin{align}\label{Case 2 bias difference}
           \frac{\mu_{11}}{q}[h_N(s+1)-h_N(s)]=g_N-\mu_{22}+\frac{\mu_{22}+(s-1)\theta }{q}[h_N(s)-h_N(s-1)].
       \end{align} 
       Therefore
        \begin{align*}
            &\Gamma(s,a_{22})\\
            =&\mu_{22}-\Sigma'_2+\frac{\Sigma'_2-\mu_{22}}{q}[h_N(s)-h_N(s-1)]+g_N-\mu_{22}+\frac{\mu_{22}+(s-1)\theta }{q}[h_N(s)-h_N(s-1)]\\
            =&g_N-\Sigma'_2+\frac{\Sigma'_2+(s-1)\theta}{q}[h_N(s)-h_N(s-1)]. 
        \end{align*}
        From Lemma \ref{lemma h_N1} we have 
        \begin{align*}
            &\Gamma(s,a_{22})\\
            =&g_N-\Sigma'_2+[\Sigma'_2+(s-1)\theta]\bigg[\sum_{i=0}^{s-2} \frac{g_N-\mu_{22}}{\mu_{11}^{i+1}}f(s-i,s)+\frac{g_N}{\Sigma'_1\mu_{11}^{s-1}}f(1,s) \bigg]\\
            =&g_N \bigg(1+\frac{\Sigma'_2+(s-1)\theta}{\Sigma'_1 \mu_{11}^{s-1}}f(1,s)+[\Sigma'_2+(s-1)\theta]\sum_{i=0}^{s-2}\frac{1}{\mu_{11}^{i+1}}f(s-i,s) \bigg)\\
            &-\bigg(\Sigma'_2+\mu_{22}[\Sigma'_2+(s-1)\theta]\sum_{i=0}^{s-2}\frac{1}{\mu_{11}^{i+1}}f(s-i,s)  \bigg)\\
            =&g_N \bigg(\frac{\Sigma'_1 \mu_{11}^{s-1}+[\Sigma'_2+(s-1)\theta]f(1,s)+\Sigma'_1[\Sigma'_2+(s-1)\theta]\sum_{i=0}^{s-2}\mu_{11}^{s-i-2}f(s-i,s)}{\Sigma'_1 \mu_{11}^{s-1}} \bigg)\\
            &-\bigg(\frac{\Sigma'_1 \Sigma'_2 \mu_{11}^{s-1}+\mu_{22}\Sigma'_1[\Sigma'_2+(s-1)\theta]\sum_{i=0}^{s-2}\mu_{11}^{s-i-2}f(s-i,s)}{\Sigma'_1 \mu_{11}^{s-1}} \bigg).
        \end{align*}
        Let $s-i=i'$  and $j'=j+1$. We have
        \begin{align*}
         &\Gamma(s,a_{22})\\
            =&g_N \bigg(\frac{\Sigma'_1 \mu_{11}^{s-1}+[\Sigma'_2+(s-1)\theta]f(1,s)+\Sigma'_1[\Sigma'_2+(s-1)\theta]\sum_{i'=2}^{s}\mu_{11}^{i'-2}f(i',s)}{\Sigma'_1 \mu_{11}^{s-1}} \bigg)\\
            &-\bigg(\frac{\Sigma'_1 \Sigma'_2 \mu_{11}^{s-1}+\mu_{22}\Sigma'_1[\Sigma'_2+(s-1)\theta]\sum_{i'=2}^{s}\mu_{11}^{i'-2}f(i',s)}{\Sigma'_1 \mu_{11}^{s-1}} \bigg)\\
             =&g_N \bigg(\frac{\Sigma'_1 \mu_{11}^{s-1}+[\Sigma'_2+(s-1)\theta]f(1,s)+\Sigma'_1[\Sigma'_2+(s-1)\theta]\alpha(s)}{\Sigma'_1 \mu_{11}^{s-1}} \bigg)\\
            &-\bigg(\frac{\Sigma'_1 \Sigma'_2 \mu_{11}^{s-1}+\mu_{22}\Sigma'_1[\Sigma'_2+(s-1)\theta]\alpha(s)}{\Sigma'_1 \mu_{11}^{s-1}} \bigg)\\
           \stackrel{(*)}{\geq}& 
            \frac{\beta(s+1)}{\Delta(s+1)} \bigg(\frac{\Sigma'_1 \mu_{11}^{s-1}+[\Sigma'_2+(s-1)\theta][f(1,s)+\Sigma'_1\alpha(s)]}{\Sigma'_1 \mu_{11}^{s-1}} \bigg)\\
            &-\Delta(s+1)\bigg(\frac{\Sigma'_1 \Sigma'_2 \mu_{11}^{s-1}+\mu_{22}\Sigma'_1[\Sigma'_2+(s-1)\theta]\alpha(s)}{\Delta(s+1)\Sigma'_1 \mu_{11}^{s-1}} \bigg),
        \end{align*}
        where $(*)$ follows from $g_N \geq g_{s+1}$ (by the definition of $N$) and \eqref{g_n}.
From \eqref{beta(n)}-\eqref{Delta(n)} we have
\begin{align*}
  \Gamma(s,a_{22})  \geq&\frac{\beta(s+1)\Delta(s)-\Delta(s+1)\beta(s)}{\Delta(s+1)\Sigma'_1 \mu_{11}^{s-1}}=\frac{\Sigma'_1\mu_{11}^{s-1} \tau'(s+1)}{\Delta(s+1)\Sigma'_1 \mu_{11}^{s-1}}\geq 0
\end{align*}
from \eqref{g difference} and the definition of $N$.

\item[Case 2, Part (ii) :] $a=a_{21}$.

We have
\begin{align*}
    \Gamma(s,a_{21})=&\mu_{22}+\frac{\mu_{22}+(s-1)\theta}{q}h_N(s-1)+\left(1-\frac{\mu_{11}+\mu_{22}+(s-1)\theta}{q} \right)h_N(s)\\
        &+\frac{\mu_{11}}{q}h_N(s+1)-\mu_{12}-\frac{\mu_{12}+(s-1)\theta}{q}h_N(s-1)-\frac{\mu_{21}}{q}h_N(s+1)\\
        &-\left(1-\frac{\mu_{21}+\mu_{12}+(s-1)\theta}{q} \right)h_N(s)\\
        =&\mu_{22}-\mu_{12}+\frac{\mu_{22}}{q}[h_N(s-1)-h_N(s)]+\frac{\mu_{11}}{q}[h_N(s+1)-h_N(s)]\\
        &+\frac{\mu_{12}}{q}[h_N(s)-h_N(s-1)]+\frac{\mu_{21}}{q}[h_N(s)-h_N(s+1)]\\
        =&\mu_{22}-\mu_{12}+\frac{\mu_{12}-\mu_{22}}{q}[h_N(s)-h_N(s-1))]+\frac{\mu_{11}-\mu_{21}}{q}[h_N(s+1)-h_N(s)]\\
         =&\mu_{22}-\Sigma'_2+\frac{\Sigma'_2-\mu_{22}}{q}[h_N(s)-h_N(s-1))]+\frac{\mu_{11}}{q}[h_N(s+1)-h_N(s)]\\
         &+\Sigma'_2-\mu_{12}+\frac{\Sigma'_2-\mu_{12}}{q}[h_N(s-1)-h_N(s)]+\frac{\mu_{21}}{q}[h_N(s)-h_N(s+1)]\\
         =&\Gamma(s,a_{22})+\Gamma'(s,a_{21})
\end{align*}
from \eqref{Gamma(s,a_22)}, where
\begin{align*}
    \Gamma'(s,a_{21})=\Sigma'_2-\mu_{12}-\frac{\Sigma'_2-\mu_{12}}{q}[h_N(s)-h_N(s-1)]-\frac{\mu_{21}}{q}[h_N(s+1)-h_N(s)].
\end{align*}
Since we have shown that  $\Gamma(s,a_{22})$ is non-negative, it is sufficient to show $\Gamma'(s,a_{21})$ is non-negative.

Using \eqref{Case 2 bias difference}, we have
\begin{align*}
    &\Gamma'(s,a_{21})\\
    =&\Sigma'_2-\mu_{12}-\frac{\Sigma'_2-\mu_{12}}{q}[h_N(s)-h_N(s-1)]-\frac{g_N\mu_{21}-\mu_{22}\mu_{21}}{\mu_{11}}\\
    &-\frac{\mu_{21}\mu_{22}+\mu_{21}(s-1)\theta }{q\mu_{11}}  \times [h_N(s)-h_N(s-1)]\\
    \stackrel{(*)}{=}&\frac{\mu_{11}(\Sigma'_2-\mu_{12})+\mu_{22}\mu_{21}-\mu_{21}g_N}{\mu_{11}}+\frac{-\mu_{11}(\Sigma'_2-\mu_{12})-\mu_{21}\mu_{22}-\mu_{21}(s-1)\theta}{\Delta(N)\mu_{11}} \\
    &\times \bigg[\Sigma'_1\mu_{11}^{N-s}[\Sigma'_2-\mu_{22}]\alpha(s)+[\Sigma'_2+(N-1)\theta]\mu_{22}\sum_{k=s+1}^{N} \mu_{11}^{k-s-1}f(k,N)f(1,s)+\Sigma'_2\mu_{11}^{N-s}f(1,s) \bigg]
\end{align*}
where $(*)$ follows from Lemma \ref{lemma h_N2}. This yields
\begin{align*}
    & \Gamma'(s,a_{21})\Delta(N)\mu_{11}\\
    =&[\mu_{11}(\Sigma'_2-\mu_{12})+\mu_{22}\mu_{21}]\Delta(N)-\mu_{21}\beta(N)+[-\mu_{11}(\Sigma'_2-\mu_{12})-\mu_{21}\mu_{22}-\mu_{21}(s-1)\theta]\\
    &\times \bigg[\Sigma'_1\mu_{11}^{N-s}[\Sigma'_2-\mu_{22}]\alpha(s)+[\Sigma'_2+(N-1)\theta]\mu_{22}\sum_{k=s+1}^{N} \mu_{11}^{k-s-1}f(k,N)f(1,s)+\Sigma'_2\mu_{11}^{N-s}f(1,s) \bigg]\\
    =&[\mu_{11}(\Sigma'_2-\mu_{12})+\mu_{22}\mu_{21}][\Sigma'_2+(N-1)\theta][f(1,N)+\Sigma'_1\alpha(N)]+[\mu_{11}(\Sigma'_2-\mu_{12})+\mu_{22}\mu_{21}]\Sigma'_1\mu_{11}^{N-1}\\
    &-\mu_{21}[\Sigma'_2+(N-1)\theta]\Sigma'_1\mu_{22}\alpha(N)-\mu_{21}\Sigma'_1 \Sigma'_2\mu_{11}^{N-1}\\
    &+[-\mu_{11}(\Sigma'_2-\mu_{12})-\mu_{21}\mu_{22}-\mu_{21}(s-1)\theta]\Sigma'_1\mu_{11}^{N-s}[\Sigma'_2-\mu_{22}]\alpha(s)\\
    &+[-\mu_{11}(\Sigma'_2-\mu_{12})-\mu_{21}\mu_{22}-\mu_{21}(s-1)\theta][\Sigma'_2+(N-1)\theta]\mu_{22}\sum_{k=s+1}^{N} \mu_{11}^{k-s-1}f(k,N)f(1,s)\\
    &+[-\mu_{11}(\Sigma'_2-\mu_{12})-\mu_{21}\mu_{22}-\mu_{21}(s-1)\theta]\Sigma'_2\mu_{11}^{N-s}f(1,s)\\
    =&[\mu_{11}(\Sigma'_2-\mu_{12})+\mu_{22}\mu_{21}]\Sigma'_1\mu_{11}^{N-1}-\mu_{21}\Sigma'_1 \Sigma'_2\mu_{11}^{N-1}\\
    &+[-\mu_{11}(\Sigma'_2-\mu_{12})-\mu_{21}\mu_{22}-\mu_{21}(s-1)\theta]\Sigma'_1\mu_{11}^{N-s}[\Sigma'_2-\mu_{22}]\alpha(s)\\
    &+[\mu_{11}(\Sigma'_2-\mu_{12})+\mu_{22}\mu_{21}][\Sigma'_2+(N-1)\theta][f(1,N)+\Sigma'_1\alpha(N)]-\mu_{21}[\Sigma'_2+(N-1)\theta]\Sigma'_1\mu_{22}\alpha(N)\\
    &+[-\mu_{11}(\Sigma'_2-\mu_{12})-\mu_{21}\mu_{22}-\mu_{21}(s-1)\theta][\Sigma'_2+(N-1)\theta]\mu_{22}\sum_{k=s+1}^{N} \mu_{11}^{k-s-1}f(k,N)f(1,s)\\
    &+[-\mu_{11}(\Sigma'_2-\mu_{12})-\mu_{21}\mu_{22}-\mu_{21}(s-1)\theta]\Sigma'_2\mu_{11}^{N-s}f(1,s)\\
    =&[\mu_{11}(\Sigma'_2-\mu_{12})+\mu_{22}\mu_{21}-\mu_{21}\Sigma'_2]\Sigma'_1\mu_{11}^{N-1}-[\mu_{11}(\Sigma'_2-\mu_{12})+\mu_{21}\mu_{22}]\Sigma'_1\mu_{11}^{N-s}[\Sigma'_2-\mu_{22}]\alpha(s)\\
    &+[\mu_{11}(\Sigma'_2-\mu_{12})+\mu_{22}\mu_{21}-\mu_{21}\Sigma'_2](s-1)\theta\Sigma'_1\mu_{11}^{N-s}\alpha(s)-(\mu_{11}(\Sigma'_2-\mu_{12}))(s-1)\theta \Sigma'_1 \mu_{11}^{N-s} \alpha(s)\\
     &+[\mu_{11}(\Sigma'_2-\mu_{12})+\mu_{22}\mu_{21}][\Sigma'_2+(N-1)\theta][f(1,N)+\Sigma'_1\alpha(N)]-\mu_{21}[\Sigma'_2+(N-1)\theta]\Sigma'_1\mu_{22}\alpha(N)\\
    &+[-\mu_{11}(\Sigma'_2-\mu_{12})-\mu_{21}\mu_{22}-\mu_{21}(s-1)\theta][\Sigma'_2+(N-1)\theta]\mu_{22}\sum_{k=s+1}^{N} \mu_{11}^{k-s-1}f(k,N)f(1,s)\\
    &+[-\mu_{11}(\Sigma'_2-\mu_{12})-\mu_{21}\mu_{22}-\mu_{21}(s-1)\theta]\Sigma'_2\mu_{11}^{N-s}f(1,s)\\
    =&[\mu_{11}(\Sigma'_2-\mu_{12})+\mu_{22}\mu_{21}-\mu_{21}\Sigma'_2]\bigg(\Sigma'_1\mu_{11}^{N-1}+(s-1)\theta\Sigma'_1\mu_{11}^{N-s}\alpha(s)\bigg)\\
    &-[\mu_{11}(\Sigma'_2-\mu_{12})+\mu_{21}\mu_{22}]\Sigma'_1\mu_{11}^{N-s}[\Sigma'_2-\mu_{22}]\alpha(s)\\
    &-(\mu_{11}(\Sigma'_2-\mu_{12}))(s-1)\theta \Sigma'_1 \mu_{11}^{N-s} \alpha(s)\\
     &+[\mu_{11}(\Sigma'_2-\mu_{12})+\mu_{22}\mu_{21}][\Sigma'_2+(N-1)\theta][f(1,N)+\Sigma'_1\alpha(N)]-\mu_{21}[\Sigma'_2+(N-1)\theta]\Sigma'_1\mu_{22}\alpha(N)\\
    &+[-\mu_{11}(\Sigma'_2-\mu_{12})-\mu_{21}\mu_{22}-\mu_{21}(s-1)\theta][\Sigma'_2+(N-1)\theta]\mu_{22}\sum_{k=s+1}^{N} \mu_{11}^{k-s-1}f(k,N)f(1,s)\\
    &+[-\mu_{11}(\Sigma'_2-\mu_{12})-\mu_{21}\mu_{22}-\mu_{21}(s-1)\theta]\Sigma'_2\mu_{11}^{N-s}f(1,s)\\
    =&[\mu_{11}(\Sigma'_2-\mu_{12})+\mu_{22}\mu_{21}-\mu_{21}\Sigma'_2]\bigg(\Sigma'_1\mu_{11}^{N-1}+(s-1)\theta\Sigma'_1\mu_{11}^{N-s}\alpha(s)\bigg)\\
    &-[\mu_{11}(\Sigma'_2-\mu_{12})+\mu_{21}\mu_{22}]\Sigma'_1\mu_{11}^{N-s}[\Sigma'_2-\mu_{22}]\alpha(s)\\
    &-(\mu_{11}(\Sigma'_2-\mu_{12}))(s-1)\theta \Sigma'_1 \mu_{11}^{N-s} \alpha(s)\\
    &+(-\mu_{11}(\Sigma'_2-\mu_{12})-\mu_{21}\mu_{22}-\Sigma'_1(s-1)\theta)[\Sigma'_2+(N-1)\theta]\mu_{22}\sum_{k=s+1}^{N} \mu_{11}^{k-s-1}f(k,N)f(1,s)\\
    &+(\Sigma'_1-\mu_{21})(s-1)\theta[\Sigma'_2+(N-1)\theta]\mu_{22}\sum_{k=s+1}^{N} \mu_{11}^{k-s-1}f(k,N)f(1,s)\\
    &+[\mu_{11}(\Sigma'_2-\mu_{12})+\mu_{22}\mu_{21}][\Sigma'_2+(N-1)\theta][f(1,N)+\Sigma'_1\alpha(N)]-\mu_{21}[\Sigma'_2+(N-1)\theta]\Sigma'_1\mu_{22}\alpha(N)\\
    &+[-\mu_{11}(\Sigma'_2-\mu_{12})-\mu_{21}\mu_{22}-\mu_{21}(s-1)\theta]\Sigma'_2\mu_{11}^{N-s}f(1,s)\\  
    =&[\mu_{11}(\Sigma'_2-\mu_{12})+\mu_{22}\mu_{21}-\mu_{21}\Sigma'_2]\bigg(\Sigma'_1\mu_{11}^{N-1}+(s-1)\theta\Sigma'_1\mu_{11}^{N-s}\alpha(s)\bigg)\\
    &-[\mu_{11}(\Sigma'_2-\mu_{12})+\mu_{21}\mu_{22}]\Sigma'_1\mu_{11}^{N-s}[\Sigma'_2-\mu_{22}]\alpha(s)\\
    &-(\mu_{11}(\Sigma'_2-\mu_{12}))(s-1)\theta \Sigma'_1 \mu_{11}^{N-s} \alpha(s)\\
    &+(-\mu_{11}(\Sigma'_2-\mu_{12})-\mu_{21}\mu_{22}-\Sigma'_1(s-1)\theta)[\Sigma'_2+(N-1)\theta]\mu_{22}\sum_{k=s+1}^{N} \mu_{11}^{k-s-1}f(k,N)f(1,s)\\
    &+(\Sigma'_1-\mu_{21})(s-1)\theta[\Sigma'_2+(N-1)\theta]\mu_{22}\sum_{k=s+1}^{N} \mu_{11}^{k-s-1}f(k,N)f(1,s)\\
    &+[\mu_{11}(\Sigma'_2-\mu_{12})+\mu_{22}\mu_{21}][\Sigma'_2+(N-1)\theta][f(1,N)+\mu_{11}\alpha(N)]\\
    &+[\mu_{11}(\Sigma'_2-\mu_{12})+\mu_{22}\mu_{21}][\Sigma'_2+(N-1)\theta][\lambda \mu_{11}+(\lambda+1)\mu_{21}]\alpha(N)-\mu_{21}[\Sigma'_2+(N-1)\theta]\Sigma'_1\mu_{22}\alpha(N)\\
    &+[-\mu_{11}(\Sigma'_2-\mu_{12})-\mu_{21}\mu_{22}-\mu_{21}(s-1)\theta]\Sigma'_2\mu_{11}^{N-s}f(1,s)\\  
    =&[\mu_{11}(\Sigma'_2-\mu_{12})+\mu_{22}\mu_{21}-\mu_{21}\Sigma'_2]\bigg(\Sigma'_1\mu_{11}^{N-1}+(s-1)\theta\Sigma'_1\mu_{11}^{N-s}\alpha(s)\bigg)\\
    &-[\mu_{11}(\Sigma'_2-\mu_{12})+\mu_{21}\mu_{22}]\Sigma'_1\mu_{11}^{N-s}[\Sigma'_2-\mu_{22}]\alpha(s)\\
    &-(\mu_{11}(\Sigma'_2-\mu_{12}))(s-1)\theta \Sigma'_1 \mu_{11}^{N-s} \alpha(s)\\
    &+(-\mu_{11}(\Sigma'_2-\mu_{12})-\mu_{21}\mu_{22}-\Sigma'_1(s-1)\theta)[\Sigma'_2+(N-1)\theta]\mu_{22}\sum_{k=s+1}^{N} \mu_{11}^{k-s-1}f(k,N)f(1,s)\\
    &+(\Sigma'_1-\mu_{21})(s-1)\theta[\Sigma'_2+(N-1)\theta]\mu_{22}\sum_{k=s+1}^{N} \mu_{11}^{k-s-1}f(k,N)f(1,s)\\
    &+[\mu_{11}(\Sigma'_2-\mu_{12})+\mu_{22}\mu_{21}][\Sigma'_2+(N-1)\theta][f(1,N)]\\
    &+[\mu_{11}(\Sigma'_2-\mu_{12})+\mu_{22}\mu_{21}][\Sigma'_2+(N-1)\theta]\bigg[\sum_{k=2}^{N+1-s}\mu_{11}^{k-1}f(k,N)+\sum_{k=N+2-s}^N\mu_{11}^{k-1}f(k,N) \bigg] \\
    &+[\Sigma'_2+(N-1)\theta][[\mu_{11}(\Sigma'_2-\mu_{12})+\mu_{22}\mu_{21}][\lambda \mu_{11}+(\lambda+1)\mu_{21}]-\mu_{21}\Sigma'_1\mu_{22}]\alpha(N)\\
    &+[-\mu_{11}(\Sigma'_2-\mu_{12})-\mu_{21}\mu_{22}-\mu_{21}(s-1)\theta]\Sigma'_2\mu_{11}^{N-s}f(1,s)\\  
     =&[\mu_{11}(\Sigma'_2-\mu_{12})+\mu_{22}\mu_{21}-\mu_{21}\Sigma'_2]\bigg(\Sigma'_1\mu_{11}^{N-1}+(s-1)\theta\Sigma'_1\mu_{11}^{N-s}\alpha(s)\bigg)\\
    &-[\mu_{11}(\Sigma'_2-\mu_{12})+\mu_{21}\mu_{22}]\Sigma'_1\mu_{11}^{N-s}[\Sigma'_2-\mu_{22}]\alpha(s)\\
    &-(\mu_{11}(\Sigma'_2-\mu_{12}))(s-1)\theta \Sigma'_1 \mu_{11}^{N-s} \alpha(s)\\
    &+(-\mu_{11}(\Sigma'_2-\mu_{12})-\mu_{21}\mu_{22}-\Sigma'_1(s-1)\theta)[\Sigma'_2+(N-1)\theta]\mu_{22}\sum_{k=s+1}^{N} \mu_{11}^{k-s-1}f(k,N)f(1,s)\\
    &+(\Sigma'_1-\mu_{21})(s-1)\theta[\Sigma'_2+(N-1)\theta]\mu_{22}\sum_{k=s+1}^{N} \mu_{11}^{k-s-1}f(k,N)f(1,s)\\
    &+[\mu_{11}(\Sigma'_2-\mu_{12})+\mu_{22}\mu_{21}][\Sigma'_2+(N-1)\theta][f(1,N)]\\
    &+[\mu_{11}(\Sigma'_2-\mu_{12})+\mu_{22}\mu_{21}][\Sigma'_2+(N-1)\theta]\bigg[\sum_{k=2}^{N+1-s}\mu_{11}^{k-1}f(k,N)\bigg] \\
    &+[\mu_{11}(\Sigma'_2-\mu_{12})+\mu_{22}\mu_{21}][\Sigma'_2+(N-1)\theta]\bigg[\sum_{k=N+2-s}^N\mu_{11}^{k-1}f(k,N) \bigg] \\
    &+[\Sigma'_2+(N-1)\theta][[\mu_{11}(\Sigma'_2-\mu_{12})+\mu_{22}\mu_{21}][\lambda \mu_{11}+(\lambda+1)\mu_{21}]-\mu_{21}\Sigma'_1\mu_{22}]\alpha(N)\\
    &+[-\mu_{11}(\Sigma'_2-\mu_{12})-\mu_{21}\mu_{22}-\mu_{21}(s-1)\theta]\Sigma'_2\mu_{11}^{N-s}f(1,s)\\  
\end{align*}
Let $k'+N-s=k$
\begin{align*}
     =&[\mu_{11}(\Sigma'_2-\mu_{12})+\mu_{22}\mu_{21}-\mu_{21}\Sigma'_2]\bigg(\Sigma'_1\mu_{11}^{N-1}+(s-1)\theta\Sigma'_1\mu_{11}^{N-s}\alpha(s)\bigg)\\
    &-[\mu_{11}(\Sigma'_2-\mu_{12})+\mu_{21}\mu_{22}]\Sigma'_1\mu_{11}^{N-s}[\Sigma'_2-\mu_{22}]\alpha(s)\\
    &-(\mu_{11}(\Sigma'_2-\mu_{12}))(s-1)\theta \Sigma'_1 \mu_{11}^{N-s} \alpha(s)\\
    &+(-\mu_{11}(\Sigma'_2-\mu_{12})-\mu_{21}\mu_{22}-\Sigma'_1(s-1)\theta)[\Sigma'_2+(N-1)\theta]\mu_{22}\sum_{k=s+1}^{N} \mu_{11}^{k-s-1}f(k,N)f(1,s)\\
    &+(\Sigma'_1-\mu_{21})(s-1)\theta[\Sigma'_2+(N-1)\theta]\mu_{22}\sum_{k=s+1}^{N} \mu_{11}^{k-s-1}f(k,N)f(1,s)\\
    &+[\mu_{11}(\Sigma'_2-\mu_{12})+\mu_{22}\mu_{21}][\Sigma'_2+(N-1)\theta][f(1,N)]\\
    &+[\mu_{11}(\Sigma'_2-\mu_{12})+\mu_{22}\mu_{21}][\Sigma'_2+(N-1)\theta]\bigg[\sum_{k=2}^{N+1-s}\mu_{11}^{k-1}f(k,N)\bigg] \\
    &+[\mu_{11}(\Sigma'_2-\mu_{12})+\mu_{22}\mu_{21}][\Sigma'_2+(N-1)\theta]\bigg[\sum_{k'=2}^s\mu_{11}^{k'+N-1-s}f(k'+N-s,N) \bigg] \\
    &+[\Sigma'_2+(N-1)\theta][[\mu_{11}(\Sigma'_2-\mu_{12})+\mu_{22}\mu_{21}][\lambda \mu_{11}+(\lambda+1)\mu_{21}]-\mu_{21}\Sigma'_1\mu_{22}]\alpha(N)\\
    &+[-\mu_{11}(\Sigma'_2-\mu_{12})-\mu_{21}\mu_{22}-\mu_{21}(s-1)\theta]\Sigma'_2\mu_{11}^{N-s}f(1,s)\\  
    =&[\mu_{11}(\Sigma'_2-\mu_{12})+\mu_{22}\mu_{21}-\mu_{21}\Sigma'_2]\bigg(\Sigma'_1\mu_{11}^{N-1}+(s-1)\theta\Sigma'_1\mu_{11}^{N-s}\alpha(s)\bigg)\\
    &-[\mu_{11}(\Sigma'_2-\mu_{12})+\mu_{21}\mu_{22}]\Sigma'_1\mu_{11}^{N-s}[\Sigma'_2-\mu_{22}]\alpha(s)\\
    &-(\mu_{11}(\Sigma'_2-\mu_{12}))(s-1)\theta \Sigma'_1 \mu_{11}^{N-s} \alpha(s)\\
    &+(-\mu_{11}(\Sigma'_2-\mu_{12})-\mu_{21}\mu_{22}-\Sigma'_1(s-1)\theta)[\Sigma'_2+(N-1)\theta]\mu_{22}\sum_{k=s+1}^{N} \mu_{11}^{k-s-1}f(k,N)f(1,s)\\
    &+(\Sigma'_1-\mu_{21})(s-1)\theta[\Sigma'_2+(N-1)\theta]\mu_{22}\sum_{k=s+1}^{N} \mu_{11}^{k-s-1}f(k,N)f(1,s)\\
    &+[\mu_{11}(\Sigma'_2-\mu_{12})+\mu_{22}\mu_{21}][\Sigma'_2+(N-1)\theta]\bigg[\sum_{k=2}^{N+1-s}\mu_{11}^{k-1}f(k,N)+f(1,N)\bigg] \\
    &+[\mu_{11}(\Sigma'_2-\mu_{12})+\mu_{22}\mu_{21}][\Sigma'_2+(N-1)\theta]\bigg[\sum_{k'=2}^s\mu_{11}^{k'+N-1-s}f(k'+N-s,N) \bigg] \\
    &+[\Sigma'_2+(N-1)\theta][[\mu_{11}(\Sigma'_2-\mu_{12})+\mu_{22}\mu_{21}][\lambda \mu_{11}+(\lambda+1)\mu_{21}]-\mu_{21}\Sigma'_1\mu_{22}]\alpha(N)\\
    &+[-\mu_{11}(\Sigma'_2-\mu_{12})-\mu_{21}\mu_{22}-\mu_{21}(s-1)\theta]\Sigma'_2\mu_{11}^{N-s}f(1,s)\\  
     =&[\mu_{11}(\Sigma'_2-\mu_{12})+\mu_{22}\mu_{21}-\mu_{21}\Sigma'_2]\bigg(\Sigma'_1\mu_{11}^{N-1}+(s-1)\theta\Sigma'_1\mu_{11}^{N-s}\alpha(s)\bigg)\\
    &-[\mu_{11}(\Sigma'_2-\mu_{12})+\mu_{21}\mu_{22}]\Sigma'_1\mu_{11}^{N-s}[\Sigma'_2-\mu_{22}]\alpha(s)\\
    &-(\mu_{11}(\Sigma'_2-\mu_{12}))(s-1)\theta \Sigma'_1 \mu_{11}^{N-s} \alpha(s)\\
    &+(-\mu_{11}(\Sigma'_2-\mu_{12})-\mu_{21}\mu_{22}-\Sigma'_1(s-1)\theta)[\Sigma'_2+(N-1)\theta]\mu_{22}\sum_{k=s+1}^{N} \mu_{11}^{k-s-1}f(k,N)f(1,s)\\
    &+(\Sigma'_1-\mu_{21})(s-1)\theta[\Sigma'_2+(N-1)\theta]\mu_{22}\sum_{k=s+1}^{N} \mu_{11}^{k-s-1}f(k,N)f(1,s)\\
    &+[\mu_{11}(\Sigma'_2-\mu_{12})+\mu_{22}\mu_{21}][\Sigma'_2+(N-1)\theta]\bigg[\sum_{k=1}^{N+1-s}\mu_{11}^{k-1}f(k,N)\bigg] \\
    &+[\mu_{11}(\Sigma'_2-\mu_{12})+\mu_{22}\mu_{21}][\Sigma'_2+(N-1)\theta]\bigg[\sum_{k'=2}^s\mu_{11}^{k'+N-1-s}f(k'+N-s,N) \bigg] \\
    &+[\Sigma'_2+(N-1)\theta][[\mu_{11}(\Sigma'_2-\mu_{12})+\mu_{22}\mu_{21}][\lambda \mu_{11}+(\lambda+1)\mu_{21}]-\mu_{21}\Sigma'_1\mu_{22}]\alpha(N)\\
    &+[-\mu_{11}(\Sigma'_2-\mu_{12})-\mu_{21}\mu_{22}-\mu_{21}(s-1)\theta]\Sigma'_2\mu_{11}^{N-s}f(1,s)\\  
     =&[\mu_{11}(\Sigma'_2-\mu_{12})+\mu_{22}\mu_{21}-\mu_{21}\Sigma'_2]\bigg(\Sigma'_1\mu_{11}^{N-1}+(s-1)\theta\Sigma'_1\mu_{11}^{N-s}\alpha(s)\bigg)\\
    &-[\mu_{11}(\Sigma'_2-\mu_{12})+\mu_{21}\mu_{22}]\Sigma'_1\mu_{11}^{N-s}[\Sigma'_2-\mu_{22}]\alpha(s)\\
    &-(\mu_{11}(\Sigma'_2-\mu_{12}))(s-1)\theta \Sigma'_1 \mu_{11}^{N-s} \alpha(s)\\
    &+(-\mu_{11}(\Sigma'_2-\mu_{12})-\mu_{21}\mu_{22}-\Sigma'_1(s-1)\theta)[\Sigma'_2+(N-1)\theta]\mu_{22}\sum_{k=s+1}^{N} \mu_{11}^{k-s-1}f(k,N)f(1,s)\\
    &+(\Sigma'_1-\mu_{21})(s-1)\theta[\Sigma'_2+(N-1)\theta]\mu_{22}\sum_{k=s+1}^{N} \mu_{11}^{k-s-1}f(k,N)f(1,s)\\
    &+[\mu_{11}(\Sigma'_2-\mu_{12})+\mu_{22}\mu_{21}][\Sigma'_2+(N-1)\theta]\bigg[\sum_{k=1}^{N-s}\mu_{11}^{k-1}f(k,N)+\mu_{11}^{N-s}f(N+1-s,N)\bigg] \\
    &+[\mu_{11}(\Sigma'_2-\mu_{12})+\mu_{22}\mu_{21}][\Sigma'_2+(N-1)\theta]\bigg[\sum_{k'=2}^s\mu_{11}^{k'+N-1-s}f(k'+N-s,N) \bigg] \\
    &+[\Sigma'_2+(N-1)\theta][[\mu_{11}(\Sigma'_2-\mu_{12})+\mu_{22}\mu_{21}][\lambda \mu_{11}+(\lambda+1)\mu_{21}]-\mu_{21}\Sigma'_1\mu_{22}]\alpha(N)\\
    &+[-\mu_{11}(\Sigma'_2-\mu_{12})-\mu_{21}\mu_{22}-\mu_{21}(s-1)\theta]\Sigma'_2\mu_{11}^{N-s}f(1,s)
\end{align*}
Let $k''-s=k$
\begin{align*}
     =&[\mu_{11}(\Sigma'_2-\mu_{12})+\mu_{22}\mu_{21}-\mu_{21}\Sigma'_2]\bigg(\Sigma'_1\mu_{11}^{N-1}+(s-1)\theta\Sigma'_1\mu_{11}^{N-s}\alpha(s)\bigg)\\
    &-[\mu_{11}(\Sigma'_2-\mu_{12})+\mu_{21}\mu_{22}]\Sigma'_1\mu_{11}^{N-s}[\Sigma'_2-\mu_{22}]\alpha(s)\\
    &-(\mu_{11}(\Sigma'_2-\mu_{12}))(s-1)\theta \Sigma'_1 \mu_{11}^{N-s} \alpha(s)\\
    &+[-\mu_{11}(\Sigma'_2-\mu_{12})-\mu_{21}\mu_{22}-\mu_{21}(s-1)\theta][\Sigma'_2+(N-1)\theta]\mu_{22}\sum_{k=s+1}^{N} \mu_{11}^{k-s-1}f(k,N)f(1,s)\\
    &+[\mu_{11}(\Sigma'_2-\mu_{12})+\mu_{22}\mu_{21}][\Sigma'_2+(N-1)\theta]\bigg[\sum_{k''=s+1}^{N}\mu_{11}^{k''-s-1}f(k''-s,N)+\mu_{11}^{N-s}f(N+1-s,N)\bigg] \\
    &+[\mu_{11}(\Sigma'_2-\mu_{12})+\mu_{22}\mu_{21}][\Sigma'_2+(N-1)\theta]\bigg[\sum_{k'=2}^s\mu_{11}^{k'+N-1-s}f(k'+N-s,N) \bigg] \\
    &+[\Sigma'_2+(N-1)\theta][[\mu_{11}(\Sigma'_2-\mu_{12})+\mu_{22}\mu_{21}][\lambda \mu_{11}+(\lambda+1)\mu_{21}]-\mu_{21}\Sigma'_1\mu_{22}]\alpha(N)\\
    &+[-\mu_{11}(\Sigma'_2-\mu_{12})-\mu_{21}\mu_{22}-\mu_{21}(s-1)\theta]\Sigma'_2\mu_{11}^{N-s}f(1,s)\\  
      =&[\mu_{11}(\Sigma'_2-\mu_{12})+\mu_{22}\mu_{21}-\mu_{21}\Sigma'_2]\bigg(\Sigma'_1\mu_{11}^{N-1}+(s-1)\theta\Sigma'_1\mu_{11}^{N-s}\alpha(s)\bigg)\\
    &-[\mu_{11}(\Sigma'_2-\mu_{12})+\mu_{21}\mu_{22}]\Sigma'_1\mu_{11}^{N-s}[\Sigma'_2-\mu_{22}]\alpha(s)\\
    &-(\mu_{11}(\Sigma'_2-\mu_{12}))(s-1)\theta \Sigma'_1 \mu_{11}^{N-s} \alpha(s)\\
    &+\frac{[-\mu_{11}(\Sigma'_2-\mu_{12})-\mu_{21}\mu_{22}-\mu_{21}(s-1)\theta]}{\mu_{22}+(s-1)\theta}[\Sigma'_2+(N-1)\theta]\mu_{22}\sum_{k=s+1}^{N} \mu_{11}^{k-s-1}f(k,N)f(1,s+1)\\
    &+[\mu_{11}(\Sigma'_2-\mu_{12})+\mu_{22}\mu_{21}][\Sigma'_2+(N-1)\theta]\bigg[\sum_{k''=s+1}^{N}\mu_{11}^{k''-s-1}f(k''-s,N)+\mu_{11}^{N-s}f(N+1-s,N)\bigg] \\
    &+\mu_{22}\frac{[\mu_{11}(\Sigma'_2-\mu_{12})+\mu_{22}\mu_{21}+\mu_{21}(s-1)\theta]}{\mu_{22}+(s-1)\theta}[\Sigma'_2+(N-1)\theta]\sum_{k''=s+1}^{N}\mu_{11}^{k''-s-1}f(k''-s,N)\\
     &-\mu_{22}\frac{[\mu_{11}(\Sigma'_2-\mu_{12})+\mu_{22}\mu_{21}+\mu_{21}(s-1)\theta]}{\mu_{22}+(s-1)\theta}[\Sigma'_2+(N-1)\theta]\sum_{k''=s+1}^{N}\mu_{11}^{k''-s-1}f(k''-s,N)\\
    &+[\mu_{11}(\Sigma'_2-\mu_{12})+\mu_{22}\mu_{21}][\Sigma'_2+(N-1)\theta]\bigg[\sum_{k'=2}^s\mu_{11}^{k'+N-1-s}f(k'+N-s,N) \bigg] \\
    &+[\Sigma'_2+(N-1)\theta][[\mu_{11}(\Sigma'_2-\mu_{12})+\mu_{22}\mu_{21}][\lambda \mu_{11}+(\lambda+1)\mu_{21}]-\mu_{21}\Sigma'_1\mu_{22}]\alpha(N)\\
    &+[-\mu_{11}(\Sigma'_2-\mu_{12})-\mu_{21}\mu_{22}-\mu_{21}(s-1)\theta]\Sigma'_2\mu_{11}^{N-s}f(1,s)\\  
    =&[\mu_{11}(\Sigma'_2-\mu_{12})+\mu_{22}\mu_{21}-\mu_{21}\Sigma'_2]\bigg(\Sigma'_1\mu_{11}^{N-1}+(s-1)\theta\Sigma'_1\mu_{11}^{N-s}\alpha(s)\bigg)\\
    &-[\mu_{11}(\Sigma'_2-\mu_{12})+\mu_{21}\mu_{22}]\Sigma'_1\mu_{11}^{N-s}[\Sigma'_2-\mu_{22}]\alpha(s)\\
    &-(\mu_{11}(\Sigma'_2-\mu_{12}))(s-1)\theta \Sigma'_1 \mu_{11}^{N-s} \alpha(s)\\
    +&\frac{[\mu_{11}(\Sigma'_2-\mu_{12})+\mu_{22}\mu_{21}+\mu_{21}(s-1)\theta]}{\mu_{22}+(s-1)\theta}[\Sigma'_2+(N-1)\theta]\mu_{22}\times\\
    &\bigg[\sum_{k=s+1}^{N}\mu_{11}^{k-s-1}[f(k-s,N)-f(k,N)f(1,s+1)]\bigg]\\
    &+[\mu_{11}(\Sigma'_2-\mu_{12})+\mu_{22}\mu_{21}][\Sigma'_2+(N-1)\theta]\bigg[\mu_{11}^{N-s}f(N+1-s,N)\bigg] \\
     +&\bigg[\mu_{11}(\Sigma'_2-\mu_{12})+\mu_{22}\mu_{21}-\mu_{22}\frac{[\mu_{11}(\Sigma'_2-\mu_{12})+\mu_{22}\mu_{21}+\mu_{21}(s-1)\theta]}{\mu_{22}+(s-1)\theta}\bigg][\Sigma'_2+(N-1)\theta]\times \\
     &\sum_{k''=s+1}^{N}\mu_{11}^{k''-s-1}f(k''-s,N)\\
    &+[\mu_{11}(\Sigma'_2-\mu_{12})+\mu_{22}\mu_{21}][\Sigma'_2+(N-1)\theta]\bigg[\sum_{k'=2}^s\mu_{11}^{k'+N-1-s}f(k'+N-s,N) \bigg] \\
    &+[\Sigma'_2+(N-1)\theta]\bigg[[\mu_{11}(\Sigma'_2-\mu_{12})+\mu_{22}\mu_{21}][\lambda \mu_{11}+(\lambda+1)\mu_{21}]-\mu_{21}\Sigma'_1\mu_{22}\bigg]\alpha(N)\\
    &+[-\mu_{11}(\Sigma'_2-\mu_{12})-\mu_{21}\mu_{22}-\mu_{21}(s-1)\theta]\Sigma'_2\mu_{11}^{N-s}f(1,s)\\  
  =&[\mu_{11}(\Sigma'_2-\mu_{12})+\mu_{22}\mu_{21}-\mu_{21}\Sigma'_2]\bigg(\Sigma'_1\mu_{11}^{N-1}+(s-1)\theta\Sigma'_1\mu_{11}^{N-s}\alpha(s)\bigg)\\
  +&\frac{[\mu_{11}(\Sigma'_2-\mu_{12})+\mu_{22}\mu_{21}+\mu_{21}(s-1)\theta]}{\mu_{22}+(s-1)\theta}[\Sigma'_2+(N-1)\theta]\mu_{22}\times\\
  &\bigg[\sum_{k=s+1}^{N}\mu_{11}^{k-s-1}[f(k-s,N)-f(k,N)f(1,s+1)]\bigg]\\
  +&\bigg[\mu_{11}(\Sigma'_2-\mu_{12})+\mu_{22}\mu_{21}-\mu_{22}\frac{[\mu_{11}(\Sigma'_2-\mu_{12})+\mu_{22}\mu_{21}+\mu_{21}(s-1)\theta]}{\mu_{22}+(s-1)\theta}\bigg][\Sigma'_2+(N-1)\theta]\times \\
     &\sum_{k''=s+1}^{N}\mu_{11}^{k''-s-1}f(k''-s,N)\\
    &+ [\mu_{11}(\Sigma'_2-\mu_{12})+\mu_{22}\mu_{21}]\mu_{11}^{N-s}\bigg[\Sigma'_2\sum_{k'=2}^s\mu_{11}^{k'-1}f(k'+N-s,N)-\Sigma_1'[\Sigma'_2-\mu_{22}]\alpha(s) \bigg]\\
    &-(\mu_{11}(\Sigma'_2-\mu_{12}))(s-1)\theta \Sigma'_1 \mu_{11}^{N-s} \alpha(s)\\
    &+[\mu_{11}(\Sigma'_2-\mu_{12})+\mu_{22}\mu_{21}][\Sigma'_2+(N-1)\theta]\bigg[\mu_{11}^{N-s}f(N+1-s,N)\bigg] \\  
     &+[\mu_{11}(\Sigma'_2-\mu_{12})+\mu_{22}\mu_{21}](N-1)\theta\bigg[\sum_{k'=2}^s\mu_{11}^{k'+N-1-s}f(k'+N-s,N) \bigg] \\
    &+[\Sigma'_2+(N-1)\theta]\bigg[[\mu_{11}(\Sigma'_2-\mu_{12})+\mu_{22}\mu_{21}][\lambda \mu_{11}+(\lambda+1)\mu_{21}]-\mu_{21}\Sigma'_1\mu_{22}\bigg]\alpha(N)\\
    &+[-\mu_{11}(\Sigma'_2-\mu_{12})-\mu_{21}\mu_{22}-\mu_{21}(s-1)\theta]\Sigma'_2\mu_{11}^{N-s}f(1,s)\\  
    =&[\mu_{11}(\Sigma'_2-\mu_{12})+\mu_{22}\mu_{21}-\mu_{21}\Sigma'_2]\Sigma'_1\mu_{11}^{N-1}\\
  +&\frac{[\mu_{11}(\Sigma'_2-\mu_{12})+\mu_{22}\mu_{21}+\mu_{21}(s-1)\theta]}{\mu_{22}+(s-1)\theta}[\Sigma'_2+(N-1)\theta]\mu_{22}\times\\
  &\bigg[\sum_{k=s+1}^{N}\mu_{11}^{k-s-1}[f(k-s,N)-f(k,N)f(1,s+1)]\bigg]\\
  +&\bigg[\mu_{11}(\Sigma'_2-\mu_{12})+\mu_{22}\mu_{21}-\mu_{22}\frac{[\mu_{11}(\Sigma'_2-\mu_{12})+\mu_{22}\mu_{21}+\mu_{21}(s-1)\theta]}{\mu_{22}+(s-1)\theta}\bigg][\Sigma'_2+(N-1)\theta]\times \\
     &\sum_{k''=s+1}^{N}\mu_{11}^{k''-s-1}f(k''-s,N)\\
    &+ [\mu_{11}(\Sigma'_2-\mu_{12})+\mu_{22}\mu_{21}]\mu_{11}^{N-s}\bigg[\Sigma'_2\sum_{k'=2}^s\mu_{11}^{k'-1}f(k'+N-s,N)-\Sigma_1'[\Sigma'_2-\mu_{22}]\alpha(s) \bigg]\\
    &-(\mu_{21}\Sigma_2'-\mu_{22}\mu_{21})(s-1)\theta \Sigma'_1 \mu_{11}^{N-s} \alpha(s)\\
    &+[\mu_{11}(\Sigma'_2-\mu_{12})+\mu_{22}\mu_{21}][\Sigma'_2+(N-1)\theta]\bigg[\mu_{11}^{N-s}f(N+1-s,N)\bigg] \\  
     &+(\mu_{11}(\Sigma'_2-\mu_{12})+\mu_{22}\mu_{21}-\mu_{21}\Sigma_2'+\mu_{22}\mu_{21})(N-1)\theta\bigg[\sum_{k'=2}^s\mu_{11}^{k'+N-1-s}f(k'+N-s,N) \bigg] \\
     &+(\mu_{21}\Sigma'_2-\mu_{22}\mu_{21})(N-1)\theta\bigg[\sum_{k'=2}^s\mu_{11}^{k'+N-1-s}f(k'+N-s,N) \bigg] \\
    &+[\Sigma'_2+(N-1)\theta]\bigg[[\mu_{11}(\Sigma'_2-\mu_{12})+\mu_{22}\mu_{21}][\lambda \mu_{11}+(\lambda+1)\mu_{21}]-\mu_{21}\Sigma'_1\mu_{22}\bigg]\alpha(N)\\
    &+[-\mu_{11}(\Sigma'_2-\mu_{12})-\mu_{21}\mu_{22}-\mu_{21}(s-1)\theta]\Sigma'_2\mu_{11}^{N-s}f(1,s)\\  
\end{align*}
\begin{align*}
     =&[\mu_{11}(\Sigma'_2-\mu_{12})+\mu_{22}\mu_{21}-\mu_{21}\Sigma'_2]\Sigma'_1\mu_{11}^{N-1}\\
  +&\frac{[\mu_{11}(\Sigma'_2-\mu_{12})+\mu_{22}\mu_{21}+\mu_{21}(s-1)\theta]}{\mu_{22}+(s-1)\theta}[\Sigma'_2+(N-1)\theta]\mu_{22}\times\\
  &\bigg[\sum_{k=s+1}^{N}\mu_{11}^{k-s-1}[f(k-s,N)-f(k,N)f(1,s+1)]\bigg]\\
  +&\bigg[\mu_{11}(\Sigma'_2-\mu_{12})+\mu_{22}\mu_{21}-\mu_{22}\frac{[\mu_{11}(\Sigma'_2-\mu_{12})+\mu_{22}\mu_{21}+\mu_{21}(s-1)\theta]}{\mu_{22}+(s-1)\theta}\bigg][\Sigma'_2+(N-1)\theta]\times \\
     &\sum_{k''=s+1}^{N}\mu_{11}^{k''-s-1}f(k''-s,N)\\
    &+ [\mu_{11}(\Sigma'_2-\mu_{12})+\mu_{22}\mu_{21}]\mu_{11}^{N-s}\bigg[\Sigma'_2\sum_{k'=2}^s\mu_{11}^{k'-1}f(k'+N-s,N)-\Sigma_1'[\Sigma'_2-\mu_{22}]\alpha(s) \bigg]\\
    &-(\mu_{21}\Sigma_2'-\mu_{22}\mu_{21})(s-1)\theta\mu_{11}\mu_{11}^{N-s} \alpha(s)\\
    &-(\mu_{21}\Sigma_2'-\mu_{22}\mu_{21})(s-1)\theta[\lambda \mu_{11}+(\lambda+1)\mu_{21}]\mu_{11}^{N-s} \alpha(s)\\
    &+[\mu_{11}(\Sigma'_2-\mu_{12})+\mu_{22}\mu_{21}][\Sigma'_2+(N-1)\theta]\bigg[\mu_{11}^{N-s}f(N+1-s,N)\bigg] \\  
     &+(\mu_{11}(\Sigma'_2-\mu_{12})+\mu_{22}\mu_{21}-\mu_{21}\Sigma_2'+\mu_{22}\mu_{21})(N-1)\theta\bigg[\sum_{k'=2}^s\mu_{11}^{k'+N-1-s}f(k'+N-s,N) \bigg] \\
     &+(\mu_{21}\Sigma'_2-\mu_{22}\mu_{21})(N-1)\theta\bigg[\sum_{k'=2}^s\mu_{11}^{k'+N-1-s}f(k'+N-s,N) \bigg] \\
    &+[\Sigma'_2+(N-1)\theta]\bigg[[\mu_{11}(\Sigma'_2-\mu_{12})+\mu_{22}\mu_{21}][\lambda \mu_{11}+(\lambda+1)\mu_{21}]-\mu_{21}\Sigma'_1\mu_{22}\bigg]\alpha(N)\\
    &+[-\mu_{11}(\Sigma'_2-\mu_{12})-\mu_{21}\mu_{22}-\mu_{21}(s-1)\theta]\Sigma'_2\mu_{11}^{N-s}f(1,s)\\  
     =&[\mu_{11}(\Sigma'_2-\mu_{12})+\mu_{22}\mu_{21}-\mu_{21}\Sigma'_2]\Sigma'_1\mu_{11}^{N-1}\\
  +&\frac{[\mu_{11}(\Sigma'_2-\mu_{12})+\mu_{22}\mu_{21}+\mu_{21}(s-1)\theta]}{\mu_{22}+(s-1)\theta}[\Sigma'_2+(N-1)\theta]\mu_{22}\times \\
  &\bigg[\sum_{k=s+1}^{N}\mu_{11}^{k-s-1}[f(k-s,N)-f(k,N)f(1,s+1)]\bigg]\\
  +&\bigg[\mu_{11}(\Sigma'_2-\mu_{12})+\mu_{22}\mu_{21}-\mu_{22}\frac{[\mu_{11}(\Sigma'_2-\mu_{12})+\mu_{22}\mu_{21}+\mu_{21}(s-1)\theta]}{\mu_{22}+(s-1)\theta}\bigg][\Sigma'_2+(N-1)\theta]\times \\
     &\sum_{k''=s+1}^{N}\mu_{11}^{k''-s-1}f(k''-s,N)\\
    &+ [\mu_{11}(\Sigma'_2-\mu_{12})+\mu_{22}\mu_{21}]\mu_{11}^{N-s}\bigg[\Sigma'_2\sum_{k'=2}^s\mu_{11}^{k'-1}f(k'+N-s,N)-\Sigma_1'[\Sigma'_2-\mu_{22}]\alpha(s) \bigg]\\
    &+(\mu_{21}\Sigma_2'-\mu_{22}\mu_{21})\mu_{11}^{N-s+1}\bigg[(N-1)\theta \sum_{k'=2}^s\mu_{11}^{k'-2}f(k'+N-s,N)- (s-1)\theta\alpha(s)\bigg]\\
    &-(\mu_{21}\Sigma_2'-\mu_{22}\mu_{21})(s-1)\theta[\lambda \mu_{11}+(\lambda+1)\mu_{21}]\mu_{11}^{N-s} \alpha(s)\\
    &+[\mu_{11}(\Sigma'_2-\mu_{12})+\mu_{22}\mu_{21}][\Sigma'_2+(N-1)\theta]\bigg[\mu_{11}^{N-s}f(N+1-s,N)\bigg] \\  
     &+(\mu_{11}(\Sigma'_2-\mu_{12})+\mu_{22}\mu_{21}-\mu_{21}\Sigma_2'+\mu_{22}\mu_{21})(N-1)\theta\bigg[\sum_{k'=2}^s\mu_{11}^{k'+N-1-s}f(k'+N-s,N) \bigg] \\
    &+[\Sigma'_2+(N-1)\theta]\bigg[[\mu_{11}(\Sigma'_2-\mu_{12})+\mu_{22}\mu_{21}][\lambda \mu_{11}+(\lambda+1)\mu_{21}]-\mu_{21}\Sigma'_1\mu_{22}\bigg]\alpha(N)\\
    &+[-\mu_{11}(\Sigma'_2-\mu_{12})-\mu_{21}\mu_{22}-\mu_{21}(s-1)\theta]\Sigma'_2\mu_{11}^{N-s}f(1,s)\\  
     =&[\mu_{11}(\Sigma'_2-\mu_{12})+\mu_{22}\mu_{21}-\mu_{21}\Sigma'_2]\Sigma'_1\mu_{11}^{N-1}\\
  +&\frac{[\mu_{11}(\Sigma'_2-\mu_{12})+\mu_{22}\mu_{21}+\mu_{21}(s-1)\theta]}{\mu_{22}+(s-1)\theta}[\Sigma'_2+(N-1)\theta]\mu_{22}\times\\
  &\bigg[\sum_{k=s+1}^{N}\mu_{11}^{k-s-1}[f(k-s,N)-f(k,N)f(1,s+1)]\bigg]\\
  +&\bigg[\mu_{11}(\Sigma'_2-\mu_{12})+\mu_{22}\mu_{21}-\mu_{22}\frac{[\mu_{11}(\Sigma'_2-\mu_{12})+\mu_{22}\mu_{21}+\mu_{21}(s-1)\theta]}{\mu_{22}+(s-1)\theta}\bigg][\Sigma'_2+(N-1)\theta]\times \\
     &\sum_{k''=s+1}^{N}\mu_{11}^{k''-s-1}f(k''-s,N)\\
    &+ [\mu_{11}(\Sigma'_2-\mu_{12})+\mu_{22}\mu_{21}]\mu_{11}^{N-s}\bigg[\Sigma'_2\sum_{k'=2}^s\mu_{11}^{k'-1}f(k'+N-s,N)-\Sigma_1'[\Sigma'_2-\mu_{22}]\alpha(s) \bigg]\\
    &+(\mu_{21}\Sigma_2'-\mu_{22}\mu_{21})\mu_{11}^{N-s+1}\bigg[(N-1)\theta \sum_{k'=2}^s\mu_{11}^{k'-2}f(k'+N-s,N)- (s-1)\theta\alpha(s)\bigg]\\
    &-(\mu_{21}\Sigma_2'-\mu_{22}\mu_{21})(s-1)\theta[\lambda \mu_{11}+(\lambda+1)\mu_{21}]\mu_{11}^{N-s} \alpha(s)\\
    &+(N-1)\theta[[\mu_{11}(\Sigma'_2-\mu_{12})+\mu_{22}\mu_{21}][\lambda \mu_{11}+(\lambda+1)\mu_{21}]-\mu_{21}\Sigma'_1\mu_{22}]\alpha(N)\\
    &+\mu_{22}\mu_{21}(N-1)\theta\bigg[\sum_{k'=2}^s\mu_{11}^{k'+N-1-s}f(k'+N-s,N) \bigg] \\
    &+[\mu_{11}(\Sigma'_2-\mu_{12})+\mu_{22}\mu_{21}][\Sigma'_2+(N-1)\theta]\bigg[\mu_{11}^{N-s}f(N+1-s,N)\bigg] \\  
     &+(\mu_{11}(\Sigma'_2-\mu_{12})+\mu_{22}\mu_{21}-\mu_{21}\Sigma_2')(N-1)\theta\bigg[\sum_{k'=2}^s\mu_{11}^{k'+N-1-s}f(k'+N-s,N) \bigg] \\
    &+\Sigma'_2\bigg[[\mu_{11}(\Sigma'_2-\mu_{12})+\mu_{22}\mu_{21}][\lambda \mu_{11}+(\lambda+1)\mu_{21}]-\mu_{21}\Sigma'_1\mu_{22}\bigg]\alpha(N)\\
    &+[-\mu_{11}(\Sigma'_2-\mu_{12})-\mu_{21}\mu_{22}-\mu_{21}(s-1)\theta]\Sigma'_2\mu_{11}^{N-s}f(1,s)\\  
     =&[\mu_{11}(\Sigma'_2-\mu_{12})+\mu_{22}\mu_{21}-\mu_{21}\Sigma'_2]\Sigma'_1\mu_{11}^{N-1}\\
  +&\frac{[\mu_{11}(\Sigma'_2-\mu_{12})+\mu_{22}\mu_{21}+\mu_{21}(s-1)\theta]}{\mu_{22}+(s-1)\theta}[\Sigma'_2+(N-1)\theta]\mu_{22}\times\\
  &\bigg[\sum_{k=s+1}^{N}\mu_{11}^{k-s-1}[f(k-s,N)-f(k,N)f(1,s+1)]\bigg]\\
  +&\bigg[\mu_{11}(\Sigma'_2-\mu_{12})+\mu_{22}\mu_{21}-\mu_{22}\frac{[\mu_{11}(\Sigma'_2-\mu_{12})+\mu_{22}\mu_{21}+\mu_{21}(s-1)\theta]}{\mu_{22}+(s-1)\theta}\bigg][\Sigma'_2+(N-1)\theta]\times \\
     &\sum_{k''=s+1}^{N}\mu_{11}^{k''-s-1}f(k''-s,N)\\
    &+ [\mu_{11}(\Sigma'_2-\mu_{12})+\mu_{22}\mu_{21}]\mu_{11}^{N-s}\bigg[\Sigma'_2\sum_{k'=2}^s\mu_{11}^{k'-1}f(k'+N-s,N)-\Sigma_1'[\Sigma'_2-\mu_{22}]\alpha(s) \bigg]\\
    &+(\mu_{21}\Sigma_2'-\mu_{22}\mu_{21})\mu_{11}^{N-s+1}\bigg[(N-1)\theta \sum_{k'=2}^s\mu_{11}^{k'-2}f(k'+N-s,N)- (s-1)\theta\alpha(s)\bigg]\\
    &-(\mu_{21}\Sigma_2'-\mu_{22}\mu_{21})(s-1)\theta[\lambda \mu_{11}+(\lambda+1)\mu_{21}]\mu_{11}^{N-s} \alpha(s)\\
    &+(\mu_{11}(\Sigma'_2-\mu_{12})+\mu_{22}\mu_{21}-\mu_{21}\Sigma'_2)(s-1)\theta[\lambda \mu_{11}+(\lambda+1)\mu_{21}]\mu_{11}^{N-s} \alpha(s)\\
    &-[\mu_{11}(\Sigma'_2-\mu_{12})+\mu_{22}\mu_{21}-\mu_{21}\Sigma'_2](s-1)\theta[\lambda \mu_{11}+(\lambda+1)\mu_{21}]\mu_{11}^{N-s} \alpha(s)\\
    &+(N-1)\theta[[\mu_{11}(\Sigma'_2-\mu_{12})+\mu_{22}\mu_{21}][\lambda \mu_{11}+(\lambda+1)\mu_{21}]-\mu_{21}\Sigma'_1\mu_{22}]\alpha(N)\\
    &+\mu_{22}\mu_{21}(N-1)\theta\bigg[\sum_{k'=2}^s\mu_{11}^{k'+N-1-s}f(k'+N-s,N) \bigg] \\
    &+[\mu_{11}(\Sigma'_2-\mu_{12})+\mu_{22}\mu_{21}][\Sigma'_2+(N-1)\theta]\bigg[\mu_{11}^{N-s}f(N+1-s,N)\bigg] \\  
     &+(\mu_{11}(\Sigma'_2-\mu_{12})+\mu_{22}\mu_{21}-\mu_{21}\Sigma_2')(N-1)\theta\bigg[\sum_{k'=2}^s\mu_{11}^{k'+N-1-s}f(k'+N-s,N) \bigg] \\
    &+\Sigma'_2\bigg[[\mu_{11}(\Sigma'_2-\mu_{12})+\mu_{22}\mu_{21}][\lambda \mu_{11}+(\lambda+1)\mu_{21}]-\mu_{21}\Sigma'_1\mu_{22}\bigg]\alpha(N)\\
    &+[-\mu_{11}(\Sigma'_2-\mu_{12})-\mu_{21}\mu_{22}-\mu_{21}(s-1)\theta]\Sigma'_2\mu_{11}^{N-s}f(1,s)\\  
     =&[\mu_{11}(\Sigma'_2-\mu_{12})+\mu_{22}\mu_{21}-\mu_{21}\Sigma'_2]\Sigma'_1\mu_{11}^{N-1}\\
  +&\frac{[\mu_{11}(\Sigma'_2-\mu_{12})+\mu_{22}\mu_{21}+\mu_{21}(s-1)\theta]}{\mu_{22}+(s-1)\theta}[\Sigma'_2+(N-1)\theta]\mu_{22}\times\\
  &\bigg[\sum_{k=s+1}^{N}\mu_{11}^{k-s-1}[f(k-s,N)-f(k,N)f(1,s+1)]\bigg]\\
  +&\bigg[\mu_{11}(\Sigma'_2-\mu_{12})+\mu_{22}\mu_{21}-\mu_{22}\frac{[\mu_{11}(\Sigma'_2-\mu_{12})+\mu_{22}\mu_{21}+\mu_{21}(s-1)\theta]}{\mu_{22}+(s-1)\theta}\bigg][\Sigma'_2+(N-1)\theta]\times \\
     &\sum_{k''=s+1}^{N}\mu_{11}^{k''-s-1}f(k''-s,N)\\
    &+ [\mu_{11}(\Sigma'_2-\mu_{12})+\mu_{22}\mu_{21}]\mu_{11}^{N-s}\bigg[\Sigma'_2\sum_{k'=2}^s\mu_{11}^{k'-1}f(k'+N-s,N)-\Sigma_1'[\Sigma'_2-\mu_{22}]\alpha(s) \bigg]\\
    &+(\mu_{21}\Sigma_2'-\mu_{22}\mu_{21})\mu_{11}^{N-s+1}\bigg[(N-1)\theta \sum_{k'=2}^s\mu_{11}^{k'-2}f(k'+N-s,N)- (s-1)\theta\alpha(s)\bigg]\\
    &+[\mu_{11}(\Sigma'_2-\mu_{12})+\mu_{22}\mu_{21}-\mu_{21}\Sigma'_2](s-1)\theta[\lambda \mu_{11}+(\lambda+1)\mu_{21}]\mu_{11}^{N-s} \alpha(s)\\
    &-(\mu_{11}(\Sigma'_2-\mu_{12}))(s-1)\theta[\lambda \mu_{11}+(\lambda+1)\mu_{21}]\mu_{11}^{N-s} \alpha(s)\\
    &+(N-1)\theta[[\mu_{11}(\Sigma'_2-\mu_{12})+\mu_{22}\mu_{21}][\lambda \mu_{11}+(\lambda+1)\mu_{21}]-\mu_{21}\Sigma'_1\mu_{22}]\alpha(N)\\
    &+\mu_{22}\mu_{21}(N-1)\theta\bigg[\sum_{k'=2}^s\mu_{11}^{k'+N-1-s}f(k'+N-s,N) \bigg] \\
    &+[\mu_{11}(\Sigma'_2-\mu_{12})+\mu_{22}\mu_{21}][\Sigma'_2+(N-1)\theta]\bigg[\mu_{11}^{N-s}f(N+1-s,N)\bigg] \\  
     &+(\mu_{11}(\Sigma'_2-\mu_{12})+\mu_{22}\mu_{21}-\mu_{21}\Sigma_2')(N-1)\theta\bigg[\sum_{k'=2}^s\mu_{11}^{k'+N-1-s}f(k'+N-s,N) \bigg] \\
    &+\Sigma'_2\bigg[[\mu_{11}(\Sigma'_2-\mu_{12})+\mu_{22}\mu_{21}][\lambda \mu_{11}+(\lambda+1)\mu_{21}]-\mu_{21}\Sigma'_1\mu_{22}\bigg]\alpha(N)\\
    &+[-\mu_{11}(\Sigma'_2-\mu_{12})-\mu_{21}\mu_{22}-\mu_{21}(s-1)\theta]\Sigma'_2\mu_{11}^{N-s}f(1,s)\\  
     =&[\mu_{11}(\Sigma'_2-\mu_{12})+\mu_{22}\mu_{21}-\mu_{21}\Sigma'_2]\Sigma'_1\mu_{11}^{N-1}\\
  +&\frac{[\mu_{11}(\Sigma'_2-\mu_{12})+\mu_{22}\mu_{21}+\mu_{21}(s-1)\theta]}{\mu_{22}+(s-1)\theta}[\Sigma'_2+(N-1)\theta]\mu_{22}\times\\
  &\bigg[\sum_{k=s+1}^{N}\mu_{11}^{k-s-1}[f(k-s,N)-f(k,N)f(1,s+1)]\bigg]\\
  +&\bigg[\mu_{11}(\Sigma'_2-\mu_{12})+\mu_{22}\mu_{21}-\mu_{22}\frac{[\mu_{11}(\Sigma'_2-\mu_{12})+\mu_{22}\mu_{21}+\mu_{21}(s-1)\theta]}{\mu_{22}+(s-1)\theta}\bigg][\Sigma'_2+(N-1)\theta]\times \\
     &\sum_{k''=s+1}^{N}\mu_{11}^{k''-s-1}f(k''-s,N)\\
    &+ [\mu_{11}(\Sigma'_2-\mu_{12})+\mu_{22}\mu_{21}]\mu_{11}^{N-s}\bigg[\Sigma'_2\sum_{k'=2}^s\mu_{11}^{k'-1}f(k'+N-s,N)-\Sigma_1'[\Sigma'_2-\mu_{22}]\alpha(s) \bigg]\\
    &+(\mu_{21}\Sigma_2'-\mu_{22}\mu_{21})\mu_{11}^{N-s+1}\bigg[(N-1)\theta \sum_{k'=2}^s\mu_{11}^{k'-2}f(k'+N-s,N)- (s-1)\theta\alpha(s)\bigg]\\
    &+[\mu_{11}(\Sigma'_2-\mu_{12})+\mu_{22}\mu_{21}-\mu_{21}\Sigma'_2](s-1)\theta[\lambda \mu_{11}+(\lambda+1)\mu_{21}]\mu_{11}^{N-s} \alpha(s)\\
    &-(\mu_{11}(\Sigma'_2-\mu_{12}))(s-1)\theta[\lambda \mu_{11}+(\lambda+1)\mu_{21}]\mu_{11}^{N-s} \alpha(s)\\
    &+(N-1)\theta\bigg[[\mu_{11}(\Sigma'_2-\mu_{12})+\mu_{22}\mu_{21}][\lambda \mu_{11}+(\lambda+1)\mu_{21}]-\mu_{21}\Sigma'_1\mu_{22}\bigg]\alpha(N)\\
    &+\mu_{22}\mu_{21}(N-1)\theta\bigg[\sum_{k'=2}^s\mu_{11}^{k'+N-1-s}f(k'+N-s,N) \bigg] \\
     &+(\mu_{11}(\Sigma'_2-\mu_{12})+\mu_{22}\mu_{21}-\mu_{21}\Sigma_2')(N-1)\theta\bigg[\sum_{k'=2}^s\mu_{11}^{k'+N-1-s}f(k'+N-s,N) \bigg] \\
    &+\Sigma'_2\bigg[[\mu_{11}(\Sigma'_2-\mu_{12})+\mu_{22}\mu_{21}][\lambda \mu_{11}+(\lambda+1)\mu_{21}]-\mu_{21}\Sigma'_1\mu_{22}\bigg]\alpha(N)\\
    &+[-\mu_{11}(\Sigma'_2-\mu_{12})-\mu_{21}\mu_{22}-\mu_{21}(s-1)\theta]\Sigma'_2\mu_{11}^{N-s}f(1,s)\\  
    &+[\mu_{11}(\Sigma'_2-\mu_{12})+\mu_{22}\mu_{21}][\Sigma'_2+(N-1)\theta]\bigg[\mu_{11}^{N-s}f(N+1-s,N)\bigg] \\  
     =&[\mu_{11}(\Sigma'_2-\mu_{12})+\mu_{22}\mu_{21}-\mu_{21}\Sigma'_2]\Sigma'_1\mu_{11}^{N-1}\\
  +&\frac{[\mu_{11}(\Sigma'_2-\mu_{12})+\mu_{22}\mu_{21}+\mu_{21}(s-1)\theta]}{\mu_{22}+(s-1)\theta}[\Sigma'_2+(N-1)\theta]\mu_{22}\times\\
  &\bigg[\sum_{k=s+1}^{N}\mu_{11}^{k-s-1}[f(k-s,N)-f(k,N)f(1,s+1)]\bigg]\\
  +&\bigg[\mu_{11}(\Sigma'_2-\mu_{12})+\mu_{22}\mu_{21}-\mu_{22}\frac{[\mu_{11}(\Sigma'_2-\mu_{12})+\mu_{22}\mu_{21}+\mu_{21}(s-1)\theta]}{\mu_{22}+(s-1)\theta}\bigg][\Sigma'_2+(N-1)\theta]\times \\
     &\sum_{k''=s+1}^{N}\mu_{11}^{k''-s-1}f(k''-s,N)\\
    &+ [\mu_{11}(\Sigma'_2-\mu_{12})+\mu_{22}\mu_{21}]\mu_{11}^{N-s}\bigg[\Sigma'_2\sum_{k'=2}^s\mu_{11}^{k'-1}f(k'+N-s,N)-\Sigma_1'[\Sigma'_2-\mu_{22}]\alpha(s) \bigg]\\
    &+(\mu_{21}\Sigma_2'-\mu_{22}\mu_{21})\mu_{11}^{N-s+1}\bigg[(N-1)\theta \sum_{k'=2}^s\mu_{11}^{k'-2}f(k'+N-s,N)- (s-1)\theta\alpha(s)\bigg]\\
    &+[\mu_{11}(\Sigma'_2-\mu_{12})+\mu_{22}\mu_{21}-\mu_{21}\Sigma'_2](s-1)\theta[\lambda \mu_{11}+(\lambda+1)\mu_{21}]\mu_{11}^{N-s} \alpha(s)\\
    &-(\mu_{11}(\Sigma'_2-\mu_{12}))(s-1)\theta[\lambda \mu_{11}+(\lambda+1)\mu_{21}]\mu_{11}^{N-s} \alpha(s)\\
    &+(N-1)\theta\bigg[[\mu_{11}(\Sigma'_2-\mu_{12})+\mu_{22}\mu_{21}][\lambda \mu_{11}+(\lambda+1)\mu_{21}]-\mu_{21}\Sigma'_1\mu_{22}\bigg]\alpha(N)\\
    &+\mu_{22}\mu_{21}(N-1)\theta\bigg[\sum_{k'=2}^s\mu_{11}^{k'+N-1-s}f(k'+N-s,N) \bigg] \\
     &+(\mu_{11}(\Sigma'_2-\mu_{12})+\mu_{22}\mu_{21}-\mu_{21}\Sigma_2')(N-1)\theta\bigg[\sum_{k'=2}^s\mu_{11}^{k'+N-1-s}f(k'+N-s,N) \bigg] \\
    &+\Sigma'_2\bigg[[\mu_{11}(\Sigma'_2-\mu_{12})+\mu_{22}\mu_{21}][\lambda \mu_{11}+(\lambda+1)\mu_{21}]-\mu_{21}\Sigma'_1\mu_{22}\bigg]\alpha(N)\\
    &-[\mu_{11}(\Sigma'_2-\mu_{12})+\mu_{21}\mu_{22}]\Sigma'_2\mu_{11}^{N-s}f(1,s)\\  
    &-\mu_{21}(s-1)\theta\Sigma'_2\mu_{11}^{N-s}f(1,s)\\
    &+[\mu_{11}(\Sigma'_2-\mu_{12})+\mu_{22}\mu_{21}]\Sigma'_2\bigg[\mu_{11}^{N-s}f(N+1-s,N)\bigg] \\  
    &+[\mu_{11}(\Sigma'_2-\mu_{12})+\mu_{22}\mu_{21}](N-1)\theta\bigg[\mu_{11}^{N-s}f(N+1-s,N)\bigg] \\
     =&[\mu_{11}(\Sigma'_2-\mu_{12})+\mu_{22}\mu_{21}-\mu_{21}\Sigma'_2]\Sigma'_1\mu_{11}^{N-1}\\
  +&\frac{[\mu_{11}(\Sigma'_2-\mu_{12})+\mu_{22}\mu_{21}+\mu_{21}(s-1)\theta]}{\mu_{22}+(s-1)\theta}[\Sigma'_2+(N-1)\theta]\mu_{22}\times \\
  &\bigg[\sum_{k=s+1}^{N}\mu_{11}^{k-s-1}[f(k-s,N)-f(k,N)f(1,s+1)]\bigg]\\
  +&\bigg[\mu_{11}(\Sigma'_2-\mu_{12})+\mu_{22}\mu_{21}-\mu_{22}\frac{[\mu_{11}(\Sigma'_2-\mu_{12})+\mu_{22}\mu_{21}+\mu_{21}(s-1)\theta]}{\mu_{22}+(s-1)\theta}\bigg][\Sigma'_2+(N-1)\theta]\times \\
     &\sum_{k''=s+1}^{N}\mu_{11}^{k''-s-1}f(k''-s,N)\\
    &+ [\mu_{11}(\Sigma'_2-\mu_{12})+\mu_{22}\mu_{21}]\mu_{11}^{N-s}\bigg[\Sigma'_2\sum_{k'=2}^s\mu_{11}^{k'-1}f(k'+N-s,N)-\Sigma_1'[\Sigma'_2-\mu_{22}]\alpha(s) \bigg]\\
    &+(\mu_{21}\Sigma_2'-\mu_{22}\mu_{21})\mu_{11}^{N-s+1}\bigg[(N-1)\theta \sum_{k'=2}^s\mu_{11}^{k'-2}f(k'+N-s,N)- (s-1)\theta\alpha(s)\bigg]\\
    &+[\mu_{11}(\Sigma'_2-\mu_{12})+\mu_{22}\mu_{21}-\mu_{21}\Sigma'_2](s-1)\theta[\lambda \mu_{11}+(\lambda+1)\mu_{21}]\mu_{11}^{N-s} \alpha(s)\\
    &-(\mu_{11}(\Sigma'_2-\mu_{12}))(s-1)\theta[\lambda \mu_{11}+(\lambda+1)\mu_{21}]\mu_{11}^{N-s} \alpha(s)\\
    &+(N-1)\theta\bigg[[\mu_{11}(\Sigma'_2-\mu_{12})+\mu_{22}\mu_{21}][\lambda \mu_{11}+(\lambda+1)\mu_{21}]-\mu_{21}\Sigma'_1\mu_{22}\bigg]\alpha(N)\\
    &+\mu_{22}\mu_{21}(N-1)\theta\bigg[\sum_{k'=2}^s\mu_{11}^{k'+N-1-s}f(k'+N-s,N) \bigg] \\
     &+(\mu_{11}(\Sigma'_2-\mu_{12})+\mu_{22}\mu_{21}-\mu_{21}\Sigma_2')(N-1)\theta\bigg[\sum_{k'=2}^s\mu_{11}^{k'+N-1-s}f(k'+N-s,N) \bigg] \\
    &+\Sigma'_2\bigg[[\mu_{11}(\Sigma'_2-\mu_{12})+\mu_{22}\mu_{21}][\lambda \mu_{11}+(\lambda+1)\mu_{21}]-\mu_{21}\Sigma'_1\mu_{22}\bigg]\alpha(N)\\
    &-\mu_{21}(s-1)\theta\Sigma'_2\mu_{11}^{N-s}f(1,s)\\
    &+[\mu_{11}(\Sigma'_2-\mu_{12})+\mu_{22}\mu_{21}]\Sigma'_2\mu_{11}^{N-s}\bigg[f(N+1-s,N)-f(1,s)\bigg] \\  
    &+[\mu_{11}(\Sigma'_2-\mu_{12})+\mu_{22}\mu_{21}-\mu_{21}\Sigma'_2](N-1)\theta\bigg[\mu_{11}^{N-s}f(N+1-s,N)\bigg]\\
    &+\mu_{21}\Sigma'_2 (N-1)\theta\bigg[\mu_{11}^{N-s}f(N+1-s,N)\bigg]\\
     =&[\mu_{11}(\Sigma'_2-\mu_{12})+\mu_{22}\mu_{21}-\mu_{21}\Sigma'_2]\Sigma'_1\mu_{11}^{N-1}\refstepcounter{equation}\tag{\theequation} \label{t1}\\
  +&\frac{[\mu_{11}(\Sigma'_2-\mu_{12})+\mu_{22}\mu_{21}+\mu_{21}(s-1)\theta]}{\mu_{22}+(s-1)\theta}[\Sigma'_2+(N-1)\theta]\mu_{22} \nonumber\\
  &\times \bigg[\sum_{k=s+1}^{N}\mu_{11}^{k-s-1}[f(k-s,N)-f(k,N)f(1,s+1)]\bigg]\refstepcounter{equation}\tag{\theequation} \label{t2}\\
  +&\bigg[\mu_{11}(\Sigma'_2-\mu_{12})+\mu_{22}\mu_{21}-\mu_{22}\frac{[\mu_{11}(\Sigma'_2-\mu_{12})+\mu_{22}\mu_{21}+\mu_{21}(s-1)\theta]}{\mu_{22}+(s-1)\theta}\bigg][\Sigma'_2+(N-1)\theta]\nonumber \\
     &\times \sum_{k=s+1}^{N}\mu_{11}^{k-s-1}f(k-s,N)\refstepcounter{equation}\tag{\theequation} \label{t3}\\
    &+ [\mu_{11}(\Sigma'_2-\mu_{12})+\mu_{22}\mu_{21}]\mu_{11}^{N-s}\bigg[\Sigma'_2\sum_{k=2}^s\mu_{11}^{k-1}f(k+N-s,N)-\Sigma_1'[\Sigma'_2-\mu_{22}]\alpha(s) \bigg] \refstepcounter{equation}\tag{\theequation}\label{t4}\\
    &+\Sigma'_2\big[[\mu_{11}(\Sigma'_2-\mu_{12})+\mu_{22}\mu_{21}](\Sigma_1'-\mu_{11})-\mu_{21}\Sigma'_1\mu_{22}\big]\alpha(N) \refstepcounter{equation}\tag{\theequation}\label{t5}\\
    &+\mu_{21}(\Sigma_2'-\mu_{22})\mu_{11}^{N-s+1}\bigg[(N-1)\theta \sum_{k=2}^s\mu_{11}^{k-2}f(k+N-s,N)- (s-1)\theta\alpha(s)\bigg] \refstepcounter{equation}\tag{\theequation}\label{t6}\\
    &+[\mu_{11}(\Sigma'_2-\mu_{12})+\mu_{22}\mu_{21}-\mu_{21}\Sigma'_2](s-1)\theta(\Sigma_1'-\mu_{11})\mu_{11}^{N-s} \alpha(s) \refstepcounter{equation}\tag{\theequation}\label{t7}\\
    &-\mu_{11}(\Sigma'_2-\mu_{12})(s-1)\theta(\Sigma_1'-\mu_{11})\mu_{11}^{N-s} \alpha(s) \refstepcounter{equation}\tag{\theequation}\label{t8}\\
    &+(N-1)\theta\big[[\mu_{11}(\Sigma'_2-\mu_{12})+\mu_{22}\mu_{21}](\Sigma_1'-\mu_{11})-\mu_{21}\Sigma'_1\mu_{22}\big]\alpha(N)\refstepcounter{equation}\tag{\theequation} \label{t9}\\
    &+\mu_{22}\mu_{21}(N-1)\theta\bigg[\sum_{k=2}^s\mu_{11}^{N+k-s-1}f(k+N-s,N) \bigg] \refstepcounter{equation}\tag{\theequation}\label{t10}\\
     &+[\mu_{11}(\Sigma'_2-\mu_{12})+\mu_{22}\mu_{21}-\mu_{21}\Sigma_2'](N-1)\theta\bigg[\sum_{k=2}^s\mu_{11}^{N+k-s-1}f(k+N-s,N) \bigg] \refstepcounter{equation}\tag{\theequation}\label{t11} \\
    &+\mu_{21}\Sigma'_2\mu_{11}^{N-s}\theta \bigg[(N-1)f(N+1-s,N)-(s-1)f(1,s) \bigg] \refstepcounter{equation}\tag{\theequation}\label{t12}\\
    &+[\mu_{11}(\Sigma'_2-\mu_{12})+\mu_{22}\mu_{21}]\Sigma'_2\mu_{11}^{N-s}\bigg[f(N+1-s,N)-f(1,s)\bigg] \refstepcounter{equation}\tag{\theequation}\label{t13}\\  
    &+[\mu_{11}(\Sigma'_2-\mu_{12})+\mu_{22}\mu_{21}-\mu_{21}\Sigma'_2](N-1)\theta\mu_{11}^{N-s}f(N+1-s,N). \refstepcounter{equation}\tag{\theequation}\label{t14}
\end{align*}
We now show that the sum of the terms \eqref{t1}-\eqref{t14} is positive. Note that as $N \geq 2$ for this case to be possible, we have $\tau'(2) \geq 0$. Note that
\begin{align}
    \tau'(2) 
    &=[\Sigma'_2+\theta]\Sigma'_2\mu_{22}-[\Sigma'_2+\theta][\Sigma'_1(\Sigma'_2-\mu_{22})+\Sigma'_2\mu_{22}]+\mu_{11}(\Sigma'_2)^2 \nonumber \\
     &=-[\Sigma'_2+\theta][\Sigma'_1(\Sigma'_2-\mu_{22})]+\mu_{11}(\Sigma'_2)^2 \nonumber\\
   &= -\Sigma'_2 \Sigma'_1(\Sigma'_2-\mu_{22})+(\Sigma'_2)^2\mu_{11} -\theta \Sigma'_1(\Sigma'_2-\mu_{22}) \label{tau'(2.2)supp}\\
   &=-\Sigma'_2 \Sigma'_1(\lambda \Sigma_2+\mu_{12})+(\Sigma'_2)^2\mu_{11} -\theta \Sigma'_1(\Sigma'_2-\mu_{22}) \nonumber\\
   &=\Sigma'_2(1+\lambda)[-\lambda\Sigma_1\Sigma_2-\mu_{12}\Sigma_1+\Sigma_2\mu_{11}] -\theta \Sigma'_1(\Sigma'_2-\mu_{22}) \nonumber\\
   &=\Sigma'_2(1+\lambda)[-\lambda\Sigma_1\Sigma_2+\mu_{11}\mu_{22}-\mu_{21}\mu_{12}] -\theta \Sigma'_1(\Sigma'_2-\mu_{22}). \label{tau'(2)supp}
\end{align}
That is, $ \tau'(2) \geq 0$ implies that 
    $-\lambda\Sigma_1 \Sigma_2+\mu_{11}\mu_{22}-\mu_{21}\mu_{12} \geq 0$, and hence 
\begin{align}\label{lambda}
    \lambda \leq \frac{\mu_{11}\mu_{22}-\mu_{21}\mu_{12}}{\Sigma_1 \Sigma_2}.
\end{align}

For \eqref{t1} we have
\begin{align}\label{sub}
    &\mu_{11}(\Sigma'_2-\mu_{12})+\mu_{22}\mu_{21}-\mu_{21}\Sigma'_2  \nonumber \\
    =&\mu_{11}[(\lambda+1)\mu_{22}+\lambda \mu_{12}]+\mu_{22}\mu_{21}-\mu_{21}[(\lambda+1)\mu_{22}+(\lambda+1)\mu_{12}] \nonumber \\
    =&(\lambda+1)\mu_{11}\mu_{22}+\lambda \mu_{11}\mu_{12}+\mu_{22}\mu_{21}-(\lambda+1)\mu_{21}\mu_{22}-(\lambda+1)\mu_{21}\mu_{12} \nonumber \\
    =&(\lambda+1)(\mu_{11}\mu_{22}-\mu_{21}\mu_{12})+\lambda(\mu_{11}\mu_{12}-\mu_{21}\mu_{22}) \nonumber \\
    =&(\mu_{11}\mu_{22}-\mu_{21}\mu_{12})+\lambda(\mu_{11}\mu_{22}-\mu_{21}\mu_{12}+\mu_{11}\mu_{12}-\mu_{21}\mu_{22}) \nonumber \\
    =&(\mu_{11}\mu_{22}-\mu_{21}\mu_{12})+\lambda(\mu_{11}\mu_{22}+\mu_{11}\mu_{12})-\lambda\mu_{21}\Sigma_2 \nonumber \\
    \stackrel{(*)}{\geq} &(\mu_{11}\mu_{22}-\mu_{21}\mu_{12})+\lambda(\mu_{11}\mu_{22}+\mu_{11}\mu_{12})-\frac{\mu_{21}}{\Sigma_1}(\mu_{11}\mu_{22}-\mu_{21}\mu_{12}) \nonumber \\
    =&(\mu_{11}\mu_{22}-\mu_{21}\mu_{12}) \bigg[1-\frac{\mu_{21}}{\Sigma_1} \bigg]+\lambda(\mu_{11}\mu_{22}+\mu_{11}\mu_{12}) \nonumber \\
    =&(\mu_{11}\mu_{22}-\mu_{21}\mu_{12})\frac{\mu_{11}}{\Sigma_1}+\lambda(\mu_{11}\mu_{22}+\mu_{11}\mu_{12}) \nonumber \\
    \stackrel{(**)}{\geq} &0,
\end{align}
where in $(*)$ we have used \eqref{lambda} and in $(**)$ we have used $\mu_{11}\mu_{22}-\mu_{21}\mu_{12} \geq 0$. It follows that  
\eqref{t1} is non-negative.

The non-negativity of \eqref{t2} follows because for $s+1 \leq k \leq N$, we have 
\begin{align*}
    &f(k-s,N)-f(1,s+1)f(k,N)\\
    &=f(k-s,k)f(k,N)-f(1,s+1)f(k,N) \\
   & \geq f(1,s+1) f(k,N)-f(1,s+1)f(k,N) = 0.
\end{align*}

The non-negativity of \eqref{t3} follows from
\begin{align*}
    &\mu_{11}(\Sigma'_2-\mu_{12})+\mu_{22}\mu_{21}-\mu_{22}\frac{[\mu_{11}(\Sigma'_2-\mu_{12})+\mu_{22}\mu_{21}+\mu_{21}(s-1)\theta]}{\mu_{22}+(s-1)\theta}\\
   \geq &\frac{\mu_{22}\mu_{21}(s-1)\theta-\mu_{22}\mu_{21}(s-1)\theta+\mu_{22}[\mu_{11}(\Sigma'_2-\mu_{12})+\mu_{22}\mu_{21}]-\mu_{22}[\mu_{11}(\Sigma'_2-\mu_{12})+\mu_{22}\mu_{21}]}{\mu_{22}+(s-1)\theta}\\
   =&0.
\end{align*}

Next we show \eqref{t4}+\eqref{t5} is positive. From \eqref{t4} we have
\begin{align*}
     &[\mu_{11}(\Sigma'_2-\mu_{12})+\mu_{22}\mu_{21}]\mu_{11}^{N-s}\bigg[\Sigma'_2\sum_{k=2}^s\mu_{11}^{k-1}f(k+N-s,N)-\Sigma_1'[\Sigma'_2-\mu_{22}]\alpha(s) \bigg]\\
    \stackrel{(*)}{\geq} & [\mu_{11}(\Sigma'_2-\mu_{12})+\mu_{22}\mu_{21}]\mu_{11}^{N-s}\bigg[\Sigma'_2\mu_{11}\sum_{k=2}^s\mu_{11}^{k-2}f(k,s)-\Sigma_1'[\Sigma'_2-\mu_{22}]\alpha(s) \bigg]\\
     =& [\mu_{11}(\Sigma'_2-\mu_{12})+\mu_{22}\mu_{21}]\mu_{11}^{N-s}\bigg[\Sigma'_2\mu_{11}\alpha(s)-\Sigma_1'[\Sigma'_2-\mu_{22}]\alpha(s)\bigg]\\
     =& [\mu_{11}(\Sigma'_2-\mu_{12})+\mu_{22}\mu_{21}]\mu_{11}^{N-s}\alpha(s)[\Sigma_2'\mu_{11}-\Sigma'_1\Sigma'_2+\Sigma_1'\mu_{22}],
\end{align*}
where $(*)$ follows since 
\begin{align}\label{upper bound alpha}
    \sum_{k=2}^s\mu_{11}^{k-2}f(k+N-s,N) \geq \sum_{k=2}^s\mu_{11}^{k-2}f(k,s)=\alpha(s) \text{ for } s \leq N-1.
\end{align}
And for \eqref{t5} we have 
\begin{align}\label{eqalphan}
    &[\mu_{11}(\Sigma'_2-\mu_{12})+\mu_{22}\mu_{21}](\Sigma_1'-\mu_{11})-\mu_{21}\Sigma'_1\mu_{22} \nonumber \\
    =&[\mu_{11}(1+\lambda)\mu_{22}+\lambda \mu_{11}\mu_{12}+\mu_{22}\mu_{21}][\lambda \mu_{11}+(\lambda+1)\mu_{21}]-(1+\lambda)\mu_{21}^2\mu_{22}-(1+\lambda)\mu_{21}\mu_{22}\mu_{11} \nonumber\\
    \geq & (1+\lambda)\mu_{21}^2\mu_{22}+(1+\lambda)^2\mu_{22}\mu_{11}\mu_{21}-(1+\lambda)\mu_{21}^2\mu_{22}-(1+\lambda)\mu_{21}\mu_{22}\mu_{11} \nonumber\\
    = &0.
\end{align}
It now follows from \eqref{alphaineq} that \eqref{t5} satisfies
\begin{align*}
    &\Sigma'_2\big[[\mu_{11}(\Sigma'_2-\mu_{12})+\mu_{22}\mu_{21}](\Sigma_1'-\mu_{11})-\mu_{21}\Sigma'_1\mu_{22}\big]\alpha(N)\\
    \geq & \Sigma'_2\big[[\mu_{11}(\Sigma'_2-\mu_{12})+\mu_{22}\mu_{21}](\Sigma_1'-\mu_{11})-\mu_{21}\Sigma'_1\mu_{22}\big]\mu_{11}^{N-s} \alpha(s).
\end{align*}
That is 
\begin{align*}
    \eqref{t4}+\eqref{t5} \geq& \mu_{11}^{N-s} \alpha(s)\bigg[\Sigma_2'\big[[\mu_{11}(\Sigma'_2-\mu_{12})+\mu_{22}\mu_{21}](\Sigma_1'-\mu_{11})-\mu_{21}\Sigma'_1\mu_{22}\big]\\
    &+[\mu_{11}(\Sigma'_2-\mu_{12})+\mu_{21}\mu_{22 
 }][\Sigma_2'\mu_{11}-\Sigma'_1\Sigma'_2+\Sigma_1'\mu_{22}]\bigg]\\
 =& \mu_{11}^{N-s} \alpha(s)\bigg[\Sigma_2'\big[[\mu_{11}(\Sigma'_2-\mu_{12})+\mu_{22}\mu_{21}](\Sigma'_1-\mu_{11})-\mu_{21}\Sigma'_1\mu_{22}\big]\\
    &+[\mu_{11}(\Sigma'_2-\mu_{12})+\mu_{21}\mu_{22 
 }][\Sigma_2'(\mu_{11}-\Sigma'_1)+\Sigma_1'\mu_{22}] \bigg]\\  
 =& \mu_{11}^{N-s} \alpha(s)\bigg[\Sigma_2'[\mu_{11}(\Sigma'_2-\mu_{12})+\mu_{22}\mu_{21}](\Sigma'_1-\mu_{11})-\Sigma_2'\Sigma'_1\mu_{22}\mu_{21}\\
 &-\Sigma_2'[\mu_{11}(\Sigma'_2-\mu_{12})+\mu_{22}\mu_{21}](\Sigma'_1-\mu_{11})+[\mu_{11}(\Sigma'_2-\mu_{12})+\mu_{21}\mu_{22 }]\Sigma_1'\mu_{22}\bigg]\\
 =& \mu_{11}^{N-s} \alpha(s)\bigg[-\Sigma_2'\Sigma'_1\mu_{22}\mu_{21}+[\mu_{11}(\Sigma'_2-\mu_{12})+\mu_{21}\mu_{22 }]\Sigma_1'\mu_{22}\bigg]\\
 =& \mu_{11}^{N-s} \alpha(s)\Sigma_1'\mu_{22}[\mu_{11}(\Sigma'_2-\mu_{12})+\mu_{21}\mu_{22 }-\Sigma_2'\mu_{21}] \geq 0,
\end{align*}
where the last inequality follows from \eqref{sub}. Hence $\eqref{t4}+\eqref{t5} \geq 0$.

The non-negativity of  \eqref{t6} follows since $1 \leq s \leq N-1$ implies that
\begin{align*}
   & (N-1)\theta \sum_{k=2}^s\mu_{11}^{k-2}f(k+N-s,N)- (s-1)\theta\alpha(s)\\
  \geq & (s-1)\theta\bigg[\sum_{k=2}^s\mu_{11}^{k-2}f(k+N-s,N)-\alpha(s)\bigg] \geq  0,
\end{align*}
where the last inequality follows from \eqref{upper bound alpha}.

The non-negativity of \eqref{t7} follows directly from \eqref{sub}.

Next we show \eqref{t8}+\eqref{t9}+\eqref{t10} is positive. We have 
\begin{align*}
    \eqref{t8}+\eqref{t9}+\eqref{t10}
    =&-[\mu_{11}(\Sigma'_2-\mu_{12}))(s-1)\theta(\Sigma_1'-\mu_{11})\mu_{11}^{N-s} \alpha(s) \\
    &+(N-1)\theta\big[[\mu_{11}(\Sigma'_2-\mu_{12})+\mu_{22}\mu_{21}](\Sigma_1'-\mu_{11})-\mu_{21}\Sigma'_1\mu_{22}\big]\alpha(N  )\\
    &+\mu_{22}\mu_{21}(N-1)\theta\bigg[\sum_{k=2}^s\mu_{11}^{N+k-s-1}f(k+N-s,N) \bigg] \\
    \stackrel{(*)}{\geq}&-[\mu_{11}(\Sigma'_2-\mu_{12})](s-1)\theta(\Sigma_1'-\mu_{11})\mu_{11}^{N-s} \alpha(s) \\
    &+(s-1)\theta\big[[\mu_{11}(\Sigma'_2-\mu_{12})+\mu_{22}\mu_{21}](\Sigma_1'-\mu_{11})-\mu_{21}\Sigma'_1\mu_{22}\big]\mu_{11}^{N-s} \alpha(s)\\
    &+\mu_{22}\mu_{21}(N-1)\theta\mu_{11}^{N-s+1}\bigg[\sum_{k=2}^s\mu_{11}^{k-2}f(k+N-s,N) \bigg] \\
    \stackrel{(**)}{\geq}&-[\mu_{11}(\Sigma'_2-\mu_{12})](s-1)\theta(\Sigma_1'-\mu_{11})\mu_{11}^{N-s} \alpha(s) \\
    &+(s-1)\theta \big[[\mu_{11}(\Sigma'_2-\mu_{12})+\mu_{22}\mu_{21}](\Sigma_1'-\mu_{11})-\mu_{21}\Sigma'_1\mu_{22}\big]\mu_{11}^{N-s} \alpha(s)\\
    &+\mu_{22}\mu_{21}(s-1)\theta\mu_{11}^{N-s+1}\alpha(s) \\
    =&(s-1)\theta\mu_{11}^{N-s} \alpha(s) \bigg[ -[\mu_{11}(\Sigma'_2-\mu_{12})](\Sigma_1'-\mu_{11})\\
    &+[\mu_{11}(\Sigma'_2-\mu_{12})+\mu_{22}\mu_{21}](\Sigma_1'-\mu_{11})-\mu_{21}\Sigma'_1\mu_{22}+\mu_{22}\mu_{21}\mu_{11}\bigg]\\
     =&(s-1)\theta\mu_{11}^{N-s} \alpha(s) \bigg[\mu_{22}\mu_{21}(\Sigma'_1-\mu_{11})-\mu_{21}\Sigma'_1\mu_{22}+\mu_{22}\mu_{21}\mu_{11}\bigg]\\
     =&0,
\end{align*}
where $(*)$ follows from \eqref{alphaineq} and \eqref{eqalphan} and  $(**)$ follows from \eqref{upper bound alpha} and $s \leq N-1$. Hence $\eqref{t8}+\eqref{t9}+\eqref{t10}\geq 0$.

The non-negativity of \eqref{t11}  follows from \eqref{sub}.

The non-negativity of \eqref{t12} follows since $s \leq N-1$ and hence
\begin{align*}
   (N-1)f(N+1-s,N)-(s-1)f(1,s) \geq (s-1) \bigg[f(N+1-s,N)-f(1,s)\bigg] \geq 0,
\end{align*}
where the last inequality follows from $f(N+1-s,N)-f(1,s) \geq 0$ for $s \leq N-1$.

Similarly, the non-negativity of \eqref{t13} also follows from $f(N+1-s,N)-f(1,s) \geq 0$ for $s \leq N-1$.

Finally, \eqref{t14} is non-negative follows from \eqref{sub}. Hence we have shown $ \Gamma'(s,a_{21}) \geq 0$ and therefore $ \Gamma(s,a_{21}) \geq 0$.
\item[Case 2, Part (iii) :]  $a=a_{11}$.

We have
\begin{align*}
    \Gamma(s,a_{11})=&\mu_{22}+\frac{\mu_{11}}{q}h_N(s+1)+\frac{\mu_{22}+(s-1)\theta}{q}h_N(s-1)+\bigg(1-\frac{\mu_{11}+\mu_{22}+(s-1)\theta}{q} \bigg)h_N(s)\\
    &-\frac{\Sigma'_1}{q}h_N(s+1)-\frac{s\theta}{q}h_N(s-1)-\bigg(1-\frac{\Sigma'_1+s\theta}{q} \bigg)h_N(s)\\
    =&\mu_{22}+\frac{\mu_{11}-\Sigma'_1}{q}h_N(s+1)+\frac{\mu_{22}-\theta}{q}h_N(s-1)+\frac{\Sigma'_1+s\theta-\mu_{11}-\mu_{22}-(s-1)\theta}{q}h_N(s)\\
     =&\mu_{22}-\frac{\Sigma'_1-\mu_{11}}{q}[h_N(s+1)-h_N(s)]-\frac{\mu_{22}}{q}[h_N(s)-h_N(s-1)]+\frac{\theta}{q}[h_N(s)-h_N(s-1)]\\
     =&\Gamma''(s,a_{11})+\frac{\theta}{q}[h_N(s)-h_N(s-1)],
\end{align*}
where 
\begin{align*}
    \Gamma''(s,a_{11})=\mu_{22}-  \frac{\Sigma_1'-\mu_{11}}{q}[h_N(s+1)-h_N(s)]-\frac{\mu_{22}}{q}[h_N(s)-h_N(s-1)].
\end{align*}
Since $h_N(s)-h_N(s-1)\geq 0$ from Lemma \ref{lemma h_N2}, it is sufficient to show that $\Gamma''(s,a_{11})$ is non negative. Using \eqref{Case 2 bias difference},
   we have
   \begin{align*}
        \frac{\Sigma_1'-\mu_{11}}{q}[h_N(s+1)-h_N(s)]=&\frac{(\Sigma'_1-\mu_{11})(g_N-\mu_{22})}{\mu_{11}}\\
        &+\frac{(\Sigma'_1-\mu_{11})[\mu_{22}+(s-1)\theta] }{\mu_{11}q}[h_N(s)-h_N(s-1)].
   \end{align*}
   Hence,
   \begin{align*}
       & \Gamma''(s,a_{11})\\
       =&\mu_{22}- \frac{(\Sigma'_1-\mu_{11})(g_N-\mu_{22})}{\mu_{11}}-\frac{(\Sigma'_1-\mu_{11})[\mu_{22}+(s-1)\theta] }{\mu_{11}q}[h_N(s)-h_N(s-1)]\\
       &-\frac{\mu_{22}}{q}[h_N(s)-h_N(s-1)]\\
       =&\frac{\mu_{22}\mu_{11}+(\mu_{11}-\Sigma'_1)g_N+\Sigma'_1\mu_{22}-\mu_{22}\mu_{11}}{\mu_{11}}\\
       &+\frac{(\mu_{11}-\Sigma'_1)\mu_{22}-\mu_{22}\mu_{11}+(\mu_{11}-\Sigma'_1)(s-1)\theta}{q\mu_{11}}[h_N(s)-h_N(s-1)]\\
       =&\frac{\Sigma'_1\mu_{22}+(\mu_{11}-\Sigma'_1)g_N}{\mu_{11}}+\frac{-\Sigma'_1\mu_{22}+(\mu_{11}-\Sigma'_1)(s-1)\theta}{q\mu_{11}}[h_N(s)-h_N(s-1)]\\
        =&\frac{\Sigma'_1\mu_{22}+(\mu_{11}-\Sigma'_1)g_N}{\mu_{11}}+\frac{-\Sigma'_1\mu_{22}+(\mu_{11}-\Sigma'_1)(s-1)\theta}{\mu_{11}\Delta(N)}\bigg[ \Sigma'_1\mu_{11}^{N-s}[\Sigma'_2-\mu_{22}]\alpha(s)\\
        &+[\Sigma'_2+(N-1)\theta]\mu_{22}\sum_{k=s+1}^{N} \mu_{11}^{k-s-1}f(k,N)f(1,s)+\Sigma'_2\mu_{11}^{N-s}f(1,s)\bigg ],
        \end{align*}
        where in the last equality we have used Lemma \ref{lemma h_N2}. Therefore, \eqref{g_n} yields
\begin{align*}
      &\Gamma''(s,a_{11})\Delta(N)\mu_{11}  \\
      =&\Sigma'_1\mu_{22}\Delta(N)+(\mu_{11}-\Sigma'_1)\beta(N)+[-\Sigma'_1\mu_{22}+(\mu_{11}-\Sigma'_1)(s-1)\theta]\bigg[ \Sigma'_1\mu_{11}^{N-s}[\Sigma'_2-\mu_{22}]\alpha(s)\\
        &+[\Sigma'_2+(N-1)\theta]\mu_{22}\sum_{k=s+1}^{N} \mu_{11}^{k-s-1}f(k,N)f(1,s)+\Sigma'_2\mu_{11}^{N-s}f(1,s)\bigg ].
   \end{align*}

   Now \eqref{beta(n)}-\eqref{Delta(n)} yield
\begin{align*}
     &\Gamma''(s,a_{11})\Delta(N)\mu_{11} \\
    =&\Sigma'_1\mu_{22}[\Sigma_2'+(N-1)\theta][f(1,N)+\Sigma_1'\alpha(N)]+\Sigma'_1\mu_{22}\Sigma_1'\mu_{11}^{N-1}+(\mu_{11}-\Sigma'_1)[\Sigma_2'+(N-1)\theta]\Sigma_1'\mu_{22}\alpha(N)\\
    &+(\mu_{11}-\Sigma'_1)\Sigma_1' \Sigma_2'\mu_{11}^{N-1}+(-\Sigma'_1\mu_{22}+(\mu_{11}-\Sigma'_1)(s-1)\theta)\Sigma'_1\mu_{11}^{N-s}[\Sigma'_2-\mu_{22}]\alpha(s)\\
    &+(-\Sigma'_1\mu_{22}+(\mu_{11}-\Sigma'_1)(s-1)\theta)[\Sigma'_2+(N-1)\theta]\mu_{22}\sum_{k=s+1}^{N} \mu_{11}^{k-s-1}f(k,N)f(1,s)\\
    &+(-\Sigma'_1\mu_{22}+(\mu_{11}-\Sigma'_1)(s-1)\theta)\Sigma'_2\mu_{11}^{N-s}f(1,s)\\
=&\Sigma'_1\mu_{22}\Sigma_1'\mu_{11}^{N-1}+(\mu_{11}-\Sigma'_1)\Sigma_1' \Sigma_2'\mu_{11}^{N-1}+(-\Sigma'_1\mu_{22}+(\mu_{11}-\Sigma'_1)(s-1)\theta)\Sigma'_1\mu_{11}^{N-s}[\Sigma'_2-\mu_{22}]\alpha(s)\\
    &+\Sigma'_1\mu_{22}[\Sigma_2'+(N-1)\theta][f(1,N)+\Sigma_1'\alpha(N)]+(\mu_{11}-\Sigma'_1)[\Sigma_2'+(N-1)\theta]\Sigma_1'\mu_{22}\alpha(N)\\
     &+(-\Sigma'_1\mu_{22}+(\mu_{11}-\Sigma'_1)(s-1)\theta)[\Sigma'_2+(N-1)\theta]\mu_{22}\sum_{k=s+1}^{N} \mu_{11}^{k-s-1}f(k,N)f(1,s)\\
    &+(-\Sigma'_1\mu_{22}+(\mu_{11}-\Sigma'_1)(s-1)\theta)\Sigma'_2\mu_{11}^{N-s}f(1,s)\\
    =&(\Sigma'_1\mu_{22}+\mu_{11}\Sigma_2'-\Sigma'_1\Sigma_2')\Sigma'_1\mu_{11}^{N-1}+(\mu_{11}-\Sigma'_1)[\Sigma'_2-\mu_{22}](s-1)\theta\Sigma'_1\mu_{11}^{N-s}\alpha(s)\\
    &-\Sigma'_1\mu_{22}\Sigma'_1\mu_{11}^{N-s}[\Sigma'_2-\mu_{22}]\alpha(s)\\
     &+\Sigma'_1\mu_{22}[\Sigma_2'+(N-1)\theta][f(1,N)+\Sigma_1'\alpha(N)]+(\mu_{11}-\Sigma'_1)[\Sigma_2'+(N-1)\theta]\Sigma_1'\mu_{22}\alpha(N)\\
     &+(-\Sigma'_1\mu_{22}+(\mu_{11}-\Sigma'_1)(s-1)\theta)[\Sigma'_2+(N-1)\theta]\mu_{22}\sum_{k=s+1}^{N} \mu_{11}^{k-s-1}f(k,N)f(1,s)\\
    &+(-\Sigma'_1\mu_{22}+(\mu_{11}-\Sigma'_1)(s-1)\theta)\Sigma'_2\mu_{11}^{N-s}f(1,s)\\
    =&(\Sigma'_1\mu_{22}+\mu_{11}\Sigma_2'-\Sigma'_1\Sigma_2')\Sigma'_1\mu_{11}^{N-1}+(\Sigma'_1\mu_{22}+\mu_{11}\Sigma_2'-\Sigma'_1\Sigma_2')(s-1)\theta\Sigma'_1\mu_{11}^{N-s}\alpha(s)\\
    &-\mu_{22}(s-1)\theta\Sigma'_1\mu_{11}^{N-s+1}\alpha(s)-\Sigma'_1\mu_{22}\Sigma'_1\mu_{11}^{N-s}[\Sigma'_2-\mu_{22}]\alpha(s)\\
     &+\Sigma'_1\mu_{22}[\Sigma_2'+(N-1)\theta][f(1,N)+\Sigma_1'\alpha(N)]+(\mu_{11}-\Sigma'_1)[\Sigma_2'+(N-1)\theta]\Sigma_1'\mu_{22}\alpha(N)\\
     &+(-\Sigma'_1\mu_{22}+(\mu_{11}-\Sigma'_1)(s-1)\theta)[\Sigma'_2+(N-1)\theta]\mu_{22}\sum_{k=s+1}^{N} \mu_{11}^{k-s-1}f(k,N)f(1,s)\\
    &+(-\Sigma'_1\mu_{22}+(\mu_{11}-\Sigma'_1)(s-1)\theta)\Sigma'_2\mu_{11}^{N-s}f(1,s)\\
    =&(\Sigma'_1\mu_{22}+\mu_{11}\Sigma_2'-\Sigma'_1\Sigma_2')\bigg( \Sigma'_1\mu_{11}^{N-1}+(s-1)\theta\Sigma'_1\mu_{11}^{N-s}\alpha(s)\bigg)\\
     &-\mu_{22}(s-1)\theta\Sigma'_1\mu_{11}^{N-s+1}\alpha(s)-\Sigma'_1\mu_{22}\Sigma'_1\mu_{11}^{N-s}[\Sigma'_2-\mu_{22}]\alpha(s)\\
     &+\Sigma'_1\mu_{22}[\Sigma_2'+(N-1)\theta][f(1,N)+\Sigma_1'\alpha(N)]+(\mu_{11}-\Sigma'_1)[\Sigma_2'+(N-1)\theta]\Sigma_1'\mu_{22}\alpha(N)\\
     &+(-\Sigma'_1\mu_{22}+(\mu_{11}-\Sigma'_1)(s-1)\theta)[\Sigma'_2+(N-1)\theta]\mu_{22}\sum_{k=s+1}^{N} \mu_{11}^{k-s-1}f(k,N)f(1,s)\\
    &+(-\Sigma'_1\mu_{22}+(\mu_{11}-\Sigma'_1)(s-1)\theta)\Sigma'_2\mu_{11}^{N-s}f(1,s)\\
    =&(\Sigma'_1\mu_{22}+\mu_{11}\Sigma_2'-\Sigma'_1\Sigma_2')\bigg( \Sigma'_1\mu_{11}^{N-1}+(s-1)\theta\Sigma'_1\mu_{11}^{N-s}\alpha(s)\bigg)\\
     &-\mu_{22}(s-1)\theta\Sigma'_1\mu_{11}^{N-s+1}\alpha(s)-\Sigma'_1\mu_{22}\Sigma'_1\mu_{11}^{N-s}[\Sigma'_2-\mu_{22}]\alpha(s)\\
     &+(-\Sigma'_1\mu_{22}+(\mu_{11}-\Sigma'_1)(s-1)\theta)[\Sigma'_2+(N-1)\theta]\mu_{22}\sum_{k=s+1}^{N} \mu_{11}^{k-s-1}f(k,N)f(1,s)\\
     &+\Sigma'_1\mu_{22}[\Sigma_2'+(N-1)\theta][f(1,N)+\Sigma_1'\alpha(N)]+(\mu_{11}-\Sigma'_1)[\Sigma_2'+(N-1)\theta]\Sigma_1'\mu_{22}\alpha(N)\\
    &+(-\Sigma'_1\mu_{22}+(\mu_{11}-\Sigma'_1)(s-1)\theta)\Sigma'_2\mu_{11}^{N-s}f(1,s)\\
    =&(\Sigma'_1\mu_{22}+\mu_{11}\Sigma_2'-\Sigma'_1\Sigma_2')\bigg( \Sigma'_1\mu_{11}^{N-1}+(s-1)\theta\Sigma'_1\mu_{11}^{N-s}\alpha(s)\bigg)\\
     &-\mu_{22}(s-1)\theta\Sigma'_1\mu_{11}^{N-s+1}\alpha(s)-\Sigma'_1\mu_{22}\Sigma'_1\mu_{11}^{N-s}[\Sigma'_2-\mu_{22}]\alpha(s)\\
     &+(-\Sigma'_1\mu_{22}-\Sigma'_1(s-1)\theta)[\Sigma'_2+(N-1)\theta]\mu_{22}\sum_{k=s+1}^{N} \mu_{11}^{k-s-1}f(k,N)f(1,s)\\
     &+\mu_{11}(s-1)\theta[\Sigma'_2+(N-1)\theta]\mu_{22}\sum_{k=s+1}^{N} \mu_{11}^{k-s-1}f(k,N)f(1,s)\\
     &+\Sigma'_1\mu_{22}[\Sigma_2'+(N-1)\theta][f(1,N)+\Sigma_1'\alpha(N)]+(\mu_{11}-\Sigma'_1)[\Sigma_2'+(N-1)\theta]\Sigma_1'\mu_{22}\alpha(N)\\
    &+(-\Sigma'_1\mu_{22}+(\mu_{11}-\Sigma'_1)(s-1)\theta)\Sigma'_2\mu_{11}^{N-s}f(1,s)\\
     =&(\Sigma'_1\mu_{22}+\mu_{11}\Sigma_2'-\Sigma'_1\Sigma_2')\bigg( \Sigma'_1\mu_{11}^{N-1}+(s-1)\theta\Sigma'_1\mu_{11}^{N-s}\alpha(s)\bigg)\\
     &-\mu_{22}(s-1)\theta\Sigma'_1\mu_{11}^{N-s+1}\alpha(s)-\Sigma'_1\mu_{22}\Sigma'_1\mu_{11}^{N-s}[\Sigma'_2-\mu_{22}]\alpha(s)\\
     &-\Sigma_1'(\mu_{22}+(s-1)\theta)[\Sigma'_2+(N-1)\theta]\mu_{22}\sum_{k=s+1}^{N} \mu_{11}^{k-s-1}f(k,N)f(1,s)\\
     &+\mu_{11}(s-1)\theta[\Sigma'_2+(N-1)\theta]\mu_{22}\sum_{k=s+1}^{N} \mu_{11}^{k-s-1}f(k,N)f(1,s)\\
     &+\Sigma'_1\mu_{22}[\Sigma_2'+(N-1)\theta][f(1,N)+\Sigma_1'\alpha(N)]+(\mu_{11}-\Sigma'_1)[\Sigma_2'+(N-1)\theta]\Sigma_1'\mu_{22}\alpha(N)\\
    &+(-\Sigma'_1\mu_{22}+(\mu_{11}-\Sigma'_1)(s-1)\theta)\Sigma'_2\mu_{11}^{N-s}f(1,s)\\
     =&(\Sigma'_1\mu_{22}+\mu_{11}\Sigma_2'-\Sigma'_1\Sigma_2')\bigg( \Sigma'_1\mu_{11}^{N-1}+(s-1)\theta\Sigma'_1\mu_{11}^{N-s}\alpha(s)\bigg)\\
     &-\mu_{22}(s-1)\theta\Sigma'_1\mu_{11}^{N-s+1}\alpha(s)-\Sigma'_1\mu_{22}\Sigma'_1\mu_{11}^{N-s}[\Sigma'_2-\mu_{22}]\alpha(s)\\
     &-\Sigma_1'[\Sigma'_2+(N-1)\theta]\mu_{22}\sum_{k=s+1}^{N} \mu_{11}^{k-s-1}f(k,N)f(1,s+1)\\
     &+\mu_{11}(s-1)\theta[\Sigma'_2+(N-1)\theta]\mu_{22}\sum_{k=s+1}^{N} \mu_{11}^{k-s-1}f(k,N)f(1,s)\\
     &+\Sigma'_1\mu_{22}[\Sigma_2'+(N-1)\theta][f(1,N)+\Sigma_1'\alpha(N)]+(\mu_{11}-\Sigma'_1)[\Sigma_2'+(N-1)\theta]\Sigma_1'\mu_{22}\alpha(N)\\
    &+(-\Sigma'_1\mu_{22}+(\mu_{11}-\Sigma'_1)(s-1)\theta)\Sigma'_2\mu_{11}^{N-s}f(1,s)\\
     =&(\Sigma'_1\mu_{22}+\mu_{11}\Sigma_2'-\Sigma'_1\Sigma_2')\bigg( \Sigma'_1\mu_{11}^{N-1}+(s-1)\theta\Sigma'_1\mu_{11}^{N-s}\alpha(s)\bigg)\\
     &-\mu_{22}(s-1)\theta\Sigma'_1\mu_{11}^{N-s+1}\alpha(s)-\Sigma'_1\mu_{22}\Sigma'_1\mu_{11}^{N-s}[\Sigma'_2-\mu_{22}]\alpha(s)\\
     &-\Sigma_1'[\Sigma'_2+(N-1)\theta]\mu_{22}\sum_{k=s+1}^{N} \mu_{11}^{k-s-1}f(k,N)f(1,s+1)\\
     &+\mu_{11}(s-1)\theta[\Sigma'_2+(N-1)\theta]\mu_{22}\sum_{k=s+1}^{N} \mu_{11}^{k-s-1}f(k,N)f(1,s)\\
     &+\Sigma'_1\mu_{22}[\Sigma_2'+(N-1)\theta][f(1,N)]+\Sigma'_1\mu_{22}[\Sigma_2'+(N-1)\theta][\Sigma'_1\alpha(N)+(\mu_{11}-\Sigma'_1)\alpha(N)]\\
    &+(-\Sigma'_1\mu_{22}+(\mu_{11}-\Sigma'_1)(s-1)\theta)\Sigma'_2\mu_{11}^{N-s}f(1,s)\\
    =&(\Sigma'_1\mu_{22}+\mu_{11}\Sigma_2'-\Sigma'_1\Sigma_2')\bigg( \Sigma'_1\mu_{11}^{N-1}+(s-1)\theta\Sigma'_1\mu_{11}^{N-s}\alpha(s)\bigg)\\
     &-\mu_{22}(s-1)\theta\Sigma'_1\mu_{11}^{N-s+1}\alpha(s)-\Sigma'_1\mu_{22}\Sigma'_1\mu_{11}^{N-s}[\Sigma'_2-\mu_{22}]\alpha(s)\\
     &-\Sigma_1'[\Sigma'_2+(N-1)\theta]\mu_{22}\sum_{k=s+1}^{N} \mu_{11}^{k-s-1}f(k,N)f(1,s+1)\\
     &+\mu_{11}(s-1)\theta[\Sigma'_2+(N-1)\theta]\mu_{22}\sum_{k=s+1}^{N} \mu_{11}^{k-s-1}f(k,N)f(1,s)\\
     &+\Sigma'_1\mu_{22}[\Sigma_2'+(N-1)\theta][f(1,N)]+\Sigma'_1\mu_{22}[\Sigma_2'+(N-1)\theta]\mu_{11}\alpha(N)\\
    &+(-\Sigma'_1\mu_{22}+(\mu_{11}-\Sigma'_1)(s-1)\theta)\Sigma'_2\mu_{11}^{N-s}f(1,s)\\
    =&(\Sigma'_1\mu_{22}+\mu_{11}\Sigma_2'-\Sigma'_1\Sigma_2')\bigg( \Sigma'_1\mu_{11}^{N-1}+(s-1)\theta\Sigma'_1\mu_{11}^{N-s}\alpha(s)\bigg)\\
     &-\mu_{22}(s-1)\theta\Sigma'_1\mu_{11}^{N-s+1}\alpha(s)-\Sigma'_1\mu_{22}\Sigma'_1\mu_{11}^{N-s}[\Sigma'_2-\mu_{22}]\alpha(s)\\
     &-\Sigma_1'[\Sigma'_2+(N-1)\theta]\mu_{22}\sum_{k=s+1}^{N} \mu_{11}^{k-s-1}f(k,N)f(1,s+1)\\
     &+\mu_{11}(s-1)\theta[\Sigma'_2+(N-1)\theta]\mu_{22}\sum_{k=s+1}^{N} \mu_{11}^{k-s-1}f(k,N)f(1,s)\\
     &+\Sigma'_1\mu_{22}[\Sigma_2'+(N-1)\theta][f(1,N)]\\
     &+\Sigma'_1\mu_{22}[\Sigma_2'+(N-1)\theta]\mu_{11}\bigg[\sum_{k=2}^{N+1-s}\mu_{11}^{k-2}f(k,N)+\sum_{k=N+2-s}^N\mu_{11}^{k-2}f(k,N) \bigg]\\
    &+(-\Sigma'_1\mu_{22}+(\mu_{11}-\Sigma'_1)(s-1)\theta)\Sigma'_2\mu_{11}^{N-s}f(1,s)\\
    =&(\Sigma'_1\mu_{22}+\mu_{11}\Sigma_2'-\Sigma'_1\Sigma_2')\bigg( \Sigma'_1\mu_{11}^{N-1}+(s-1)\theta\Sigma'_1\mu_{11}^{N-s}\alpha(s)\bigg)\\
     &-\mu_{22}(s-1)\theta\Sigma'_1\mu_{11}^{N-s+1}\alpha(s)-\Sigma'_1\mu_{22}\Sigma'_1\mu_{11}^{N-s}[\Sigma'_2-\mu_{22}]\alpha(s)\\
     &-\Sigma_1'[\Sigma'_2+(N-1)\theta]\mu_{22}\sum_{k=s+1}^{N} \mu_{11}^{k-s-1}f(k,N)f(1,s+1)\\
     &+\mu_{11}(s-1)\theta[\Sigma'_2+(N-1)\theta]\mu_{22}\sum_{k=s+1}^{N} \mu_{11}^{k-s-1}f(k,N)f(1,s)\\
     &+\Sigma'_1\mu_{22}[\Sigma_2'+(N-1)\theta][f(1,N)]+\Sigma'_1\mu_{22}[\Sigma_2'+(N-1)\theta]\bigg[\sum_{k=2}^{N+1-s}\mu_{11}^{k-1}f(k,N)\bigg]\\
     &+\Sigma'_1\mu_{22}[\Sigma_2'+(N-1)\theta]\sum_{k=N+2-s}^N\mu_{11}^{k-1}f(k,N)\\
    &+(-\Sigma'_1\mu_{22}+(\mu_{11}-\Sigma'_1)(s-1)\theta)\Sigma'_2\mu_{11}^{N-s}f(1,s)
\end{align*}
Let $k'+N-s=k$
\begin{align*}
    =&(\Sigma'_1\mu_{22}+\mu_{11}\Sigma_2'-\Sigma'_1\Sigma_2')\bigg( \Sigma'_1\mu_{11}^{N-1}+(s-1)\theta\Sigma'_1\mu_{11}^{N-s}\alpha(s)\bigg)\\
     &-\mu_{22}(s-1)\theta\Sigma'_1\mu_{11}^{N-s+1}\alpha(s)-\Sigma'_1\mu_{22}\Sigma'_1\mu_{11}^{N-s}[\Sigma'_2-\mu_{22}]\alpha(s)\\
     &-\Sigma_1'[\Sigma'_2+(N-1)\theta]\mu_{22}\sum_{k=s+1}^{N} \mu_{11}^{k-s-1}f(k,N)f(1,s+1)\\
     &+\mu_{11}(s-1)\theta[\Sigma'_2+(N-1)\theta]\mu_{22}\sum_{k=s+1}^{N} \mu_{11}^{k-s-1}f(k,N)f(1,s)\\
     &+\Sigma'_1\mu_{22}[\Sigma_2'+(N-1)\theta][f(1,N)]+\Sigma'_1\mu_{22}[\Sigma_2'+(N-1)\theta]\bigg[\sum_{k=2}^{N+1-s}\mu_{11}^{k-1}f(k,N)\bigg]\\
     &+\Sigma'_1\mu_{22}[\Sigma_2'+(N-1)\theta]\sum_{k'=2}^s\mu_{11}^{k'-1+N-s}f(k'+N-s,N)\\
    &+(-\Sigma'_1\mu_{22}+(\mu_{11}-\Sigma'_1)(s-1)\theta)\Sigma'_2\mu_{11}^{N-s}f(1,s)\\
    =&(\Sigma'_1\mu_{22}+\mu_{11}\Sigma_2'-\Sigma'_1\Sigma_2')\bigg( \Sigma'_1\mu_{11}^{N-1}+(s-1)\theta\Sigma'_1\mu_{11}^{N-s}\alpha(s)\bigg)\\
     &-\mu_{22}(s-1)\theta\Sigma'_1\mu_{11}^{N-s+1}\alpha(s)-\Sigma'_1\mu_{22}\Sigma'_1\mu_{11}^{N-s}[\Sigma'_2-\mu_{22}]\alpha(s)\\
     &-\Sigma_1'[\Sigma'_2+(N-1)\theta]\mu_{22}\sum_{k=s+1}^{N} \mu_{11}^{k-s-1}f(k,N)f(1,s+1)\\
     &+\mu_{11}(s-1)\theta[\Sigma'_2+(N-1)\theta]\mu_{22}\sum_{k=s+1}^{N} \mu_{11}^{k-s-1}f(k,N)f(1,s)\\
     &+\Sigma'_1\mu_{22}[\Sigma_2'+(N-1)\theta]\bigg[\sum_{k=2}^{N+1-s}\mu_{11}^{k-1}f(k,N)+f(1,N)\bigg]\\
     &+\Sigma'_1\mu_{22}[\Sigma_2'+(N-1)\theta]\sum_{k'=2}^s\mu_{11}^{k'-1+N-s}f(k'+N-s,N)\\
    &+(-\Sigma'_1\mu_{22}+(\mu_{11}-\Sigma'_1)(s-1)\theta)\Sigma'_2\mu_{11}^{N-s}f(1,s)\\
    =&(\Sigma'_1\mu_{22}+\mu_{11}\Sigma_2'-\Sigma'_1\Sigma_2')\bigg( \Sigma'_1\mu_{11}^{N-1}+(s-1)\theta\Sigma'_1\mu_{11}^{N-s}\alpha(s)\bigg)\\
     &-\mu_{22}(s-1)\theta\Sigma'_1\mu_{11}^{N-s+1}\alpha(s)-\Sigma'_1\mu_{22}\Sigma'_1\mu_{11}^{N-s}[\Sigma'_2-\mu_{22}]\alpha(s)\\
     &-\Sigma_1'[\Sigma'_2+(N-1)\theta]\mu_{22}\sum_{k=s+1}^{N} \mu_{11}^{k-s-1}f(k,N)f(1,s+1)\\
     &+\mu_{11}(s-1)\theta[\Sigma'_2+(N-1)\theta]\mu_{22}\sum_{k=s+1}^{N} \mu_{11}^{k-s-1}f(k,N)f(1,s)\\
     &+\Sigma'_1\mu_{22}[\Sigma_2'+(N-1)\theta]\bigg[\sum_{k=1}^{N+1-s}\mu_{11}^{k-1}f(k,N)\bigg]\\
     &+\Sigma'_1\mu_{22}[\Sigma_2'+(N-1)\theta]\sum_{k'=2}^s\mu_{11}^{k'-1+N-s}f(k'+N-s,N)\\
    &+(-\Sigma'_1\mu_{22}+(\mu_{11}-\Sigma'_1)(s-1)\theta)\Sigma'_2\mu_{11}^{N-s}f(1,s)\\
    =&(\Sigma'_1\mu_{22}+\mu_{11}\Sigma_2'-\Sigma'_1\Sigma_2')\bigg( \Sigma'_1\mu_{11}^{N-1}+(s-1)\theta\Sigma'_1\mu_{11}^{N-s}\alpha(s)\bigg)\\
     &-\mu_{22}(s-1)\theta\Sigma'_1\mu_{11}^{N-s+1}\alpha(s)-\Sigma'_1\mu_{22}\Sigma'_1\mu_{11}^{N-s}[\Sigma'_2-\mu_{22}]\alpha(s)\\
     &-\Sigma_1'[\Sigma'_2+(N-1)\theta]\mu_{22}\sum_{k=s+1}^{N} \mu_{11}^{k-s-1}f(k,N)f(1,s+1)\\
     &+\mu_{11}(s-1)\theta[\Sigma'_2+(N-1)\theta]\mu_{22}\sum_{k=s+1}^{N} \mu_{11}^{k-s-1}f(k,N)f(1,s)\\
     &+\Sigma'_1\mu_{22}[\Sigma_2'+(N-1)\theta]\bigg[\sum_{k=1}^{N-s}\mu_{11}^{k-1}f(k,N)+\mu_{11}^{N-s}f(N+1-s,N)\bigg]\\
     &+\Sigma'_1\mu_{22}[\Sigma_2'+(N-1)\theta]\sum_{k'=2}^s\mu_{11}^{k'-1+N-s}f(k'+N-s,N)\\
    &+(-\Sigma'_1\mu_{22}+(\mu_{11}-\Sigma'_1)(s-1)\theta)\Sigma'_2\mu_{11}^{N-s}f(1,s)
\end{align*}
Let $k''-s=k$
\begin{align*}
     =&(\Sigma'_1\mu_{22}+\mu_{11}\Sigma_2'-\Sigma'_1\Sigma_2')\bigg( \Sigma'_1\mu_{11}^{N-1}+(s-1)\theta\Sigma'_1\mu_{11}^{N-s}\alpha(s)\bigg)\\
     &-\mu_{22}(s-1)\theta\Sigma'_1\mu_{11}^{N-s+1}\alpha(s)-\Sigma'_1\mu_{22}\Sigma'_1\mu_{11}^{N-s}[\Sigma'_2-\mu_{22}]\alpha(s)\\
     &-\Sigma_1'[\Sigma'_2+(N-1)\theta]\mu_{22}\sum_{k=s+1}^{N} \mu_{11}^{k-s-1}f(k,N)f(1,s+1)\\
     &+\mu_{11}(s-1)\theta[\Sigma'_2+(N-1)\theta]\mu_{22}\sum_{k=s+1}^{N} \mu_{11}^{k-s-1}f(k,N)f(1,s)\\
     &+\Sigma'_1\mu_{22}[\Sigma_2'+(N-1)\theta]\bigg[\sum_{k''=s+1}^{N}\mu_{11}^{k''-s-1}f(k''-s,N)+\mu_{11}^{N-s}f(N+1-s,N)\bigg]\\
     &+\Sigma'_1\mu_{22}[\Sigma_2'+(N-1)\theta]\sum_{k'=2}^s\mu_{11}^{k'-1+N-s}f(k'+N-s,N)\\
    &+(-\Sigma'_1\mu_{22}+(\mu_{11}-\Sigma'_1)(s-1)\theta)\Sigma'_2\mu_{11}^{N-s}f(1,s)\\
    =&(\Sigma'_1\mu_{22}+\mu_{11}\Sigma_2'-\Sigma'_1\Sigma_2')\bigg( \Sigma'_1\mu_{11}^{N-1}+(s-1)\theta\Sigma'_1\mu_{11}^{N-s}\alpha(s)\bigg)\\
    &+\Sigma'_1\mu_{22}[\Sigma_2'+(N-1)\theta]\bigg[\sum_{k=s+1}^{N}\mu_{11}^{k-s-1}[f(k-s,N)-f(k,n)f(1,s+1)]\bigg]\\
     &-\mu_{22}(s-1)\theta\Sigma'_1\mu_{11}^{N-s+1}\alpha(s)-\Sigma'_1\mu_{22}\Sigma'_1\mu_{11}^{N-s}[\Sigma'_2-\mu_{22}]\alpha(s)\\
     &+\mu_{11}(s-1)\theta[\Sigma'_2+(N-1)\theta]\mu_{22}\sum_{k=s+1}^{N} \mu_{11}^{k-s-1}f(k,N)f(1,s)\\
     &+\Sigma'_1\mu_{22}[\Sigma_2'+(N-1)\theta]\bigg[\mu_{11}^{N-s}f(N+1-s,N)\bigg]\\
     &+\Sigma'_1\mu_{22}[\Sigma_2'+(N-1)\theta]\sum_{k'=2}^s\mu_{11}^{k'-1+N-s}f(k'+N-s,N)\\
    &+(-\Sigma'_1\mu_{22}+(\mu_{11}-\Sigma'_1)(s-1)\theta)\Sigma'_2\mu_{11}^{N-s}f(1,s)\\
    =&(\Sigma'_1\mu_{22}+\mu_{11}\Sigma_2'-\Sigma'_1\Sigma_2')\bigg( \Sigma'_1\mu_{11}^{N-1}+(s-1)\theta\Sigma'_1\mu_{11}^{N-s}\alpha(s)\bigg)\\
    &+\Sigma'_1\mu_{22}[\Sigma_2'+(N-1)\theta]\bigg[\sum_{k=s+1}^{N}\mu_{11}^{k-s-1}[f(k-s,N)-f(k,n)f(1,s+1)]\bigg]\\
    &+\mu_{11}(s-1)\theta[\Sigma'_2+(N-1)\theta]\mu_{22}\sum_{k=s+1}^{N} \mu_{11}^{k-s-1}f(k,N)f(1,s)\\
     &-\mu_{22}(s-1)\theta\Sigma'_1\mu_{11}^{N-s+1}\alpha(s)-\Sigma'_1\mu_{22}\Sigma'_1\mu_{11}^{N-s}[\Sigma'_2-\mu_{22}]\alpha(s)\\
     &+\Sigma'_1\mu_{22}[\Sigma_2'+(N-1)\theta]\bigg[\mu_{11}^{N-s}f(N+1-s,N)\bigg]\\
     &+\Sigma'_1\mu_{22}[\Sigma_2'+(N-1)\theta]\sum_{k'=2}^s\mu_{11}^{k'-1+N-s}f(k'+N-s,N)\\
    &+(-\Sigma'_1\mu_{22}+(\mu_{11}-\Sigma'_1)(s-1)\theta)\Sigma'_2\mu_{11}^{N-s}f(1,s)\\
    =&(\Sigma'_1\mu_{22}+\mu_{11}\Sigma_2'-\Sigma'_1\Sigma_2')\bigg( \Sigma'_1\mu_{11}^{N-1}+(s-1)\theta\Sigma'_1\mu_{11}^{N-s}\alpha(s)\bigg)\\
    &+\Sigma'_1\mu_{22}[\Sigma_2'+(N-1)\theta]\bigg[\sum_{k=s+1}^{N}\mu_{11}^{k-s-1}[f(k-s,N)-f(k,n)f(1,s+1)]\bigg]\\
    &+\mu_{11}(s-1)\theta[\Sigma'_2+(N-1)\theta]\mu_{22}\sum_{k=s+1}^{N} \mu_{11}^{k-s-1}f(k,N)f(1,s)\\
     &-\mu_{22}(s-1)\theta\Sigma'_1\mu_{11}^{N-s+1}\alpha(s)-\Sigma'_1\mu_{22}\Sigma'_1\mu_{11}^{N-s}[\Sigma'_2-\mu_{22}]\alpha(s)\\
     &+\Sigma'_1\mu_{22}[\Sigma_2'+(N-1)\theta]\bigg[\mu_{11}^{N-s}f(N+1-s,N)\bigg]\\
     &+\Sigma'_1\mu_{22}\Sigma_2'\sum_{k'=2}^s\mu_{11}^{k'-1+N-s}f(k'+N-s,N)\\
     &+\Sigma'_1\mu_{22}[(N-1)\theta]\sum_{k'=2}^s\mu_{11}^{k'-1+N-s}f(k'+N-s,N)\\
    &+(-\Sigma'_1\mu_{22}+(\mu_{11}-\Sigma'_1)(s-1)\theta)\Sigma'_2\mu_{11}^{N-s}f(1,s)\\
    =&(\Sigma'_1\mu_{22}+\mu_{11}\Sigma_2'-\Sigma'_1\Sigma_2')\bigg( \Sigma'_1\mu_{11}^{N-1}+(s-1)\theta\Sigma'_1\mu_{11}^{N-s}\alpha(s)\bigg)\\
    &+\Sigma'_1\mu_{22}[\Sigma_2'+(N-1)\theta]\bigg[\sum_{k=s+1}^{N}\mu_{11}^{k-s-1}[f(k-s,N)-f(k,n)f(1,s+1)]\bigg]\\
    &+\mu_{11}(s-1)\theta[\Sigma'_2+(N-1)\theta]\mu_{22}\sum_{k=s+1}^{N} \mu_{11}^{k-s-1}f(k,N)f(1,s)\\
    &+\Sigma'_1\mu_{22}\mu_{11}^{N-s} \bigg[\Sigma_2'\sum_{k'=2}^s\mu_{11}^{k'-1}f(k'+N-s,N)-\Sigma_1' [\Sigma'_2-\mu_{22}]\alpha(s)\bigg]\\
     &-\mu_{22}(s-1)\theta\Sigma'_1\mu_{11}^{N-s+1}\alpha(s)\\
     &+\Sigma'_1\mu_{22}[\Sigma_2'+(N-1)\theta]\bigg[\mu_{11}^{N-s}f(N+1-s,N)\bigg]\\
     &+\Sigma'_1\mu_{22}[(N-1)\theta]\sum_{k'=2}^s\mu_{11}^{k'-1+N-s}f(k'+N-s,N)\\
    &+(-\Sigma'_1\mu_{22}+(\mu_{11}-\Sigma'_1)(s-1)\theta)\Sigma'_2\mu_{11}^{N-s}f(1,s)\\
    =&(\Sigma'_1\mu_{22}+\mu_{11}\Sigma_2'-\Sigma'_1\Sigma_2')\bigg( \Sigma'_1\mu_{11}^{N-1}+(s-1)\theta\Sigma'_1\mu_{11}^{N-s}\alpha(s)\bigg)\\
    &+\Sigma'_1\mu_{22}[\Sigma_2'+(N-1)\theta]\bigg[\sum_{k=s+1}^{N}\mu_{11}^{k-s-1}[f(k-s,N)-f(k,n)f(1,s+1)]\bigg]\\
    &+\mu_{11}(s-1)\theta[\Sigma'_2+(N-1)\theta]\mu_{22}\sum_{k=s+1}^{N} \mu_{11}^{k-s-1}f(k,N)f(1,s)\\
    &+\Sigma'_1\mu_{22}\mu_{11}^{N-s} \bigg[\Sigma_2'\sum_{k'=2}^s\mu_{11}^{k'-1}f(k'+N-s,N)-\Sigma_1' [\Sigma'_2-\mu_{22}]\alpha(s)\bigg]\\
    &+\Sigma'_1\mu_{22}\mu_{11}^{N-s+1} \bigg[(N-1)
    \theta\sum_{k'=2}^s\mu_{11}^{k'-2}f(k'+N-s,N)-(s-1)\theta\alpha(s)\bigg] \\
     &+\Sigma'_1\mu_{22}[\Sigma_2'+(N-1)\theta]\bigg[\mu_{11}^{N-s}f(N+1-s,N)\bigg]\\
    &+(-\Sigma'_1\mu_{22}+(\mu_{11}-\Sigma'_1)(s-1)\theta)\Sigma'_2\mu_{11}^{N-s}f(1,s)\\
    =&(\Sigma'_1\mu_{22}+\mu_{11}\Sigma_2'-\Sigma'_1\Sigma_2')\bigg( \Sigma'_1\mu_{11}^{N-1}+(s-1)\theta\Sigma'_1\mu_{11}^{N-s}\alpha(s)\bigg)\\
    &+\Sigma'_1\mu_{22}[\Sigma_2'+(N-1)\theta]\bigg[\sum_{k=s+1}^{N}\mu_{11}^{k-s-1}[f(k-s,N)-f(k,n)f(1,s+1)]\bigg]\\
    &+\mu_{11}(s-1)\theta[\Sigma'_2+(N-1)\theta]\mu_{22}\sum_{k=s+1}^{N} \mu_{11}^{k-s-1}f(k,N)f(1,s)\\
    &+\Sigma'_1\mu_{22}\mu_{11}^{N-s} \bigg[\Sigma_2'\sum_{k'=2}^s\mu_{11}^{k'-1}f(k'+N-s,N)-\Sigma_1' [\Sigma'_2-\mu_{22}]\alpha(s)\bigg]\\
    &+\Sigma'_1\mu_{22}\mu_{11}^{N-s+1} \bigg[(N-1)
    \theta\sum_{k'=2}^s\mu_{11}^{k'-2}f(k'+N-s,N)-(s-1)\theta\alpha(s)\bigg] \\
     &+\Sigma'_1\mu_{22}\Sigma_2'\bigg[\mu_{11}^{N-s}f(N+1-s,N)\bigg]+\Sigma'_1\mu_{22}[(N-1)\theta]\bigg[\mu_{11}^{N-s}f(N+1-s,N)\bigg]\\
     &-\Sigma'_1\mu_{22}\Sigma'_2\mu_{11}^{N-s}f(1,s)+(\mu_{11}-\Sigma'_1)(s-1)\theta\Sigma'_2\mu_{11}^{N-s}f(1,s)\\
    =&[\Sigma'_1\mu_{22}+\mu_{11}\Sigma_2'-\Sigma'_1\Sigma_2']\bigg( \Sigma'_1\mu_{11}^{N-1}+(s-1)\theta\Sigma'_1\mu_{11}^{N-s}\alpha(s)\bigg) \refstepcounter{equation}\tag{\theequation}\label{te1}\\
    &+\Sigma'_1\mu_{22}[\Sigma_2'+(N-1)\theta]\bigg[\sum_{k=s+1}^{N}\mu_{11}^{k-s-1}[f(k-s,N)-f(k,N)f(1,s+1)]\bigg] \refstepcounter{equation}\tag{\theequation}\label{te2}\\
    &+\mu_{11}(s-1)\theta[\Sigma'_2+(N-1)\theta]\mu_{22}\sum_{k=s+1}^{N} \mu_{11}^{k-s-1}f(k,N)f(1,s) \refstepcounter{equation}\tag{\theequation}\label{te3}\\
    &+\Sigma'_1\mu_{22}\mu_{11}^{N-s} \bigg[\Sigma_2'\sum_{k=2}^s\mu_{11}^{k-1}f(k+N-s,N)-\Sigma_1' [\Sigma'_2-\mu_{22}]\alpha(s)\bigg] \refstepcounter{equation}\tag{\theequation}\label{te4}\\
    &+\Sigma'_1\mu_{22}\mu_{11}^{N-s+1} \bigg[(N-1)
    \theta\sum_{k=2}^s\mu_{11}^{k-2}f(k+N-s,N)-(s-1)\theta\alpha(s)\bigg] \refstepcounter{equation}\tag{\theequation}\label{te5}\\
     &+\Sigma'_1\mu_{22}\Sigma_2'\mu_{11}^{N-s}\bigg[f(N+1-s,N)-f(1,s)\bigg] \refstepcounter{equation}\tag{\theequation}\label{te6}\\
     &+\mu_{11}^{N-s}\bigg[\Sigma'_1\mu_{22}(N-1)\theta f(N+1-s,N)+(\mu_{11}-\Sigma'_1)(s-1)\theta\Sigma'_2f(1,s)\bigg]. \refstepcounter{equation}\tag{\theequation}\label{te7}
\end{align*}
We show that all the terms \eqref{te1}-\eqref{te7} are non-negative.

The non-negativity of 
 \eqref{te1} follows from
 \begin{align}\label{N geq 2 implication}
     &\Sigma'_1\mu_{22}+\mu_{11}\Sigma_2'-\Sigma'_1\Sigma_2' \nonumber\\
     =&(1+\lambda)(\mu_{11}+\mu_{21})\mu_{22}+(1+\lambda)(\mu_{12}+\mu_{22})\mu_{11}-(1+\lambda)^2\Sigma_1 \Sigma_2 \nonumber\\
     =&(1+\lambda)\bigg[\mu_{11}\mu_{22}+\mu_{21}\mu_{22}+\mu_{12}\mu_{11}+\mu_{22}\mu_{11}- (\mu_{11}+\mu_{21})(\mu_{22}+\mu_{12})-\lambda \Sigma_1 \Sigma_2\bigg] \nonumber\\
     =&(1+\lambda)\bigg[\mu_{11}\mu_{22}+\mu_{21}\mu_{22}+\mu_{12}\mu_{11}+\mu_{22}\mu_{11}-\mu_{11}\mu_{22}-\mu_{11}\mu_{12}-\mu_{21}\mu_{22}-\mu_{21}\mu_{12}-\lambda\Sigma_1 \Sigma_2\bigg] \nonumber\\
     =&(1+\lambda)\bigg[\mu_{11}\mu_{22}-\mu_{21}\mu_{12}-\lambda\Sigma_1 \Sigma_2\bigg] \geq 0.
 \end{align}
from \eqref{lambda} since $\tau'(2) \geq 0$ for this case to be possible.

The non-negativity of  \eqref{te2} follows since $s+1 \leq k \leq N$ and 
 \begin{align*}
     f(k-s,N)-f(1,s+1)f(k,N)=f(k-s,k)f(k,N)-f(1,s+1)f(k,N)\geq 0
 \end{align*}
 as $f(k-s,k) \geq f(1,s+1)$ for $k \geq s+1$.

Clearly, \eqref{te3} is non-negative.

For \eqref{te4} we have
\begin{align*}
    \Sigma_2'\sum_{k=2}^s\mu_{11}^{k-1}f(k+N-s,N)-\Sigma_1' [\Sigma'_2-\mu_{22}]\alpha(s)
    \geq & \Sigma_2'\mu_{11}\sum_{k=2}^s\mu_{11}^{k-2}f(k,s)-\Sigma_1' [\Sigma'_2-\mu_{22}]\alpha(s)\\
    =&\Sigma_2'\mu_{11}\alpha(s)-\Sigma_1' [\Sigma'_2-\mu_{22}]\alpha(s)\\
    =&[\Sigma'_1\mu_{22}+\mu_{11}\Sigma'_2-\Sigma'_1\Sigma'_2]\alpha(s) \geq 0,
\end{align*}
 where the last inequality follows from \eqref{N geq 2 implication}. Therefore \eqref{te4} is non-negative.

 For \eqref{te5}, $s \leq N-1$ yields
 \begin{align*}
     &(N-1)
    \theta\sum_{k=2}^s\mu_{11}^{k-2}f(k+N-s,N)-(s-1)\theta\alpha(s)\\
    \geq &(s-1)\theta \bigg[\sum_{k=2}^s\mu_{11}^{k-2}f(k+N-s,N)-\alpha(s) \bigg]\geq 0,
 \end{align*}
where the last inequality follows from \eqref{upper bound alpha}.

The non-negativity of  \eqref{te6} follows as $f(N+1-s,N)\geq f(1,s)$ for $s=1,\ldots,N-1$.

 Finally, for \eqref{te7}, $s \leq N-1$ yields
 \begin{align*}
     &\Sigma'_1\mu_{22}(N-1)\theta f(N+1-s,N)+(\mu_{11}-\Sigma'_1)(s-1)\theta\Sigma'_2f(1,s)\\
      \geq &\Sigma'_1\mu_{22}(s-1)\theta f(1,s)+(\mu_{11}-\Sigma'_1)(s-1)\theta\Sigma'_2f(1,s)\\
      =&(s-1)\theta f(1,s) \bigg[ \Sigma'_1\mu_{22}+(\mu_{11}-\Sigma'_1)\Sigma'_2 \bigg]\\
      =&(s-1)\theta f(1,s) \bigg[ \Sigma'_1\mu_{22}+\mu_{11}\Sigma'_2-\Sigma'_1\Sigma'_2\bigg] \geq 0,
 \end{align*}
 where the last inequality again follows from \eqref{N geq 2 implication}.
 Hence we have shown $\Gamma''(s,a_{11}) \geq 0$ and therefore $\Gamma(s,a_{11}) \geq 0$.

\item[ Case 3:] $s=N, \ldots, B+2$. 
 
 In this case we have $d_N(s)=a_{22}$.
    \item[\ Case 3, Part (i) :]  $a=a_{12}$. 
    
    We consider the two sub-cases $(a)\ s=N. \ldots,B+1 $ and $(b)\  s=B+2$ separately. 
         
  \item[$(a)$] $s=N. \ldots,B+1$. We have
    \begin{align*}
        \Gamma(s,a_{12})=&\Sigma'_2+\frac{\Sigma'_2+(s-1)\theta}{q}h_N(s-1)+\bigg(1-\frac{\Sigma_2'+(s-1)\theta}{q} \bigg)h_N(s)\\
        &-\mu_{22}-\frac{\mu_{11}}{q}h_N(s+1)-\frac{\mu_{22}+(s-1)\theta}{q}h_N(s-1)-\bigg(1-\frac{\mu_{11}+\mu_{22}+(s-1)\theta}{q}\bigg)h_N(s).
    \end{align*}
    That is 
    \begin{align}\label{gamma(s,a_12).2}
          \Gamma(s,a_{12})=& \Sigma'_2-\mu_{22}-\frac{\mu_{11}}{q} h_N(s+1)+\frac{\Sigma_2'-\mu_{22}}{q}h_N(s-1)+\frac{\mu_{11}}{q}h_N(s)-\frac{\Sigma'_2-\mu_{22}}{q}h_N(s) \nonumber\\
          =&\Sigma'_2-\mu_{22}+\frac{\Sigma'_2-\mu_{22}}{q}[h_N(s-1)-h_N(s)]+\frac{\mu_{11}}{q}[h_N(s)-h_N(s+1)].
    \end{align}
    It follows from Lemma \ref{lemma h_N1} that 
    \begin{align}\label{gamma(s,a_12)}
         \Gamma(s,a_{12})=&\Sigma'_2-\mu_{22}+\frac{\Sigma'_2-\mu_{22}}{\Sigma'_2+(s-1)\theta}[g_N-\Sigma_2']+\frac{\mu_{11}}{\Sigma'_2+s\theta}[g_N-\Sigma_2'].
    \end{align}
    In particular for $s=N$,
    \begin{align*}
         \Gamma(N,a_{12})=&\Sigma'_2-\mu_{22}+\frac{\Sigma'_2-\mu_{22}}{\Sigma'_2+(N-1)\theta}[g_N-\Sigma_2']+\frac{\mu_{11}}{\Sigma'_2+N\theta}[g_N-\Sigma_2'].
    \end{align*}
    Using Proposition \ref{prop gain}, we have 
\begin{align}\label{g_N-sigma_2}
    g_N-\Sigma'_2
    =&\frac{-[\Sigma'_2+(N-1)\theta]\bigg(\Sigma'_2f(1,N)+\Sigma'_1[\Sigma'_2-\mu_{22}]\alpha(N)\bigg)}{\Delta(N)}.
\end{align}
Hence
\begin{align*}
     \Gamma(N,a_{12})
     =&\Sigma'_2-\mu_{22}-\frac{\Sigma'_2-\mu_{22}}{\Delta(N)}\bigg(\Sigma'_2f(1,N)+\Sigma'_1[\Sigma'_2-\mu_{22}]\alpha(N)\bigg)\\
     &-\frac{\mu_{11}}{\Sigma'_2+N\theta}\Bigg[\frac{[\Sigma'_2+(N-1)\theta]\bigg(\Sigma'_2f(1,N)+\Sigma'_1[\Sigma'_2-\mu_{22}]\alpha(N)\bigg)}{\Delta(N)} \Bigg].
\end{align*}
That is 
\begin{align*}
   & [\Sigma'_2+N\theta]\Delta(N) \Gamma(N,a_{12})\\
   =&[\Sigma'_2+N\theta]\Bigg[(\Sigma'_2-\mu_{22})\Delta(N)-(\Sigma'_2-\mu_{22})\bigg(\Sigma'_2f(1,N)+\Sigma'_1[\Sigma'_2-\mu_{22}]\alpha(N)\bigg) \Bigg]\\
     &-\mu_{11}[\Sigma'_2+(N-1)\theta]\bigg(\Sigma'_2f(1,N)+\Sigma'_1[\Sigma'_2-\mu_{22}]\alpha(N)\bigg).
\end{align*}
Recalling \eqref{Delta(n)}, we have
\begin{align*}
& [\Sigma'_2+N\theta]\Delta(N) \Gamma(N,a_{12})\\
    =&[\Sigma'_2+N\theta]\Bigg[(\Sigma'_2-\mu_{22})[\Sigma'_2+(N-1)\theta][f(1,N)+\Sigma'_1\alpha(N)]+(\Sigma'_2-\mu_{22})\Sigma'_1\mu_{11}^{N-1}\\
    &-(\Sigma'_2-\mu_{22})\bigg(\Sigma'_2f(1,N)+\Sigma'_1[\Sigma'_2-\mu_{22}]\alpha(N)\bigg) \Bigg]\\
     &-\mu_{11}[\Sigma'_2+(N-1)\theta]\bigg(\Sigma'_2f(1,N)+\Sigma'_1[\Sigma'_2-\mu_{22}]\alpha(N)\bigg)\\
     =&[\Sigma'_2+N\theta]\Bigg[(\Sigma'_2-\mu_{22})[\Sigma'_2+(N-1)\theta]f(1,N)+(\Sigma'_2-\mu_{22})[(\Sigma'_2-\mu_{22})+\mu_{22}+(N-1)\theta]\Sigma'_1\alpha(N)\\
     &+(\Sigma'_2-\mu_{22})\Sigma'_1\mu_{11}^{N-1}-(\Sigma'_2-\mu_{22})\bigg(\Sigma'_2f(1,N)+\Sigma'_1[\Sigma'_2-\mu_{22}]\alpha(N)\bigg) \Bigg]\\
     &-\mu_{11}[\Sigma'_2+(N-1)\theta]\bigg(\Sigma'_2f(1,N)+\Sigma'_1[\Sigma'_2-\mu_{22}]\alpha(N)\bigg)\\
     =&[\Sigma'_2+N\theta]\Bigg[(\Sigma'_2-\mu_{22})[\Sigma'_2+(N-1)\theta]f(1,N)+(\Sigma'_2-\mu_{22})^2\Sigma'_1\alpha(N)\\
     &+(\Sigma'_2-\mu_{22})\Sigma'_1\bigg((\mu_{22}+(N-1)\theta)\alpha(N)+\mu_{11}^{N-1}\bigg)-(\Sigma'_2-\mu_{22})\Sigma'_2f(1,N)-\Sigma'_1(\Sigma'_2-\mu_{22})^2\alpha(N) \Bigg]\\
     &-\mu_{11}[\Sigma'_2+(N-1)\theta]\bigg(\Sigma'_2f(1,N)+\Sigma'_1[\Sigma'_2-\mu_{22}]\alpha(N)\bigg)\\
    \stackrel{(*)}{=}&[\Sigma'_2+N\theta]\Bigg[(\Sigma'_2-\mu_{22})[\Sigma'_2+(N-1)\theta]f(1,N)+(\Sigma'_2-\mu_{22})\Sigma'_1\alpha(N+1)-(\Sigma'_2-\mu_{22})\Sigma'_2f(1,N) \Bigg]\\
     &-\mu_{11}[\Sigma'_2+(N-1)\theta]\bigg(\Sigma'_2f(1,N)+\Sigma'_1[\Sigma'_2-\mu_{22}]\alpha(N)\bigg)\\
     =&[\Sigma'_2+N\theta]\Bigg[-\mu_{22}[\Sigma'_2+(N-1)\theta]f(1,N)+\Sigma'_2[\Sigma'_2+(N-1)\theta]f(1,N)-(\Sigma'_2-\mu_{22})\Sigma'_2f(1,N)\\
     &+(\Sigma'_2-\mu_{22})\Sigma'_1\alpha(N+1) \Bigg]-\mu_{11}[\Sigma'_2+(N-1)\theta]\bigg(\Sigma'_2f(1,N)+\Sigma'_1[\Sigma'_2-\mu_{22}]\alpha(N)\bigg)\\
      =&[\Sigma'_2+N\theta]\Bigg[-\mu_{22}[\Sigma'_2+(N-1)\theta]f(1,N)+\Sigma'_2[\mu_{22}+(N-1)\theta]f(1,N)+(\Sigma'_2-\mu_{22})\Sigma'_1\alpha(N+1) \Bigg]\\
     &-\mu_{11}[\Sigma'_2+(N-1)\theta]\bigg(\Sigma'_2f(1,N)+\Sigma'_1[\Sigma'_2-\mu_{22}]\alpha(N)\bigg)\\
     =&[\Sigma'_2+N\theta]\Bigg[-\mu_{22}[\Sigma'_2+(N-1)\theta]f(1,N)+\Sigma'_2f(1,N+1)+(\Sigma'_2-\mu_{22})\Sigma'_1\alpha(N+1) \Bigg]\\
     &-\mu_{11}[\Sigma'_2+(N-1)\theta]\bigg(\Sigma'_2f(1,N)+\Sigma'_1[\Sigma'_2-\mu_{22}]\alpha(N)\bigg)\\
     =&-[\Sigma'_2+N\theta][\Sigma'_2+(N-1)\theta]\mu_{22}f(1,N)+[\Sigma'_2+N\theta][\Sigma'_1\alpha(N+1)(\Sigma'_2-\mu_{22})+\Sigma'_2f(1,N+1)]\\
     &-[\Sigma'_2+(N-1)\theta]\mu_{11}[\Sigma'_2f(1,N)+\Sigma'_1\alpha(N)(\Sigma'_2-\mu_{22})]\\
     =&-\tau'(N+1),
\end{align*}
where $(*)$ follows from \eqref{eq5} and the last equality follows from \eqref{tau'(n+1)}.
That is we have shown that 
\begin{align*}
     [\Sigma'_2+N\theta]\Delta(N) \Gamma(N,a_{12})=-\tau'(N+1) \geq 0,
\end{align*}
where the last inequality follows from the definition of $N$.

Now for $s=N, \ldots, B$, \eqref{gamma(s,a_12)} yields
\begin{align*}
   & \Gamma(s+1,a_{12})-\Gamma(s,a_{12})\\
   =& \Sigma'_2-\mu_{22}+\frac{\Sigma'_2-\mu_{22}}{\Sigma'_2+s\theta}[g_N-\Sigma_2']+\frac{\mu_{11}}{\Sigma'_2+(s+1)\theta}[g_N-\Sigma_2']\\
   &-(\Sigma'_2-\mu_{22})-\frac{\Sigma'_2-\mu_{22}}{\Sigma'_2+(s-1)\theta}[g_N-\Sigma_2']-\frac{\mu_{11}}{\Sigma'_2+s\theta}[g_N-\Sigma_2']\\
   =&(g_N-\Sigma'_2)\Bigg[(\Sigma_2'-\mu_{22})\bigg(\frac{1}{\Sigma'_2+s\theta}- \frac{1}{\Sigma'_2+(s-1)\theta}\bigg)+\mu_{11}\bigg(\frac{1}{\Sigma'_2+(s+1)\theta}-\frac{1}{\Sigma'_2+s\theta}\bigg) \Bigg] \geq 0,
\end{align*}
where the last inequality follows as $(g_N-\Sigma'_2)<0$ from \eqref{g_N-sigma_2}. This shows $\Gamma(s,a_{12}) \geq 0$ for for $s=N. \ldots,B+1$. 

\item[ $(b)$]$s=B+2$. We have

\begin{align*}
    \Gamma(B+2,a_{12})=&\Sigma'_2+\frac{\Sigma'_2+(B+1)\theta}{q}h_N(B+1)+\bigg(1-\frac{\Sigma'_2+(B+1)\theta}{q}\bigg)h_N(B+2)\\
    &-\mu_{22}-\frac{\mu_{22}+(B+1)\theta}{q}h_N(B+1)-\bigg(1-\frac{\mu_{22}+(B+1)\theta}{q}\bigg)h_N(B+2)\\
    =&\Sigma'_2-\mu_{22}+\frac{\Sigma'_2-\mu_{22}}{q}[h_N(B+1)-h_N(B+2)]\\
    =&\Sigma'_2-\mu_{22}+\frac{\Sigma'_2-\mu_{22}}{\Sigma'_2+(B+1)\theta}(g_N-\Sigma'_2),
\end{align*}
where the last equality follows from Lemma \ref{lemma h_N1}. 
Moreover, from \eqref{gamma(s,a_12)} we have
\begin{align*}
     \Gamma(B+1,a_{12})=&\Sigma'_2-\mu_{22}+\frac{\Sigma'_2-\mu_{22}}{\Sigma'_2+B\theta}[g_N-\Sigma_2']+\frac{\mu_{11}}{\Sigma'_2+(B+1)\theta}[g_N-\Sigma_2'].
\end{align*}
Therefore,
\begin{align*}
    &\Gamma(B+2,a_{12})-\Gamma(B+1,a_{12})\\
    &=(g_N-\Sigma'_2)\Bigg[\frac{\Sigma'_2-\mu_{22}}{\Sigma'_2+(B+1)\theta}-\frac{\mu_{11}}{\Sigma'_2+(B+1)\theta}-\frac{\Sigma'_2-\mu_{22}}{\Sigma'_2+B\theta}\Bigg]\\
    &=(g_N-\Sigma'_2)\Bigg[\frac{\Sigma'_2-\mu_{22}-\mu_{11}}{\Sigma'_2+(B+1)\theta}-\frac{\Sigma'_2-\mu_{22}}{\Sigma'_2+B\theta}\Bigg] \geq 0.
\end{align*}
The inequality follows as $\Sigma'_2-\mu_{22}-\mu_{11} \leq \Sigma'_2-\mu_{22}$, $\Sigma'_2+(B+1)\theta \geq \Sigma'_2+B\theta$ and $g_N-\Sigma'_2 \leq 0$ from \eqref{g_N-sigma_2}. Therefore $\Gamma(s,a_{12}) \geq 0$ for for $s=N. \ldots,B+2$.
\item[Case 3, Part (ii) :] For $a=a_{21}$. 

Again, we consider the two sub-cases $(a)\ s=N. \ldots,B+1 $ and $(b)\  s=B+2$ separately. 

\item[$(a)$]$s=N, \ldots, B+1$. Using an argument similar to that in \eqref{gamma(s,a_12)} (replace $\mu_{11}$ with $\mu_{21}$ and $\mu_{22}$ with $\mu_{12}$) we obtain
\begin{align*}
    \Gamma(s,a_{21})
    =&\Sigma'_2-\mu_{12}+\frac{\Sigma'_2-\mu_{12}}{\Sigma_2'+(s-1)\theta}(g_N-\Sigma'_2)+\frac{\mu_{21}}{\Sigma'_2+s\theta}(g_N-\Sigma'_2).
\end{align*}
Letting $s=N$ and using \eqref{g_N-sigma_2} we obtain
\begin{align*}
     \Gamma(N,a_{21})
     =&\Sigma'_2-\mu_{12}-\frac{\Sigma'_2-\mu_{12}}{\Delta(N)}\bigg(\Sigma'_2f(1,N)+\Sigma'_1[\Sigma'_2-\mu_{22}]\alpha(N)\bigg)\\
     &-\frac{\mu_{21}}{\Sigma'_2+N\theta}\frac{[\Sigma'_2+(N-1)\theta]\bigg(\Sigma'_2f(1,N)+\Sigma'_1[\Sigma'_2-\mu_{22}]\alpha(N)\bigg)}{\Delta(N)}.
\end{align*}
That is 
\begin{align*}
  &\Gamma(N,a_{21})\Delta(N)(\Sigma'_2+N\theta) \nonumber \\
  =&(\Sigma'_2+N\theta)\Delta(N)(\Sigma'_2-\mu_{12})-(\Sigma'_2+N\theta)\bigg(\Sigma'_2f(1,N)+\Sigma'_1[\Sigma'_2-\mu_{22}]\alpha(N)\bigg)(\Sigma'_2-\mu_{12}) \nonumber\\
     &-\mu_{21}[\Sigma'_2+(N-1)\theta]\bigg(\Sigma'_2f(1,N)+\Sigma'_1[\Sigma'_2-\mu_{22}]\alpha(N)\bigg).
\end{align*}
We consider the three terms in the above separately.
The first term satisfies
\begin{align}\label{Gamma(N,a21)Delta(N)(Sigma'_2+N theta) eq1}
    &(\Sigma'_2+N\theta)\Delta(N)(\Sigma'_2-\mu_{12}) \nonumber\\
    =&(\Sigma'_2-\mu_{12})\bigg[ [\Sigma'_2+(N-1)\theta][f(1,N)+\Sigma'_1\alpha(N)]+\Sigma'_1\mu_{11}^{N-1}\bigg](\Sigma'_2+N\theta) \nonumber\\
    =&(\Sigma'_2-\mu_{12})\Sigma'_2f(1,N)(\Sigma'_2+N\theta)+\Sigma_1'(\Sigma'_2-\mu_{12})\Sigma_2' \alpha(N)(\Sigma'_2+N\theta)+\Sigma'_1\mu_{11}^{N-1}(\Sigma'_2-\mu_{12})(\Sigma'_2+N\theta) \nonumber\\
    &+(\Sigma'_2+N\theta)(\Sigma'_2-\mu_{12})(N-1)\theta f(1,N)+\Sigma_1'(\Sigma'_2-\mu_{12})(N-1)\theta \alpha(N)(\Sigma'_2+N\theta) \nonumber\\
    =&(\Sigma'_2-\mu_{12})\Sigma'_2f(1,N)(\Sigma'_2+N\theta)+\Sigma_1'(\Sigma'_2-\mu_{12})\Sigma_2' \alpha(N)(\Sigma'_2+N\theta)+\Sigma'_1\mu_{11}^{N-1}(\Sigma'_2-\mu_{12})(\Sigma'_2+N\theta) \nonumber\\
    &+(\Sigma'_2+N\theta)(\Sigma'_2-\mu_{12})(N-1)\theta f(1,N)+\Sigma_1'(\Sigma'_2-\mu_{12})(N-2)\theta \alpha(N)(\Sigma'_2+N\theta) \nonumber\\
    &+\Sigma_1'(\Sigma'_2-\mu_{12})\theta \alpha(N)(\Sigma'_2+N\theta).
\end{align}
The second term satisfies
\begin{align}\label{Gamma(N,a21)Delta(N)(Sigma'_2+N theta) eq2}
   &-(\Sigma'_2+N\theta)\bigg(\Sigma'_2f(1,N)+\Sigma'_1[\Sigma'_2-\mu_{22}]\alpha(N)\bigg)(\Sigma'_2-\mu_{12}) \nonumber\\
   =&(\Sigma'_2+N\theta)\bigg(-\Sigma'_2f(1,N)+\Sigma'_1[\mu_{22}-\Sigma_2']\alpha(N)\bigg)(\Sigma'_2-\mu_{12}) \nonumber\\
   =&-(\Sigma'_2-\mu_{12})\Sigma'_2f(1,N)(\Sigma'_2+N\theta)+\Sigma'_1\mu_{22}(\Sigma'_2-\mu_{12})\alpha(N)(\Sigma'_2+N\theta)-\Sigma_1'\Sigma_2'(\Sigma'_2-\mu_{12})\alpha(N)(\Sigma'_2+N\theta)\nonumber\\
   =&-(\Sigma'_2-\mu_{12})\Sigma'_2f(1,N)(\Sigma'_2+N\theta)-\Sigma_1'\Sigma_2'(\Sigma'_2-\mu_{12})\alpha(N)(\Sigma'_2+N\theta) \nonumber\\
   &+\Sigma'_1\mu_{22}(\Sigma'_2-\mu_{12})(\Sigma'_2+N\theta)\bigg[\sum_{k=3}^N \mu_{11}^{k-2}f(k,N)+f(2,N) \bigg] \nonumber\\
   =&-(\Sigma'_2-\mu_{12})\Sigma'_2f(1,N)(\Sigma'_2+N\theta)-\Sigma_1'\Sigma_2'(\Sigma'_2-\mu_{12})\alpha(N)(\Sigma'_2+N\theta) \nonumber\\
   &+\Sigma'_1\mu_{22}(\Sigma'_2-\mu_{12})(\Sigma'_2+N\theta)\sum_{k=3}^N \mu_{11}^{k-2}f(k,N)+\Sigma'_1(\Sigma'_2-\mu_{12})(\Sigma'_2+N\theta)f(1,N) \nonumber\\
    =&-(\Sigma'_2-\mu_{12})\Sigma'_2f(1,N)(\Sigma'_2+N\theta)-\Sigma_1'\Sigma_2'(\Sigma'_2-\mu_{12})\alpha(N)(\Sigma'_2+N\theta) \nonumber\\
   &+\Sigma'_1\mu_{22}(\Sigma'_2-\mu_{12})(\Sigma'_2+N\theta)\sum_{k=3}^N \mu_{11}^{k-2}f(k,N)+\Sigma'_1(\Sigma'_2-\mu_{12})(\Sigma'_2+(N-1)\theta)f(1,N) \nonumber\\
   &+\Sigma'_1(\Sigma'_2-\mu_{12})\theta f(1,N).
\end{align}
And the third term satisfies 
\begin{align}\label{Gamma(N,a21)Delta(N)(Sigma'_2+N theta) eq3}
    &-\mu_{21}[\Sigma'_2+(N-1)\theta]\bigg(\Sigma'_2f(1,N)+\Sigma'_1[\Sigma'_2-\mu_{22}]\alpha(N)\bigg)  \nonumber \\
    =&-\mu_{21}\Sigma_2'[\Sigma'_2+(N-1)\theta]f(1,N)-\mu_{21}\Sigma'_1[\Sigma'_2-\mu_{22}][\Sigma'_2+(N-1)\theta]\alpha(N) \nonumber\\
    =&-\mu_{21}\Sigma_2'[\Sigma'_2+(N-1)\theta]f(1,N)+\bigg( \Sigma_1'(\Sigma_2'-\mu_{12})\mu_{11}-\mu_{21}\Sigma'_1[\Sigma'_2-\mu_{22}]\bigg)[\Sigma'_2+(N-1)\theta]\alpha(N) \nonumber\\
    &-\Sigma_1'(\Sigma_2'-\mu_{12})\mu_{11}[\Sigma'_2+(N-1)\theta]\alpha(N).
\end{align}
Adding \eqref{Gamma(N,a21)Delta(N)(Sigma'_2+N theta) eq1}-\eqref{Gamma(N,a21)Delta(N)(Sigma'_2+N theta) eq3} yields
\begin{align*}
  &\Gamma(N,a_{21})\Delta(N)(\Sigma'_2+N\theta)\\
  =&(\Sigma'_2-\mu_{12})\Sigma'_2f(1,N)(\Sigma'_2+N\theta)+\Sigma_1'(\Sigma'_2-\mu_{12})\Sigma_2' \alpha(N)(\Sigma'_2+N\theta)+\Sigma'_1\mu_{11}^{N-1}(\Sigma'_2-\mu_{12})(\Sigma'_2+N\theta)\\
    &+(\Sigma'_2+N\theta)(\Sigma'_2-\mu_{12})(N-1)\theta f(1,N)+\Sigma_1'(\Sigma'_2-\mu_{12})(N-2)\theta \alpha(N)(\Sigma'_2+N\theta)\\
    &+\Sigma_1'(\Sigma'_2-\mu_{12})\theta \alpha(N)(\Sigma'_2+N\theta)\\
    &-(\Sigma'_2-\mu_{12})\Sigma'_2f(1,N)(\Sigma'_2+N\theta)-\Sigma_1'\Sigma_2'(\Sigma'_2-\mu_{12})\alpha(N)(\Sigma'_2+N\theta)\\
   &+\Sigma'_1\mu_{22}(\Sigma'_2-\mu_{12})(\Sigma'_2+N\theta)\sum_{k=3}^N \mu_{11}^{k-2}f(k,N)+\Sigma'_1(\Sigma'_2-\mu_{12})(\Sigma'_2+(N-1)\theta)f(1,N)\\
   &+\Sigma'_1(\Sigma'_2-\mu_{12})\theta f(1,N)\\
   &-\mu_{21}\Sigma_2'[\Sigma'_2+(N-1)\theta]f(1,N)+\bigg( \Sigma_1'(\Sigma_2'-\mu_{12})\mu_{11}-\mu_{21}\Sigma'_1[\Sigma'_2-\mu_{22}]\bigg)[\Sigma'_2+(N-1)\theta]\alpha(N)\\
   &-\Sigma_1'(\Sigma_2'-\mu_{12})\mu_{11}[\Sigma'_2+(N-1)\theta]\alpha(N)\\
   =&\Sigma'_1\mu_{11}^{N-1}(\Sigma'_2-\mu_{12})(\Sigma'_2+N\theta)\\
    &+(\Sigma'_2+N\theta)(\Sigma'_2-\mu_{12})(N-1)\theta f(1,N)+\Sigma_1'(\Sigma'_2-\mu_{12})(N-2)\theta \alpha(N)(\Sigma'_2+N\theta)\\
    &+\Sigma_1'(\Sigma'_2-\mu_{12})\theta \alpha(N)(\Sigma'_2+N\theta)\\
   &+\Sigma'_1\mu_{22}(\Sigma'_2-\mu_{12})(\Sigma'_2+N\theta)\sum_{k=3}^N \mu_{11}^{k-2}f(k,N)+\Sigma'_1(\Sigma'_2-\mu_{12})(\Sigma'_2+(N-1)\theta)f(1,N)\\
   &+\Sigma'_1(\Sigma'_2-\mu_{12})\theta f(1,N)\\
   &-\mu_{21}\Sigma_2'[\Sigma'_2+(N-1)\theta]f(1,N)+\bigg( \Sigma_1'(\Sigma_2'-\mu_{12})\mu_{11}-\mu_{21}\Sigma'_1[\Sigma'_2-\mu_{22}]\bigg)[\Sigma'_2+(N-1)\theta]\alpha(N)\\
   &-\Sigma_1'(\Sigma_2'-\mu_{12})\mu_{11}[\Sigma'_2+(N-1)\theta]\alpha(N)\\
   =&(\Sigma'_2+N\theta)(\Sigma'_2-\mu_{12})(N-1)\theta f(1,N)+\Sigma'_1(\Sigma'_2-\mu_{12})\theta f(1,N)+\Sigma'_1(\Sigma'_2-\mu_{12})(\Sigma'_2+(N-1)\theta)f(1,N)\\
   &-\mu_{21}\Sigma_2'[\Sigma'_2+(N-1)\theta]f(1,N)+\bigg( \Sigma_1'(\Sigma_2'-\mu_{12})\mu_{11}-\mu_{21}\Sigma'_1[\Sigma'_2-\mu_{22}]\bigg)[\Sigma'_2+(N-1)\theta]\alpha(N)\\
   &+\Sigma_1'(\Sigma'_2-\mu_{12})(N-2)\theta \alpha(N)(\Sigma'_2+N\theta)+\Sigma'_1\mu_{11}^{N-1}(\Sigma'_2-\mu_{12})(\Sigma'_2+N\theta)\\
   &+\Sigma_1'(\Sigma'_2-\mu_{12})\theta \alpha(N)(\Sigma'_2+N\theta)+\Sigma'_1\mu_{22}(\Sigma'_2-\mu_{12})(\Sigma'_2+N\theta)\sum_{k=3}^N \mu_{11}^{k-2}f(k,N)\\
   &-\Sigma_1'(\Sigma_2'-\mu_{12})\mu_{11}[\Sigma'_2+(N-1)\theta]\alpha(N)\\
    =&(\Sigma'_2+N\theta)(\Sigma'_2-\mu_{12})(N-1)\theta f(1,N)+\Sigma'_1(\Sigma'_2-\mu_{12})\theta f(1,N) \refstepcounter{equation}\tag{\theequation}\label{ter1}\\
    &+\bigg[\Sigma'_1(\Sigma'_2-\mu_{12})-\mu_{21}\Sigma_2' \bigg][\Sigma'_2+(N-1)\theta]f(1,N) \refstepcounter{equation}\tag{\theequation}\label{ter2}\\
    &+\Sigma_1'\bigg( (\Sigma_2'-\mu_{12})\mu_{11}-\mu_{21}[\Sigma'_2-\mu_{22}]\bigg)[\Sigma'_2+(N-1)\theta]\alpha(N) \refstepcounter{equation}\tag{\theequation}\label{ter3}\\
    &+\Sigma_1'(\Sigma'_2-\mu_{12})\bigg[(N-2)\theta \alpha(N)(\Sigma'_2+N\theta)+\mu_{11}^{N-1}(\Sigma'_2+N\theta) +\theta \alpha(N)(\Sigma'_2+N\theta)\bigg] \refstepcounter{equation}\tag{\theequation}\label{ter4}\\
   &+\Sigma_1'(\Sigma'_2-\mu_{12})\bigg[\mu_{22}(\Sigma'_2+N\theta)\sum_{k=3}^N \mu_{11}^{k-2}f(k,N)-\mu_{11}[\Sigma'_2+(N-1)\theta]\alpha(N)\bigg]. \refstepcounter{equation}\tag{\theequation}\label{ter5} 
\end{align*}

We show that each of the term \eqref{ter1}-\eqref{ter5} is non-negative. \eqref{ter1} is clearly positive.

For \eqref{ter2} we have 
\begin{align*}
    &\Sigma'_1(\Sigma'_2-\mu_{12})-\mu_{21}\Sigma_2'\\
    =&(1+\lambda)(\mu_{11}+\mu_{21})(\lambda \mu_{12}+(1+\lambda)\mu_{22})-\mu_{21}(1+\lambda)(\mu_{22}+\mu_{12})\\
    =&(1+\lambda)\bigg[\lambda \mu_{11}\mu_{12}+\lambda \mu_{21}\mu_{12}+(1+\lambda)\mu_{11}\mu_{22}+(1+\lambda)\mu_{21}\mu_{22}-\mu_{21}\mu_{22}-\mu_{21}\mu_{12} \bigg]\\
     =&(1+\lambda)\bigg[\lambda \Sigma_1 \Sigma_2+(\mu_{11}\mu_{22}-\mu_{21}\mu_{12}) \bigg]\\
     \geq & 0,
\end{align*}
where the inequality follows as $\mu_{11}\mu_{22}-\mu_{21}\mu_{12} \geq 0$.

For \eqref{ter3}, note that $\alpha(1)=0$ and hence \eqref{ter3} equals zero when $N=1$. For $N \geq 2$, we have 
\begin{align*}
    (\Sigma_2'-\mu_{12})\mu_{11}-\mu_{21}[\Sigma'_2-\mu_{22}]
    =& (\mu_{22}+\lambda \Sigma_2)\mu_{11}-\mu_{21}(\mu_{12}+\lambda\Sigma_2)\\
    =&(\mu_{11}\mu_{22}-\mu_{21}\mu_{12})+\lambda\mu_{11}\Sigma_2-\lambda\mu_{21}\Sigma_2 \\
    \stackrel{(*)}{\geq} &(\mu_{11}\mu_{22}-\mu_{21}\mu_{12})+\lambda\mu_{11}\Sigma_2-\frac{\mu_{21}}{\Sigma_1}(\mu_{11}\mu_{22}-\mu_{21}\mu_{12}) \\
    =&(\mu_{11}\mu_{22}-\mu_{21}\mu_{12}) \bigg[1-\frac{\mu_{21}}{\Sigma_1} \bigg]+\lambda \mu_{11}\Sigma_2\\
    =&(\mu_{11}\mu_{22}-\mu_{21}\mu_{12})\frac{\mu_{11}}{\Sigma_1}+\lambda\mu_{11}\Sigma_2\\
    \stackrel{(**)}{\geq} &0,
\end{align*}
where in $(**)$ we have used $\mu_{11}\mu_{22}-\mu_{21}\mu_{12} \geq 0$ and in $(*)$ we have used \eqref{lambda} since we have $\tau'(2)\geq 0$ when $N \geq 2$. It follows that \eqref{ter3} is non-negative for all possible values of $N$.

Finally, we show that \eqref{ter4}+\eqref{ter5} is non-negative. We have 
\begin{align*}
    &\eqref{ter4}+\eqref{ter5}\\
    =&\Sigma_1'(\Sigma'_2-\mu_{12})\bigg[(N-2)\theta \alpha(N)(\Sigma'_2+N\theta)+\mu_{11}^{N-1}(\Sigma'_2+N\theta)+\theta \alpha(N)(\Sigma'_2+N\theta) \\
   &+\mu_{22}(\Sigma'_2+N\theta)\sum_{k=3}^N \mu_{11}^{k-2}f(k,N)-\mu_{11}[\Sigma'_2+(N-1)\theta]\alpha(N)\bigg] \\
   =&\Sigma_1'(\Sigma'_2-\mu_{12})\bigg[(N-2)\theta \alpha(N)(\Sigma'_2+N\theta)+\mu_{11}^{N-1}(\Sigma'_2+N\theta)+\theta \alpha(N)(\Sigma'_2+N\theta) \\
   &+\mu_{22}(\Sigma'_2+N\theta)\sum_{k=3}^N \mu_{11}^{k-2}f(k,N)-\mu_{11}[\Sigma'_2+N\theta]\alpha(N)+\mu_{11}\theta\alpha(N)\bigg] \\
   =&\Sigma_1'(\Sigma'_2-\mu_{12})\bigg[(N-2)\theta \alpha(N)(\Sigma'_2+N\theta)+(\Sigma_2'+N\theta)\bigg( \mu_{11}^{N-1}+\theta \alpha(N)+\mu_{22}\sum_{k=3}^N \mu_{11}^{k-2}f(k,N)\\
   &-\mu_{11}\alpha(N)\bigg)+\mu_{11}\theta\alpha(N)\bigg]\\
   =&\Sigma_1'(\Sigma'_2-\mu_{12})\bigg[(N-2)\theta \alpha(N)(\Sigma'_2+N\theta)+(\Sigma_2'+N\theta)\bigg(\mu_{11}^{N-1}-\sum_{k=2}^N\mu_{11}^{k-1}f(k,N)+\theta\alpha(N)\\
   &+\mu_{22}\sum_{k=3}^N\mu_{11}^{k-2}f(k,N) \bigg)+\mu_{11}\theta\alpha(N)\bigg]\\
   =&\Sigma_1'(\Sigma'_2-\mu_{12})\bigg[(N-2)\theta \alpha(N)(\Sigma'_2+N\theta)+(\Sigma_2'+N\theta)\bigg(-\sum_{k=2}^{N-1}\mu_{11}^{k-1}f(k,N)+\theta\alpha(N)\\
   &+\mu_{22}\sum_{k=3}^N\mu_{11}^{k-2}f(k,N)  \bigg)+\mu_{11}\theta\alpha(N)\bigg].
\end{align*}
Let $k'=k-1$. Then 
\begin{align*}
 &\eqref{ter4}+\eqref{ter5}\\
    =&\Sigma_1'(\Sigma'_2-\mu_{12})\bigg[(N-2)\theta \alpha(N)(\Sigma'_2+N\theta)\\
   &+(\Sigma_2'+N\theta)\bigg(-\sum_{k=2}^{N-1}\mu_{11}^{k-1}f(k,N)+\theta\alpha(N)+\mu_{22}\sum_{k'=2}^{N-1}\mu_{11}^{k'-1}f(k'+1,N) \bigg)+\mu_{11}\theta\alpha(N)\bigg]\\
   =&\Sigma_1'(\Sigma'_2-\mu_{12})\bigg[(N-2)\theta \alpha(N)(\Sigma'_2+N\theta)\\
   &+(\Sigma_2'+N\theta)\bigg(-\sum_{k=2}^{N-1}\mu_{11}^{k-1}[\mu_{22}+(k-1)\theta]f(k+1,N)+\theta\alpha(N)+\mu_{22}\sum_{k'=2}^{N-1}\mu_{11}^{k'-1}f(k'+1,N) \bigg)\\
   &+\mu_{11}\theta\alpha(N)\bigg]\\
   =&\Sigma_1'(\Sigma'_2-\mu_{12})\bigg[(N-2)\theta \alpha(N)(\Sigma'_2+N\theta)+(\Sigma_2'+N\theta)\bigg(-\sum_{k=2}^{N-1}\mu_{11}^{k-1}(k-1)\theta f(k+1,N)+\theta\alpha(N) \bigg)\\
   &+\mu_{11}\theta\alpha(N)\bigg]\\
   =&\Sigma_1'(\Sigma'_2-\mu_{12})\bigg[(\Sigma_2'+N\theta)\bigg((N-2)\theta\alpha(N)-\sum_{k=2}^{N-1}\mu_{11}^{k-1}(k-1)\theta f(k+1,N)+\theta\alpha(N) \bigg)+\mu_{11}\theta\alpha(N)\bigg]\\
   \geq &\Sigma_1'(\Sigma'_2-\mu_{12})\bigg[ (\Sigma_2'+N\theta)\bigg((N-1)\theta\alpha(N)-(N-1)\theta\sum_{k=2}^{N-1}\mu_{11}^{k-1}f(k+1,N)\bigg)+\mu_{11}\theta\alpha(N)\bigg]\\
   =&\Sigma_1'(\Sigma'_2-\mu_{12})\bigg[(\Sigma_2'+N\theta)(N-1)\theta\bigg( \sum_{k=2}^N \mu_{11}^{k-2}f(k,N)-\sum_{k=2}^{N-1}\mu_{11}^{k-1}f(k+1,N)\bigg)+\mu_{11}\theta\alpha(N)\bigg].
\end{align*}
Let $k''=k+1$, We have 
\begin{align*}
 \eqref{ter4}+\eqref{ter5}
    =&\Sigma_1'(\Sigma'_2-\mu_{12})\bigg[(\Sigma_2'+N\theta)(N-1)\theta\bigg( \sum_{k=2}^N \mu_{11}^{k-2}f(k,N)-\sum_{k''=3}^{N}\mu_{11}^{k''-2}f(k'',N)\bigg)+\mu_{11}\theta\alpha(N)\bigg]\\
    =&\Sigma_1'(\Sigma'_2-\mu_{12})\bigg[(\Sigma_2'+N\theta)(N-1)\theta f(2,N)+\mu_{11}\theta\alpha(N)\bigg]\\
    \geq & 0.
\end{align*}
This proves $\Gamma(N,a_{21}) \geq 0$.

Now for $s=N, \ldots, B$,
\begin{align*}
    &\Gamma(s+1,a_{21})-\Gamma(s,a_{21})\\
    =&(g_N-\Sigma'_2)\bigg[\mu_{21}\bigg(\frac{1}{\Sigma'_2+(s+1)\theta}-\frac{1}{\Sigma'_2+s\theta}\bigg)+(\Sigma'_2-\mu_{12})\bigg( \frac{1}{\Sigma'_2+s\theta}-\frac{1}{\Sigma_2'+(s-1)\theta}\bigg) \bigg] >0
\end{align*}
as $(g_N-\Sigma'_2) \leq 0$ from \eqref{g_N-sigma_2} . This proves $\Gamma(s,a_{21}) \geq 0$ for $s=N, \ldots, B+1$.

\item[$(b)$] For $s=B+2$, we have
\begin{align*}
     \Gamma(B+2,a_{21})=&\Sigma'_2+\frac{\Sigma'_2+(B+1)\theta}{q}h_N(B+1)+\bigg(1-\frac{\Sigma'_2+(B+1)\theta}{q} \bigg)h_N(B+2)\\
    &-\mu_{12}-\frac{\mu_{12}+(B+1)\theta}{q}h_N(B+1)-\bigg(1-\frac{\mu_{12}+(B+1)\theta}{q}\bigg)h_N(B+2)\\
     =&\Sigma'_2-\mu_{12}+\frac{\Sigma'_2-\mu_{12}}{q}[h_N(B+1)-h_N(B+2)]\\
     =&\Sigma'_2-\mu_{12}+\frac{\Sigma'_2-\mu_{12}}{\Sigma'_2+(B+1)\theta}(g_N-\Sigma'_2)\\
     =&\frac{(\Sigma'_2-\mu_{12})[(B+1)\theta+g_N]}{\Sigma'_2+(B+1)\theta}>0,
\end{align*}
where the third equality follows from Lemma \ref{lemma h_N1}. Thus $\Gamma(s,a_{21}) \geq 0$ for all $s=N, \ldots, B+2$.
\item[Case 3, Part (iii) :]  $a=a_{11}$. 

We consider the two sub-cases  $(a)\  s=B+2$ and $(b)\ s=N. \ldots,B+1 $ separately.  
\item[$(a)$]$s=B+2$. We have from Lemma \ref{lemma h_N1}
\begin{align}\label{a11B+2(1)}
 \Gamma(B+2,a_{11})=&\Sigma'_2+\frac{\Sigma'_2+(B+1)\theta}{q}h_N(B+1)+\bigg( 1-\frac{\Sigma'_2+(B+1)\theta}{q}\bigg)h_N(B+2) \nonumber \\
    &-\frac{(B+2)\theta}{q}h_N(B+1)-\bigg(1-\frac{(B+2)\theta}{q} \bigg)h_N(B+2) \nonumber \\
    =&\Sigma_2'+\frac{\Sigma'_2}{q}[h_N(B+1)-h_N(B+2)]+\frac{\theta}{q}[h_N(B+2)-h_N(B+1)] \nonumber \\
    =&\Sigma_2'+\frac{\Sigma'_2}{\Sigma'_2+(B+1)\theta}(g_N-\Sigma_2')+\frac{\theta}{q}[h_N(B+2)-h_N(B+1)] \nonumber\\
    =&\frac{\Sigma'_2(B+1)\theta+\Sigma'_2g_N}{\Sigma'_2+(B+1)\theta}+\frac{\theta}{q}[h_N(B+2)-h_N(B+1)] \geq 0,
\end{align}
where the third equality follows from Lemma \ref{lemma h_N1} and the last inequality follows from Lemma \ref{lemma h_N2}.
\item[$(b)$]$s=N, \ldots, B+1$. We have 
\begin{align}\label{eq Gamma(s,a_11)} \Gamma(s,a_{11})=&\Sigma'_2+\frac{\Sigma'_2+(s-1)\theta}{q}h_N(s-1)+\bigg( 1-\frac{\Sigma'_2+(s-1)\theta}{q}\bigg)h_N(s) \nonumber\\
    &-\frac{\Sigma'_1}{q}h_N(s+1)-\frac{s\theta}{q}h_N(s-1)-\bigg(1-\frac{\Sigma'_1+s\theta}{q} \bigg)h_N(s)  \nonumber\\
    =&\Sigma'_2+\frac{\Sigma'_2-\theta}{q}h_N(s-1)-\frac{\Sigma'_1}{q}h_N(s+1)-\frac{\Sigma'_2-\theta}{q}h_N(s)+\frac{\Sigma'_1}{q}h_N(s) \nonumber\\
    =&\Sigma'_2+\frac{\Sigma'_2}{q}[h_N(s-1)-h_N(s)]+\frac{\Sigma'_1}{q}[h_N(s)-h_N(s+1)]+\frac{\theta}{q}[h_N(s)-h_N(s-1)].
\end{align}
In order to show that $\Gamma(s,a_{11})$ is non-negative, we consider two scenarios; $(1)$$N=1$ and $(2)N \geq 2$ separately.
\\\\
\underline{\textbf{Scenario 1: $N=1$}. }

We have for $s=N, \ldots, B+1$, using Lemma \ref{lemma h_N1}
\begin{align*}
    \Gamma(s,a_{11})=&\Sigma'_2+\frac{\Sigma'_2}{\Sigma'_2+(s-1)\theta}(g_N-\Sigma'_2)+\frac{\Sigma'_1}{\Sigma'_2+s\theta}(g_N-\Sigma'_2)+\frac{\theta}{q}[h_N(s)-h_N(s-1)]\\
    =&\Gamma_1(s,a_{11})+\frac{\theta}{q}[h_N(s)-h_N(s-1)],
\end{align*}
where 
\begin{align*}
\Gamma_1(s,a_{11})=\Sigma'_2+\frac{\Sigma'_2}{\Sigma'_2+(s-1)\theta}(g_N-\Sigma'_2)+\frac{\Sigma'_1}{\Sigma'_2+s\theta}(g_N-\Sigma'_2).
\end{align*} 
Since $h_N(s)-h_N(s-1) \geq 0$ from Lemma \ref{lemma h_N2}, it is sufficient to show $\Gamma_1(s,a_{11})\geq 0$ for $s=N, \ldots,B+1$.
Now
\begin{align*}
    g_N-\Sigma'_2=g_1-\Sigma'_2=&\frac{\beta(1)}{\Delta(1)}-\Sigma'_2=\frac{\Sigma'_1\Sigma'_2}{\Sigma'_2+\Sigma'_1}-\Sigma'_2=\frac{-(\Sigma'_2)^2}{\Sigma'_1+\Sigma'_2}.
\end{align*}
Hence
\begin{align*}
   \Gamma_1(N,a_{11})=\Gamma_1(1,a_{11})
   =&\Sigma'_2-\frac{(\Sigma'_2)^2}{\Sigma'_1+\Sigma'_2}-\frac{\Sigma'_1(\Sigma'_2)^2}{(\Sigma'_2+\theta)(\Sigma'_1+\Sigma'_2)}\\
   =&\frac{\Sigma'_1\Sigma'_2}{\Sigma'_1+\Sigma'_2}-\frac{\Sigma'_1(\Sigma'_2)^2}{(\Sigma'_2+\theta)(\Sigma'_1+\Sigma'_2)}\\
   =&\frac{\Sigma'_1\Sigma'_2\theta}{(\Sigma'_2+\theta)(\Sigma'_1+\Sigma'_2)} \geq 0.
\end{align*}
Furthermore for $s=N, \ldots, B$, we have
\begin{align*}
   & \Gamma_1(s+1,a_{11})-\Gamma_1(s,a_{11})\\
   =&(g_N-\Sigma'_2)\bigg[\Sigma'_1 \bigg(\frac{1}{\Sigma'_2+(s+1)\theta}-\frac{1}{\Sigma'_2+s\theta} \bigg)+\Sigma_2' \bigg( \frac{1}{\Sigma'_2+s\theta}-\frac{1}{\Sigma'_2+(s-1)\theta}\bigg) \bigg] \geq 0,
\end{align*}
where the last inequality follows as $g_N-\Sigma_2 \leq 0$ from \eqref{g_N-sigma_2}.
Therefore $\Gamma_1(s,a_{11}) \geq 0$ for $s=N, \ldots, B+1$ 
and hence $\Gamma(s,a_{11}) \geq 0$ for $s=N, \ldots, B+1$ (if $N=1$).

\underline{\textbf{Scenario 2: $N \geq 2$}.}

We have from \eqref{eq Gamma(s,a_11)} for $s=N, \ldots, B+1$
\begin{align*}
    \Gamma(s,a_{11})=&\Sigma'_2+\frac{\Sigma'_2}{q}[h_N(s-1)-h_N(s)]+\frac{\Sigma'_1}{q}[h_N(s)-h_N(s+1)]+\frac{\theta}{q}[h_N(s)-h_N(s-1)]\\
    =&\Sigma'_2-\mu_{22}+\frac{\Sigma'_2-\mu_{22}}{q}[h_N(s-1)-h_N(s)]+\frac{\mu_{11}}{q}[h_N(s)-h_N(s+1)]\\
          &+\mu_{22}+\frac{\mu_{22}}{q}[h_N(s-1)-h_N(s)]+\frac{\Sigma_1'-\mu_{11}}{q}[h_N(s)-h_N(s+1)]\\
          &+\frac{\theta}{q}[h_N(s)-h_N(s-1)]\\
        =&\Gamma(s,a_{12})+\Gamma_2(s,a_{11})+\frac{\theta}{q}[h_N(s)-h_N(s-1)],
\end{align*}
where $\Gamma_2(s,a_{11})=\mu_{22}+\frac{\mu_{22}}{q}[h_N(s-1)-h_N(s)]+\frac{\Sigma_1'-\mu_{11}}{q}[h_N(s)-h_N(s+1)]$ and $\Gamma(s,a_{12})$ as defined in \eqref{gamma(s,a_12).2}. Since we have shown earlier $\Gamma(s,a_{12}) \geq 0$ and we have $h_N(s)-h_N(s-1) \geq 0$ from Lemma \ref{lemma h_N2}, it is sufficient to show $\Gamma_2(s,a_{11}) \geq 0$ for $s=N, \ldots, B+1$. Using Lemma \ref{lemma h_N1},
\begin{align*}
    \Gamma_2(s,a_{11})
    =&\mu_{22}+\frac{\mu_{22}}{\Sigma'_2+(s-1)\theta}(g_N-\Sigma'_2)+\frac{\Sigma'_1-\mu_{11}}{\Sigma'_2+s\theta}(g_N-\Sigma'_2).
\end{align*}
From \eqref{g_N-sigma_2} we obtain
\begin{align*}
    \Gamma_2(N,a_{11})=&\mu_{22}-\frac{\mu_{22}}{\Delta(N)}\bigg(\Sigma'_2f(1,N)+\Sigma'_1[\Sigma'_2-\mu_{22}]\alpha(N)\bigg)\\
    &-\frac{\Sigma'_1-\mu_{11}}{\Sigma'_2+N\theta} \times\frac{[\Sigma'_2+(N-1)\theta]\bigg(\Sigma'_2f(1,N)+\Sigma'_1[\Sigma'_2-\mu_{22}]\alpha(N)\bigg)}{\Delta(N)},
\end{align*}
and hence
\begin{align*}
    &\Gamma_2(N,a_{11})\Delta(N)(\Sigma'_2+N \theta) \nonumber\\
    &=(\Sigma'_2+N\theta)\Delta(N)\mu_{22}-(\Sigma_2'+N \theta)\bigg(\Sigma'_2f(1,N)+\Sigma'_1[\Sigma'_2-\mu_{22}]\alpha(N)\bigg)\mu_{22} \nonumber \\
    &-(\Sigma'_1-\mu_{11})[\Sigma'_2+(N-1)\theta]\bigg(\Sigma'_2f(1,N)+\Sigma'_1[\Sigma'_2-\mu_{22}]\alpha(N)\bigg).
\end{align*}

We consider the three terms in the above separately. The first term satisfies
\begin{align}\label{Gamma_2(N,a11)Delta(N)(Sigma'_2+N theta) eq1}
    &(\Sigma'_2+N\theta)\Delta(N)\mu_{22} \nonumber \\
    =&\mu_{22} \bigg[[\Sigma'_2+(N-1)\theta][f(1,N)+\Sigma'_1\alpha(N)]+\Sigma'_1\mu_{11}^{N-1} \bigg](\Sigma'_2+N\theta) \nonumber \\
    =&\mu_{22}\Sigma'_2f(1,N)(\Sigma'_2+N\theta)+\Sigma'_1\mu_{22}\Sigma'_2\alpha(N)(\Sigma'_2+N\theta)+\Sigma'_1\mu_{11}^{N-1}\mu_{22}(\Sigma'_2+N\theta) \nonumber \\
    &+(\Sigma'_2+N\theta)\mu_{22}(N-1)\theta f(1,N)+\Sigma'_1\mu_{22}(N-1)\theta\alpha(N)(\Sigma'_2+N\theta) \nonumber \\
    =&\mu_{22}\Sigma'_2f(1,N)(\Sigma'_2+N\theta)+\Sigma'_1\mu_{22}\Sigma'_2\alpha(N)(\Sigma'_2+N\theta)+\Sigma'_1\mu_{11}^{N-1}\mu_{22}(\Sigma'_2+N\theta) \nonumber \\
    &+(\Sigma'_2+N\theta)\mu_{22}(N-1)\theta f(1,N)+\Sigma'_1\mu_{22}(N-2)\theta\alpha(N)(\Sigma'_2+N\theta)+\Sigma'_1\mu_{22}\theta\alpha(N)(\Sigma'_2+N \theta).
\end{align}
The second term satisfies
\begin{align}\label{Gamma_2(N,a11)Delta(N)(Sigma'_2+N theta) eq2}
    &-(\Sigma_2'+N \theta)\bigg(\Sigma'_2f(1,N)+\Sigma'_1[\Sigma'_2-\mu_{22}]\alpha(N)\bigg)\mu_{22} \nonumber  \\
    =&(\Sigma_2'+N \theta)\bigg(-\Sigma'_2f(1,N)+\Sigma'_1[\mu_{22}-\Sigma'_2]\alpha(N)\bigg)\mu_{22}\nonumber  \\
    =&-\mu_{22}\Sigma'_2 f(1,N)(\Sigma_2'+N \theta)+\Sigma'_1 \mu_{22}^2\alpha(N)(\Sigma_2'+N \theta)-\Sigma'_1 \Sigma'_2 \mu_{22} \alpha(N)(\Sigma_2'+N \theta) \nonumber \\
    =&-\mu_{22}\Sigma'_2 f(1,N)(\Sigma_2'+N \theta)-\Sigma'_1 \Sigma'_2 \mu_{22} \alpha(N)(\Sigma_2'+N \theta)+\Sigma'_1 \mu_{22}^2(\Sigma_2'+N \theta)\sum_{k=2}^N\mu_{11}^{k-2}f(k,N) \nonumber \\
    =&-\mu_{22}\Sigma'_2 f(1,N)(\Sigma_2'+N \theta)-\Sigma'_1 \Sigma'_2 \mu_{22} \alpha(N)(\Sigma_2'+N \theta)+\Sigma'_1 \mu_{22}^2(\Sigma_2'+N \theta)\sum_{k=3}^N\mu_{11}^{k-2}f(k,N) \nonumber \\
    &+\Sigma'_1 \mu_{22}(\Sigma_2'+N \theta)f(1,N) \nonumber \\
     =&-\mu_{22}\Sigma'_2 f(1,N)(\Sigma_2'+N \theta)-\Sigma'_1 \Sigma'_2 \mu_{22} \alpha(N)(\Sigma_2'+N \theta)+\Sigma'_1 \mu_{22}^2(\Sigma_2'+N \theta)\sum_{k=3}^N\mu_{11}^{k-2}f(k,N) \nonumber \\
    &+\Sigma'_1 \mu_{22}(\Sigma_2'+(N-1) \theta)f(1,N)+\Sigma'_1 \mu_{22}\theta f(1,N).
\end{align}
And the third term satisfies
\begin{align}\label{Gamma_2(N,a11)Delta(N)(Sigma'_2+N theta) eq3}
    &-(\Sigma'_1-\mu_{11})[\Sigma'_2+(N-1)\theta]\bigg(\Sigma'_2f(1,N)+\Sigma'_1[\Sigma'_2-\mu_{22}]\alpha(N)\bigg) \nonumber \\
    =&(\Sigma'_1-\mu_{11})[\Sigma'_2+(N-1)\theta]\bigg(-\Sigma'_2f(1,N)+\Sigma'_1[\mu_{22}-\Sigma'_2]\alpha(N)\bigg) \nonumber \\
    =&(\mu_{11}\Sigma'_2-\Sigma'_1\Sigma'_2)[\Sigma'_2+(N-1)\theta]f(1,N)+[\Sigma'_2+(N-1)\theta](\Sigma'_1-\mu_{11})(\Sigma'_1\mu_{22}-\Sigma'_1\Sigma'_2)\alpha(N) \nonumber \\
    =&(\mu_{11}\Sigma'_2-\Sigma'_1\Sigma'_2)[\Sigma'_2+(N-1)\theta]f(1,N)+[\Sigma'_2+(N-1)\theta](\Sigma'_1\mu_{22}-\Sigma'_1\Sigma'_2)\Sigma'_1\alpha(N) \nonumber \\
    &-\mu_{11}[\Sigma'_2+(N-1)\theta](\Sigma'_1\mu_{22}-\Sigma'_1\Sigma'_2)\alpha(N) \nonumber \\
    =&(\mu_{11}\Sigma'_2-\Sigma'_1\Sigma'_2)[\Sigma'_2+(N-1)\theta]f(1,N)+[\Sigma'_2+(N-1)\theta](\Sigma'_1\mu_{22}+\mu_{11}\Sigma'_2-\Sigma'_1\Sigma'_2)\Sigma'_1\alpha(N) \nonumber \\
    &-\mu_{11}[\Sigma'_2+(N-1)\theta](\Sigma'_1\mu_{22}-\Sigma'_1\Sigma'_2)\alpha(N)-\mu_{11}\Sigma'_1\Sigma'_2[\Sigma'_2+(N-1)\theta]\alpha(N) \nonumber \\
    =&(\mu_{11}\Sigma'_2-\Sigma'_1\Sigma'_2)[\Sigma'_2+(N-1)\theta]f(1,N)+[\Sigma'_2+(N-1)\theta](\Sigma'_1\mu_{22}+\mu_{11}\Sigma'_2-\Sigma'_1\Sigma'_2)\Sigma'_1\alpha(N) \nonumber \\
    &-\mu_{11}\Sigma'_1 \mu_{22}[\Sigma'_2+(N-1)\theta]\alpha(N).
\end{align}
Adding \eqref{Gamma_2(N,a11)Delta(N)(Sigma'_2+N theta) eq1}-\eqref{Gamma_2(N,a11)Delta(N)(Sigma'_2+N theta) eq3} yields
\begin{align*}
   & \Gamma_2(N,a_{11})\Delta(N)(\Sigma'_2+N \theta) \\
   =&\mu_{22}\Sigma'_2f(1,N)(\Sigma'_2+N\theta)+\Sigma'_1\mu_{22}\Sigma'_2\alpha(N)(\Sigma'_2+N\theta)+\Sigma'_1\mu_{11}^{N-1}\mu_{22}(\Sigma'_2+N\theta)\\
    &+(\Sigma'_2+N\theta)\mu_{22}(N-1)\theta f(1,N)+\Sigma'_1\mu_{22}(N-2)\theta\alpha(N)(\Sigma'_2+N\theta)+\Sigma'_1\mu_{22}\theta\alpha(N)(\Sigma'_2+N \theta)\\
    &-\mu_{22}\Sigma'_2 f(1,N)(\Sigma_2'+N \theta)-\Sigma'_1 \Sigma'_2 \mu_{22} \alpha(N)(\Sigma_2'+N \theta)+\Sigma'_1 \mu_{22}^2(\Sigma_2'+N \theta)\sum_{k=3}^N\mu_{11}^{k-2}f(k,N)\\
    &+\Sigma'_1 \mu_{22}(\Sigma_2'+(N-1) \theta)f(1,N)+\Sigma'_1 \mu_{22}\theta f(1,N)\\
    &+(\mu_{11}\Sigma'_2-\Sigma'_1\Sigma'_2)[\Sigma'_2+(N-1)\theta]f(1,N)+[\Sigma'_2+(N-1)\theta](\Sigma'_1\mu_{22}+\mu_{11}\Sigma'_2-\Sigma'_1\Sigma'_2)\Sigma'_1\alpha(N)\\
    &-\mu_{11}\Sigma'_1 \mu_{22}[\Sigma'_2+(N-1)\theta]\alpha(N)\\
    =&\Sigma'_1\mu_{11}^{N-1}\mu_{22}(\Sigma'_2+N\theta)\\
    &+(\Sigma'_2+N\theta)\mu_{22}(N-1)\theta f(1,N)+\Sigma'_1\mu_{22}(N-2)\theta\alpha(N)(\Sigma'_2+N\theta)+\Sigma'_1\mu_{22}\theta\alpha(N)(\Sigma'_2+N \theta)\\
    &+\Sigma'_1 \mu_{22}^2(\Sigma_2'+N \theta)\sum_{k=3}^N\mu_{11}^{k-2}f(k,N)\\
    &+\Sigma'_1 \mu_{22}(\Sigma_2'+(N-1) \theta)f(1,N)+\Sigma'_1 \mu_{22}\theta f(1,N)\\
    &+(\mu_{11}\Sigma'_2-\Sigma'_1\Sigma'_2)[\Sigma'_2+(N-1)\theta]f(1,N)+[\Sigma'_2+(N-1)\theta](\Sigma'_1\mu_{22}+\mu_{11}\Sigma'_2-\Sigma'_1\Sigma'_2)\Sigma'_1\alpha(N)\\
    &-\mu_{11}\Sigma'_1 \mu_{22}[\Sigma'_2+(N-1)\theta]\alpha(N)\\
   = &(\Sigma'_2+N\theta)\mu_{22}(N-1)\theta f(1,N)+\Sigma'_1 \mu_{22}\theta f(1,N) \refstepcounter{equation}\tag{\theequation}\label{term1}\\
   &+(\Sigma'_1\mu_{22}+\mu_{11}\Sigma'_2-\Sigma'_1\Sigma'_2)(\Sigma_2'+(N-1) \theta)\bigg(f(1,N)+\Sigma'_1 \alpha(N)\bigg)\refstepcounter{equation}\tag{\theequation} \label{term2}\\
    &+\Sigma'_1\mu_{22}(N-2)\theta\alpha(N)(\Sigma'_2+N\theta)+\Sigma'_1\mu_{11}^{N-1}\mu_{22}(\Sigma'_2+N\theta)+\Sigma'_1\mu_{22}\theta\alpha(N)(\Sigma'_2+N \theta)\refstepcounter{equation}\tag{\theequation} \label{term3}\\
    &+\Sigma'_1 \mu_{22}^2(\Sigma_2'+N \theta)\sum_{k=3}^N\mu_{11}^{k-2}f(k,N)-\mu_{11}\Sigma'_1 \mu_{22}[\Sigma'_2+(N-1)\theta]\alpha(N).\refstepcounter{equation}\tag{\theequation} \label{term4}
\end{align*}

We show all the terms \eqref{term1}-\eqref{term4} is non-negative. \eqref{term1} is clearly non-negative.

We have \eqref{term2} is non negative as $\Sigma'_1\mu_{22}+\mu_{11}\Sigma'_2-\Sigma'_1\Sigma'_2 \geq 0$ from \eqref{N geq 2 implication} since $N \geq 2$ and hence $\tau'(2) \geq 0$.

 Finally we show \eqref{term3}+\eqref{term4} is non-negative.
 \begin{align*}
    & \eqref{term3}+\eqref{term4}\\
    =&\Sigma'_1\mu_{22}(N-2)\theta\alpha(N)(\Sigma'_2+N\theta)+\Sigma'_1\mu_{11}^{N-1}\mu_{22}(\Sigma'_2+N\theta)+\Sigma'_1\mu_{22}\theta\alpha(N)(\Sigma'_2+N \theta) \\
    &+\Sigma'_1 \mu_{22}^2(\Sigma_2'+N \theta)\sum_{k=3}^N\mu_{11}^{k-2}f(k,N)-\mu_{11}\Sigma'_1 \mu_{22}[\Sigma'_2+N\theta]\alpha(N)+\mu_{11}\Sigma'_1 \mu_{22}\theta \alpha(N) \\
    =&\Sigma'_1\mu_{22}(N-2)\theta\alpha(N)(\Sigma'_2+N\theta)\\
     &+\Sigma'_1\mu_{22}(\Sigma'_2+N\theta)\bigg[\mu_{11}^{N-1} +\theta \alpha(N)+\mu_{22}\sum_{k=3}^N\mu_{11}^{k-2}f(k,N)-\mu_{11}\alpha(N)\bigg]+\mu_{11}\Sigma'_1 \mu_{22}\theta \alpha(N) \\
    =&\Sigma'_1\mu_{22}(N-2)\theta\alpha(N)(\Sigma'_2+N\theta)\\
     &+\Sigma'_1\mu_{22}(\Sigma'_2+N\theta)\bigg[\mu_{11}^{N-1}-\sum_{k=2}^{N}\mu_{11}^{k-1}f(k,N) +\theta \alpha(N)+\mu_{22}\sum_{k=3}^N\mu_{11}^{k-2}f(k,N)\bigg]+\mu_{11}\Sigma'_1 \mu_{22}\theta \alpha(N)   \\
      =&\Sigma'_1\mu_{22}(N-2)\theta\alpha(N)(\Sigma'_2+N\theta)\\
     &+\Sigma'_1\mu_{22}(\Sigma'_2+N\theta)\bigg[-\sum_{k=2}^{N-1}\mu_{11}^{k-1}f(k,N)+\theta\alpha(N)+\mu_{22}\sum_{k=3}^N\mu_{11}^{k-2}f(k,N) \bigg]+\mu_{11}\Sigma'_1 \mu_{22}\theta \alpha(N)   .
 \end{align*}
 Let $k'=k-1$. Then
 \begin{align*}
  & \eqref{term3}+\eqref{term4}\\
      =&\Sigma'_1\mu_{22}(N-2)\theta\alpha(N)(\Sigma'_2+N\theta)\\
     &+\Sigma'_1\mu_{22}(\Sigma'_2+N\theta)\bigg[-\sum_{k=2}^{N-1}\mu_{11}^{k-1}f(k,N)+\theta\alpha(N)+\mu_{22}\sum_{k'=2}^{N-1}\mu_{11}^{k'-1}f(k'+1,N) \bigg]\\
     &+\mu_{11}\Sigma'_1 \mu_{22}\theta \alpha(N)   \\
     =&\Sigma'_1\mu_{22}(N-2)\theta\alpha(N)(\Sigma'_2+N\theta)\\
     &+\Sigma'_1\mu_{22}(\Sigma'_2+N\theta)\bigg[-\sum_{k=2}^{N-1}\mu_{11}^{k-1}[\mu_{22}+(k-1)\theta]f(k+1,N)+\theta\alpha(N)\\
     &+\mu_{22}\sum_{k'=2}^{N-1}\mu_{11}^{k'-1}f(k'+1,N)\bigg]+\mu_{11}\Sigma'_1 \mu_{22}\theta \alpha(N)   \\
     =&\Sigma'_1\mu_{22}(N-2)\theta\alpha(N)(\Sigma'_2+N\theta)\\
     &+\Sigma'_1\mu_{22}(\Sigma'_2+N\theta)\bigg[-\sum_{k=2}^{N-1}\mu_{11}^{k-1}(k-1)\theta f(k+1,N)+\theta\alpha(N)\bigg]+\mu_{11}\Sigma'_1 \mu_{22}\theta \alpha(N)   \\
      =&\Sigma'_1\mu_{22}(\Sigma'_2+N\theta)\bigg[(N-2)\theta \alpha(N)-\sum_{k=2}^{N-1}\mu_{11}^{k-1}(k-1)\theta f(k+1,N)+\theta\alpha(N)\bigg]+\mu_{11}\Sigma'_1 \mu_{22}\theta \alpha(N)  \\
       \geq &\Sigma'_1\mu_{22}(\Sigma'_2+N\theta)\bigg[(N-1)\theta \alpha(N)-(N-1)\theta\sum_{k=2}^{N-1}\mu_{11}^{k-1} f(k+1,N)\bigg]+\mu_{11}\Sigma'_1 \mu_{22}\theta \alpha(N)\\
        =&\Sigma'_1\mu_{22}(\Sigma'_2+N\theta)(N-1)\theta\bigg[ \sum_{k=2}^N \mu_{11}^{k-2}f(k,N)-\sum_{k=2}^{N-1}\mu_{11}^{k-1} f(k+1,N)\bigg]+\mu_{11}\Sigma'_1 \mu_{22}\theta \alpha(N).
 \end{align*}
 Let $k''=k+1$. Then
 \begin{align*}
  & \eqref{term3}+\eqref{term4}\\
      =&\Sigma'_1\mu_{22}(\Sigma'_2+N\theta)(N-1)\theta\bigg[ \sum_{k=2}^N \mu_{11}^{k-2}f(k,N)-\sum_{k''=3}^{N}\mu_{11}^{k''-2} f(k'',N)\bigg]+\mu_{11}\Sigma'_1 \mu_{22}\theta \alpha(N)\\
      =&\Sigma'_1\mu_{22}(\Sigma'_2+N\theta)(N-1)\theta f(2,N)+\mu_{11}\Sigma'_1 \mu_{22}\theta \alpha(N) \geq 0.
 \end{align*}
 This proves $\Gamma_2(N,a_{11}) \geq 0$. Furthermore for $s=N, \ldots, B$
 \begin{align*}
    & \Gamma_2(s+1,a_{11})-\Gamma_2(s,a_{11})\\
    &=(g_N-\Sigma'_2)\bigg[(\Sigma'_1-\mu_{11})\bigg(\frac{1}{\Sigma'_2+(s+1)\theta}-\frac{1}{\Sigma'_2+s\theta}\bigg)+\mu_{22}\bigg(\frac{1}{\Sigma'_2+s\theta}-\frac{1}{\Sigma'_2+(s-1)\theta} \bigg) \bigg] \geq 0,
 \end{align*}
 where the inequality follows as $(g_N-\Sigma'_2) \leq 0$ from \eqref{g_N-sigma_2}. Therefore $\Gamma_2(s,a_{11}) \geq 0$ for $s=N, \ldots, B+1$.
And hence $\Gamma(s,a_{11}) \geq 0$ for $s=N, \ldots, B+1$ (if $N\geq 2$). 

Hence sub-case $(a)$ and the  two scenarios for sub-case $(b)$ prove that $\Gamma(s,a_{11}) \geq 0$ for $s=N, \ldots, B+2$. 
This completes the proof of the theorem.
\end{proof}
\subsection{Details of the proof of Proposition 11}\label{SMProp11}
We have
\begin{align*}
    &\tau_\gamma'(n)=[\Sigma'_2+(n-1)\theta][\Sigma'_2+(n-2)\theta]\mu_{22}f(1,n-1)-[\Sigma'_2+(n-1)\theta][\Sigma'_1\alpha(n)(\Sigma'_2-\mu_{22})+\Sigma'_2f(1,n)]\\
     & \quad \quad \quad +[\Sigma'_2+(n-2)\theta]\mu_{11}[\Sigma'_2f(1,n-1)+\Sigma'_1\alpha(n-1)(\Sigma'_2-\mu_{22})]\\
    =&[\Sigma'_2+(n-1)\theta][\Sigma'_2+(n-2)\theta]\mu_{22}f(1,n-1)-(1+\lambda)\bigg[[\Sigma'_2+(n-1)\theta][\Sigma_1\alpha(n)(\Sigma'_2-\mu_{22})+\Sigma_2f(1,n)]\\
     & \quad \quad \quad +[\Sigma'_2+(n-2)\theta]\mu_{11}[\Sigma_2f(1,n-1)+\Sigma_1\alpha(n-1)(\Sigma'_2-\mu_{22})]\bigg]\\
     =&(1+\lambda) \bigg[[\Sigma_2+(n-1)\theta][\Sigma_2+(n-2)\theta]\mu_{22}f(1,n-1)-[\Sigma'_2+(n-1)\theta][\Sigma_1\alpha(n)(\Sigma'_2-\mu_{22})+\Sigma_2f(1,n)]\\
     & \quad \quad \quad +[\Sigma'_2+(n-2)\theta]\mu_{11}[\Sigma_2f(1,n-1)+\Sigma_1\alpha(n-1)(\Sigma'_2-\mu_{22})] \bigg]\\
     &+[\Sigma'_2+(n-1)\theta][\Sigma'_2+(n-2)\theta]\mu_{22}f(1,n-1)-(1+\lambda)[\Sigma_2+(n-1)\theta][\Sigma_2+(n-2)\theta]\mu_{22}f(1,n-1)\\
     =&(1+\lambda) \bigg[[\Sigma_2+(n-1)\theta][\Sigma_2+(n-2)\theta]\mu_{22}f(1,n-1)-[\Sigma'_2+(n-1)\theta][\Sigma_1\alpha(n)(\Sigma_2-\mu_{22})+\Sigma_2f(1,n)]\\
     & \quad \quad \quad +[\Sigma'_2+(n-2)\theta]\mu_{11}[\Sigma_2f(1,n-1)+\Sigma_1\alpha(n-1)(\Sigma_2-\mu_{22})] \bigg]\\
     &-(1+\lambda) \bigg[[\Sigma_2'+(n-1)\theta]\Sigma_1 \alpha(n)\lambda\Sigma_2 \bigg]+(1+\lambda) \bigg[ [\Sigma'_2+(n-2)\theta]\mu_{11}\Sigma_1\alpha(n-1)\lambda \Sigma_2 \bigg]\\
     &+[\Sigma'_2+(n-1)\theta][\Sigma'_2+(n-2)\theta]\mu_{22}f(1,n-1)-(1+\lambda)[\Sigma_2+(n-1)\theta][\Sigma_2+(n-2)\theta]\mu_{22}f(1,n-1)\\
     =&(1+\lambda) \bigg[[\Sigma_2+(n-1)\theta][\Sigma_2+(n-2)\theta]\mu_{22}f(1,n-1)-[\Sigma_2+(n-1)\theta][\Sigma_1\alpha(n)(\Sigma_2-\mu_{22})+\Sigma_2f(1,n)]\\
     & \quad \quad \quad +[\Sigma_2+(n-2)\theta]\mu_{11}[\Sigma_2f(1,n-1)+\Sigma_1\alpha(n-1)(\Sigma_2-\mu_{22})] \bigg]\\
     &-(1+\lambda)\lambda\Sigma_2[\Sigma_1\alpha(n)(\Sigma_2-\mu_{22})+\Sigma_2f(1,n)]+(1+\lambda) \lambda \Sigma_2\mu_{11}[\Sigma_2f(1,n-1)+\Sigma_1\alpha(n-1)(\Sigma_2-\mu_{22})]\\
     &-(1+\lambda) \bigg[[\Sigma_2'+(n-1)\theta]\Sigma_1 \alpha(n)\lambda\Sigma_2 \bigg]+(1+\lambda) \bigg[ [\Sigma'_2+(n-2)\theta]\mu_{11}\Sigma_1\alpha(n-1)\lambda \Sigma_2 \bigg]\\
     &+[\Sigma'_2+(n-1)\theta][\Sigma'_2+(n-2)\theta]\mu_{22}f(1,n-1)-(1+\lambda)[\Sigma_2+(n-1)\theta][\Sigma_2+(n-2)\theta]\mu_{22}f(1,n-1)\\
     =&(1+\lambda)\tau'_1(n)\\
      &-(1+\lambda)\lambda\Sigma_2[\Sigma_1\alpha(n)(\Sigma_2-\mu_{22})+\Sigma_2f(1,n)]+(1+\lambda) \lambda \Sigma_2\mu_{11}[\Sigma_2f(1,n-1)+\Sigma_1\alpha(n-1)(\Sigma_2-\mu_{22})]\\
     &-(1+\lambda) \bigg[[\Sigma_2'+(n-1)\theta]\Sigma_1 \alpha(n)\lambda\Sigma_2 \bigg]+(1+\lambda) \bigg[ [\Sigma'_2+(n-2)\theta]\mu_{11}\Sigma_1\alpha(n-1)\lambda \Sigma_2 \bigg]\\
     &+[\Sigma'_2+(n-1)\theta][\Sigma'_2+(n-2)\theta]\mu_{22}f(1,n-1)-(1+\lambda)[\Sigma_2+(n-1)\theta][\Sigma_2+(n-2)\theta]\mu_{22}f(1,n-1)\\
   =& (1+\lambda)\tau'_1(n)\\
   &-(1+\lambda)\lambda\Sigma_2[\Sigma_1\alpha(n)(\Sigma_2-\mu_{22})+\Sigma_2f(1,n)]+(1+\lambda) \lambda \Sigma_2\mu_{11}[\Sigma_2f(1,n-1)+\Sigma_1\alpha(n-1)(\Sigma_2-\mu_{22})]\\
     &-(1+\lambda) \bigg[[\Sigma_2'+(n-1)\theta]\Sigma_1 \alpha(n)\lambda\Sigma_2 \bigg]+(1+\lambda) \bigg[ [\Sigma'_2+(n-2)\theta]\mu_{11}\Sigma_1\alpha(n-1)\lambda \Sigma_2 \bigg]\\
     &+\mu_{22} f(1,n-1) \bigg[[\Sigma'_2+(n-1)\theta][\Sigma'_2+(n-2)\theta]-(1+\lambda)[\Sigma_2+(n-1)\theta][\Sigma_2+(n-2)\theta] \bigg]\\
      =& (1+\lambda)\tau'_1(n)\\
   &-(1+\lambda)\lambda\Sigma_2[\Sigma_1\alpha(n)(\Sigma_2-\mu_{22})+\Sigma_2f(1,n)]+(1+\lambda) \lambda \Sigma_2\mu_{11}[\Sigma_2f(1,n-1)+\Sigma_1\alpha(n-1)(\Sigma_2-\mu_{22})]\\
     &-(1+\lambda) \bigg[[\Sigma_2'+(n-1)\theta]\Sigma_1 \alpha(n)\lambda\Sigma_2 \bigg]+(1+\lambda) \bigg[ [\Sigma'_2+(n-2)\theta]\mu_{11}\Sigma_1\alpha(n-1)\lambda \Sigma_2 \bigg]\\
     &+\mu_{22} f(1,n-1) \bigg[\Sigma'_2
     \Sigma'_2+\Sigma'_2 (n-2)\theta+\Sigma'_2(n-1)\theta+(n-1)(n-2)\theta^2-(1+\lambda)\Sigma_2\Sigma_2\\
     &-(1+\lambda)\Sigma_2(n-1)\theta-(1+\lambda)\Sigma_2(n-2)\theta-(1+\lambda)(n-1)(n-2)\theta^2\bigg]\\
      =& (1+\lambda)\tau'_1(n)\\
   &-(1+\lambda)\lambda\Sigma_2[\Sigma_1\alpha(n)(\Sigma_2-\mu_{22})+\Sigma_2f(1,n)]+(1+\lambda) \lambda \Sigma_2\mu_{11}[\Sigma_2f(1,n-1)+\Sigma_1\alpha(n-1)(\Sigma_2-\mu_{22})]\\
     &-(1+\lambda) \bigg[[\Sigma_2'+(n-1)\theta]\Sigma_1 \alpha(n)\lambda\Sigma_2 \bigg]+(1+\lambda) \bigg[ [\Sigma'_2+(n-2)\theta]\mu_{11}\Sigma_1\alpha(n-1)\lambda \Sigma_2 \bigg]\\
     &+\mu_{22} f(1,n-1) \bigg[\Sigma'_2
     \Sigma'_2-(1+\lambda)\Sigma_2\Sigma_2-\lambda(n-1)(n-2)\theta^2\bigg]\\
     =&(1+\lambda)\tau'_1(n)\\
      &-(1+\lambda)\lambda\Sigma_2[\Sigma_1\alpha(n)(\Sigma_2-\mu_{22})+\Sigma_2f(1,n)]+(1+\lambda) \lambda \Sigma_2\mu_{11}[\Sigma_2f(1,n-1)+\Sigma_1\alpha(n-1)(\Sigma_2-\mu_{22})]\\
     &-(1+\lambda) \bigg[[\Sigma_2'+(n-1)\theta]\Sigma_1 \alpha(n)\lambda\Sigma_2 \bigg]+(1+\lambda) \bigg[ [\Sigma'_2+(n-2)\theta]\mu_{11}\Sigma_1\alpha(n-1)\lambda \Sigma_2 \bigg]\\
    &+\lambda \Sigma_2 \Sigma'_2 \mu_{22}f(1,n-1)-\lambda(n-1)(n-2)\theta^2\mu_{22}f(1,n-1)\\
    =& (1+\lambda)\tau'_1(n)\\
    &-(1+\lambda)\lambda\Sigma_2\bigg[\Sigma_1 \alpha(n)\mu_{12}+\Sigma_2 f(1,n)-\mu_{11}\Sigma_2 f(1,n-1)-\mu_{11}\Sigma_1 \alpha(n-1)\mu_{12}+[\Sigma_2'+(n-1)\theta]\Sigma_1 \alpha(n)\\
    &-[\Sigma'_2+(n-2)\theta]\mu_{11}\Sigma_1\alpha(n-1)-\Sigma_2\mu_{22}f(1,n-1)\bigg]-\lambda(n-1)(n-2)\theta^2\mu_{22}f(1,n-1)\\
    =& (1+\lambda)\tau'_1(n)\\
    &-(1+\lambda)\lambda\Sigma_2\bigg[\Sigma_1 \alpha(n)\mu_{12}+\Sigma_2 f(1,n)-\mu_{11}\Sigma_2 f(1,n-1)-\mu_{11}\Sigma_1 \alpha(n-1)\mu_{12}+\Sigma_2\Sigma_1 \alpha(n)\\
    &+[\lambda \Sigma_2+(n-1)\theta]\Sigma_1 \alpha(n)-\Sigma_2\Sigma_1 \mu_{11}\alpha(n-1)\\
    &-[\lambda \Sigma_2+(n-2)\theta]\mu_{11}\Sigma_1\alpha(n-1)-\Sigma_2\mu_{22}f(1,n-1)\bigg]-\lambda(n-1)(n-2)\theta^2\mu_{22}f(1,n-1)\\
    =& (1+\lambda)\tau'_1(n)\\
    &-(1+\lambda)\lambda\Sigma_2\bigg[\Sigma_1 \alpha(n)\mu_{12}+\Sigma_2 f(1,n)-\mu_{11}\Sigma_2 f(1,n-1)-\mu_{11}\Sigma_1 \alpha(n-1)\mu_{12}+\Sigma_2\Sigma_1 \alpha(n)\\
    &-\Sigma_2\Sigma_1 \mu_{11}\alpha(n-1)-\Sigma_2\mu_{22}f(1,n-1)\bigg]-(1+\lambda)\lambda \Sigma_2\bigg[[\lambda \Sigma_2+(n-1)\theta]\Sigma_1 \alpha(n)\\
    &-[\lambda \Sigma_2+(n-2)\theta]\mu_{11}\Sigma_1\alpha(n-1)\bigg]-\lambda(n-1)(n-2)\theta^2\mu_{22}f(1,n-1)\\
      =& (1+\lambda)\tau'_1(n)\\
    &-(1+\lambda)\lambda\Sigma_2\bigg[\Sigma_1 \alpha(n)\mu_{12}+\Sigma_2(\mu_{22}+(n-2)\theta) f(1,n-1)-\mu_{11}\Sigma_2 f(1,n-1)-\mu_{11}\Sigma_1 \alpha(n-1)\mu_{12}\\
    &+\Sigma_2\Sigma_1 \alpha(n)-\Sigma_2\Sigma_1 \mu_{11}\alpha(n-1)-\Sigma_2\mu_{22}f(1,n-1)\bigg]-(1+\lambda)\lambda \Sigma_2\bigg[[\lambda \Sigma_2+(n-2)\theta]\Sigma_1 \alpha(n)\\
    &-[\lambda \Sigma_2+(n-2)\theta]\mu_{11}\Sigma_1\alpha(n-1)\bigg]-(1+\lambda)\lambda \Sigma_2\theta \Sigma_1\alpha(n)-\lambda(n-1)(n-2)\theta^2\mu_{22}f(1,n-1)\\
    =& (1+\lambda)\tau'_1(n)\\
    &-(1+\lambda)\lambda\Sigma_2\bigg[\Sigma_1 \alpha(n)\mu_{12}-\mu_{11}\Sigma_2 f(1,n-1)-\mu_{11}\Sigma_1 \alpha(n-1)\mu_{12}+\Sigma_2\Sigma_1 \alpha(n)-\Sigma_2\Sigma_1 \mu_{11}\alpha(n-1)\bigg]\\
    &-(1+\lambda)\lambda\Sigma_2^2(n-2)\theta f(1,n-1)-(1+\lambda)\lambda \Sigma_2[\lambda \Sigma_2+(n-2)\theta]\Sigma_1 (\alpha(n)-\mu_{11}\alpha(n-1))\\
    &-(1+\lambda)\lambda \Sigma_2\theta \Sigma_1\alpha(n)-\lambda(n-1)(n-2)\theta^2\mu_{22}f(1,n-1).
\end{align*}
\subsection{Table for $\bar N$ Frequency Under Task Dependent Synergy}\label{SMTable}

\begin{center}
\begin{tabular}{ |c|c| } 
\hline
\ &\ \\
$\bar N$ value & Frequency\\
\hline \hline
1 & 7702\\
\hline
2 &326\\
\hline
3 & 339\\
\hline
4 & 296\\
\hline
5 & 242\\
\hline
6&176\\
\hline
7&139\\
\hline
8&101\\
\hline
9&73\\
\hline
10&70 \\
\hline
11&51\\
\hline
12&36\\
\hline
13&41\\
\hline
14&41\\
\hline
15&22\\
\hline
16&24\\
\hline
17&21\\
\hline
18&21\\
\hline
19&9\\
\hline
20&15\\
\hline
21&18\\
\hline
22&20\\
\hline
23&10\\
\hline
24&11\\
\hline
25&11\\
\hline
26&4\\
\hline
27&7\\
\hline
28&4\\
\hline
29&5\\
\hline
30&6\\
\hline
\end{tabular}
\quad
\begin{tabular}{ |c|c|}
\hline
\ &\ \\
$\bar N$ Value & Frequency\\
\hline \hline
31&6\\
\hline
32&4\\
\hline
33&8\\
\hline
34&8\\
\hline
35&7\\
\hline
36&11\\
\hline
37&4\\
\hline
38&5\\
\hline
39&1\\
\hline
40&3\\
\hline
41&5\\
\hline
42&3\\
\hline
43&3\\
\hline
44&4\\
\hline
45&3\\
\hline
46&2\\
\hline
47&2\\
\hline
48&2\\
\hline
49&4\\
\hline
50&4\\
\hline
51&3\\
\hline
52&2\\
\hline
53&4\\
\hline
54&5\\
\hline
55&1\\
\hline
56&2\\
\hline
57&2\\
\hline
58&1\\
\hline
59&2\\
\hline
60&1\\
\hline
\end{tabular}
\quad
\begin{tabular}{ |c|c|}
\hline
\ &\ \\
$\bar N$ Value & Frequency\\
\hline \hline
61&2\\
\hline
62&2\\
\hline
63&3\\
\hline
65&1\\
\hline
66&1\\
\hline
67&3\\
\hline
70&1\\
\hline
71&1\\
\hline
72&3\\
\hline
73&1\\
\hline
74&1\\
\hline
75&2\\
\hline
76&4\\
\hline
78&1\\
\hline
80&2\\
\hline
83&2\\
\hline
84&2\\
\hline
85&1\\
\hline
86&1\\
\hline
87&2\\
\hline
88&2\\
\hline
89&1\\
\hline
90&1\\
\hline
92&1\\
\hline
93&1\\
\hline
94&1\\
\hline
96&1\\
\hline
98&1\\
\hline
99&1\\
\hline
\end{tabular}
\captionof{table}{Frequency of $\bar N$ values}
\label{Nfrequency}
\end{center}
\end{document}